\documentclass[12pt]{amsart}
\font\emailfont=cmtt10

\usepackage{amsmath,amsthm,amsfonts,amscd,flafter,epsf}

\newcommand\commentable[1]{#1}
\newcommand\Rk{\mathrm{rk}}

\newcommand{\rk}{\mathrm{rk}}
\newcommand{\HF}{HF}

\newtheorem{theorem}{Theorem}[section]
\newtheorem{prop}[theorem]{Proposition}
\newtheorem{cor}[theorem]{Corollary}

\newtheorem{lemma}[theorem]{Lemma}

\newtheorem{defn}[theorem]{Definition}

\def\endproof{\relax\ifmmode\expandafter\endproofmath\else
  \unskip\nobreak\hfil\penalty50\hskip.75em\hbox{}\nobreak\hfil\bull
  {\parfillskip=0pt \finalhyphendemerits=0 \bigbreak}\fi}
\def\endproofmath$${\eqno\bull$$\bigbreak}
\def\bull{\vbox{\hrule\hbox{\vrule\kern3pt\vbox{\kern6pt}\kern3pt\vrule}\hrule}}

\newcommand{\Q}{\mathbb{Q}}
\newcommand{\R}{\mathbb{R}}

\newcommand{\Z}{\mathbb{Z}}

\newcommand{\Zmod}[1]{\Z/{#1}\Z}

\newcommand{\Ker}{\mathrm{Ker}}

\newcommand{\grad}{\vec\nabla}

\newcommand{\cm}{\cdot}

\newcommand{\Nbd}[1]{{\mathrm{nd}}(#1)}
\newcommand{\nbd}[1]{\Nbd{#1}}

\newcommand{\ModSWfour}{\mathcal{M}}
\newcommand{\ModFlow}{\ModSWfour}

\newcommand{\SpinC}{{\mathrm{Spin}}^c}

\newcommand{\goesto}{\mapsto}

\newcommand\Wedge{\Lambda}

\newcommand\abuts\Rightarrow
\newcommand\Sym{\mathrm{Sym}}

\newcommand\spinccan{\ell}

\newcommand\RelSpinC{\underline{\SpinC}}
\newcommand\relspinc{\underline{\spinc}}

\newcommand\Filt{\mathcal F}

\newcommand\x{\mathbf x}
\newcommand\w{\mathbf w}

\newcommand\y{\mathbf y}

\newcommand\ModSphere{\ModFlow\left({\mathbb S}\longrightarrow 
\Sym^{g-1}(\Sigma_{1})\times \Sym^2(\Sigma_{2})\right)}
\newcommand\ModSpheres\ModSphere
\newcommand\CF{CF}

\newcommand\CFa{\widehat{CF}}
\newcommand\CFp{\CFb}
\newcommand\CFm{\CF^-}

\newcommand\HFleq{\HF^{\leq 0}}

\newcommand\HFp{\HFb}

\newcommand\CFinf{CF^\infty}

\newcommand\CFb{CF^+}
\newcommand\HFa{\widehat{HF}}
\newcommand\HFb{HF^+}

\newcommand\Mas{\mu}
\newcommand\UnparModSp{\widehat \ModSp}
\newcommand\UnparModFlow\UnparModSp
\newcommand\Mod\ModSp

\newcommand{\cald}{{\mathcal D}}

\newcommand\PD{\mathrm{PD}}

\newcommand{\spinc}{\mathfrak s}

\newcommand{\spinct}{\mathfrak t}

\newcommand\ModMaps{\mathcal M}
\newcommand\ModSp\ModMaps

\newcommand\Ta{{\mathbb T}_{\alpha}}
\newcommand\Tb{{\mathbb T}_{\beta}}
\newcommand\Tc{{\mathbb T}_{\gamma}}
\newcommand\Td{{\mathbb T}_{\delta}}

\newcommand\alphas{\mbox{\boldmath$\alpha$}}

\newcommand\betas{\mbox{\boldmath$\beta$}}
\newcommand\gammas{\mbox{\boldmath$\gamma$}}
\newcommand\deltas{\mbox{\boldmath$\delta$}}

\newcommand\PerDom{\mathcal P}

\newcommand\uCFp{\underline{\CFp}}

\newcommand\Fp[1]{F^{+}_{#1}}

\newcommand\Field{\mathbb F}

\newcommand\Dual{\mathcal D}
\newcommand\Duality\Dual

\newcommand\TaPr{\Ta'}
\newcommand\TbPr{\Tb'}
\newcommand\TcPr{\Tc'}

\newcommand\alphaprs{\alphas'}
\newcommand\betaprs{\betas'}
\newcommand\gammaprs{\gammas'}
\newcommand\Spider{\sigma}
\newcommand\EulerMeasure{\widehat\chi}

\newcommand\InjMod[1]{{\mathcal T}^+_{#1}}

\newcommand\MCone{M}

\newcommand\CFK{CFK}
\newcommand\HFK{HFK}

\newcommand\CFKinf{\CFK^{\infty}}

\newcommand\HFKa{\widehat\HFK}

\newcommand\BK{K}
\newcommand\Fill[1]{G_{#1}}
\newcommand\FillW[1]{E_{#1}}
\newcommand\Width{b}

\newcommand\FormHFp{\mathbb{HF}^+}
\newcommand\Reflect[1]{-{#1}}

\newcommand\EquivClass{{\mathfrak T}(\xi)}
\newcommand\OrK{\overline K}

\newcommand\Xa{\widehat{\mathbb X}}
\newcommand\Xp{{\mathbb X}^+}
\newcommand\Xd{{\mathbb X}^{\delta}}

\newcommand\Aa{\widehat{A}}
\newcommand\Ba{\widehat{B}}

\newcommand\verta{\widehat v}
\newcommand\hora{\widehat h}
\newcommand\vertp{v^+}
\newcommand\horp{h^+}
\newcommand\vertd{v^{\delta}}
\newcommand\hord{h^{\delta}}

\newcommand\fp{f^+}
\newcommand\Hp[1]{H^+_{#1}}

\newcommand\spincrel{\underline\spinc}

\newcommand\Ap{{A}^+}
\newcommand\Bp{{B}^+}

\newcommand\Ad{{A}^{\delta}}
\newcommand\Bd{{B}^{\delta}}

\newcommand\fd{f^{\delta}}
\newcommand\Hd{H^{\delta}}

\newcommand\CFd{CF^{\delta}}
\newcommand\HFd{HF^{\delta}}
\newcommand\Dd{D^{\delta}}

\newcommand\BigAp{{\mathbb A}^+}
\newcommand\BigBp{{\mathbb B}^+}
\newcommand\Dp{D^+}

\newcommand\CapSurf{\widehat F}

\newcommand\Ainf{A^{\infty}}
\newcommand\Binf{B^{\infty}}
\newcommand\Dinf{D^{\infty}}

\newcommand\BigAd{{\mathbb A}^{\delta}}
\newcommand\BigBd{{\mathbb B}^{\delta}}
\newcommand\spincx{\mathfrak x}
\newcommand\spincy{\mathfrak y}

\newcommand\BigAa{\widehat{\mathbb A}}
\newcommand\BigBa{\widehat{\mathbb B}}
\newcommand\Da{{\widehat D}}

\commentable{

\title[{Knot Floer homology and rational surgeries}] 
{Knot Floer homology and rational surgeries}

\author[Peter Ozsv{\'a}th]{Peter Ozsv\'ath}
\address{Department of
Mathematics, UC/Berkeley, Berkeley, CA 94720 \newline
\indent{\emailfont{petero@math.berkeley.edu}}}
\thanks{PSO was supported by NSF grant number DMS 0234311}

\author[Zolt{\'a}n Szab{\'o}]{Zolt{\'a}n Szab{\'o}} 
\address{Department of
Mathematics, Princeton University, New Jersey 08540 \newline
\indent{\emailfont{szabo@math.princeton.edu}}}}
\thanks{ZSz was supported by NSF grant number DMS 0107792}

\begin{document}

\begin{abstract}  
Let $K$ be a rationally null-homologous knot in a three-manifold $Y$.
We construct a version of knot Floer homology in this context, including
a description of the Floer homology of a three-manifold obtained
as Morse surgery on the knot $K$. As an application, we express the
Heegaard Floer homology of rational surgeries on $Y$ along a
null-homologous knot $K$ in terms of the filtered homotopy type of the
knot invariant for $K$. This has applications to Dehn surgery
problems for knots in $S^3$.  In a different direction,
we use the techniques developed here to calculate the
Heegaard Floer homology of an arbitrary Seifert fibered
three-manifold.
\end{abstract} 

\maketitle

\section{Introduction}

Heegaard Floer homology~\cite{HolDisk} is an invariant for closed, oriented
three-manifolds $Y$, taking the form of a
collection of homology groups which are functorial under cobordisms.
In~\cite{Knots} and~\cite{RasmussenThesis}, this invariant is extended
to an invariant for null-homologous knots $K$. 
(Here, we say that a knot is null-homologous if its induced homology
class in $H_1(Y;\Z)$ is trivial. If the induced homology class of a
knot $K\subset Y$ in $H_1(Y;\Q)$ is trivial, we call the knot {\em
rationally null-homologous}.)  The knot invariant takes the form of a
$\Z\oplus\Z$-filtration of the chain complex whose homology calculates
the Heegaard Floer complex for $Y$. It is the filtered chain homotopy
type of this filtered complex which depends on the particular knot
$K$.

The knot filtration gives rise to collection of chain complexes 
$\{\Ap_s(K)\}_{s\in\Z}$ and chain maps $\{\vertp_s\colon
\Ap_s(K)\longrightarrow \Bp\}_{s\in\Z}$ and $\{\horp_s\colon
\Ap_s(K)\longrightarrow \Bp\}_{s\in\Z}$, where here $\Bp=\CFp(Y)$ is a
chain complex whose homology is the Heegaard Floer homology $\HFp(Y)$.
Indeed, the homology groups of the chain complex $\Ap_s(K)$ represents the
homology $\HFp$ of sufficiently large integer surgeries on $Y$ along
$K$, in a sense which can be made precise (c.f.
Theorem~\ref{Knots:thm:LargePosSurgeries} of~\cite{Knots} and
also~\cite{RasmussenThesis}). (These complexes are defined in a more general
setting in Section~\ref{sec:Construction}.)

Suppose that $K$ is a null-homologous knot in a three-manifold $Y$.
Given a rational number $r$, let
$Y_r(K)$ denote the three-manifold obtained by Dehn filling
$Y-\nbd{K}$ with a solid torus with slope $r$ (with respect to the
canonical Seifert framing of $K$). In the case where $r$ is an
integer, the Heegaard Floer homology of $Y_r(K)$ can described in
terms of the above-mentioned data coming from the knot filtration,
according to~\cite{IntSurg}.

The primary aim of this article is to generalize this construction to
the case of Morse surgery on a knot $K\subset Y$ which is only
rationally null-homologous. (By Morse surgery, we mean here Dehn
surgery on a knot which can be realized as the boundary of 
a single two-handle addition to $[0,1]\times Y$; in
the case where $K$ is null-homologous, this corresponds to Dehn
surgery with an integral slope). This construction has new
consequences even in the case of null-homologous knots in a
three-manifold: since the result of Dehn surgery on a null-homologous
knot $K\subset Y$ can be viewed as Morse surgery on a knot in the
connected sum of $Y$ with a lens space, we obtain a description of the
Heegaard Floer homology of $Y_r(K)$ in terms of the original knot
Floer homology of $K$.

Rather than introducing the generalization of the knot package to
knots which are only rationally null-homologous in this introduction,
which will require some additional material in its statement (c.f.
Sections~\ref{sec:Construction}, \ref{sec:Kunneth},
and~\ref{sec:Surgery} below), we focus now in the description of
the Floer homology of $Y_r(K)$ when $r$ is a rational number, and
$K\subset Y$ is null-homologous.

As a preliminary point, recall that the Heegaard Floer homology of $Y$
admits a direct sum splitting indexed by the set of $\SpinC$
structures over $Y$, which in turn is an affine space for $H^2(Y;\Z)$.
In particular, if $K\subset Y$ is a knot in an integral homology
sphere, then there is a splitting
$$\HFp(Y_{p/q}(K))\cong \bigoplus_{i\in\Zmod{p}} \HFp(Y_{p/q}(K),i).$$
Fix an integer $i$, and  consider the chain complexes
\begin{eqnarray*}
\BigAp_i=\bigoplus_{s\in\Z}(s,\Ap_{\lfloor \frac{i+ps}{q}\rfloor}(K))
&{\text{and}}&
\BigBp_i=\bigoplus_{s\in\Z}(s,\Bp),
\end{eqnarray*}
where here $\lfloor x\rfloor$ denotes the greatest integer smaller
than or equal to $x$, and all $\Bp_{\lfloor\frac{i+ps}{q}\rfloor}=\CFp(Y)$ .
We view the above chain homomorphisms $\vertp$ and $\horp$ as maps
\begin{eqnarray*}
\vertp\colon (s,\Ap_{\lfloor \frac{i+ps}{q}\rfloor}(K))\longrightarrow
(s,\Bp)
&{\text{and}}&
\horp\colon 
(s,\Ap_{\lfloor \frac{i+ps}{q}\rfloor}(K))
\longrightarrow
(s-1,\Bp).
\end{eqnarray*}
Adding these up, we obtain a chain map
$$\Dp_{i,p/q}\colon \BigAp_i \longrightarrow \BigBp_i;$$
i.e.
$$\Dp_{i,p/q} \{(s,a_s)\}_{s\in\Z}
=\{(s,b_s)\}_{s\in\Z},$$
where here
$$b_{s}
=\vertp_{\lfloor \frac{i+ps}{q}\rfloor}(a_s)
+\horp_{\lfloor \frac{i+p(s-1)}{q}\rfloor}(a_{s-1}).$$

Let $\Xp_{i,p/q}$ denote the mapping cone of $\Dp_{i,p/q}$.  Note that
$\Xp_{i,p/q}$ depends on $i$ only through its congruence class 
modulo $p$.  Note also that $\Ap_s$ and $\Bp_s$ are relatively
$\Z$-graded, and the homomorphisms $\vertp_s$ and $\horp_s$ respect
this relative grading. The mapping cone $\Xp_i$ can be endowed with a
relative grading, with the convention that $\Dp_{i,p/q}$ drops the
grading by one.

This mapping cone, whose ingredients are extracted from the knot
filtration, captures Heegaard Floer homology of $p/q$ surgeries on $Y$
along $K$, according to the following:

\begin{theorem}
\label{thm:RationalSurgeries}
Let $K\subset Y$ be a null-homologous knot, and let $p$, $q$ be a pair
of relatively prime integers.  Then, for each $i\in\Zmod{p}$,
there is a relatively graded
isomorphism of groups $$H_*(\Xp_{i,p/q})\cong \HFp(Y_{p/q}(K),i).$$
\end{theorem}

The proof of Theorem~\ref{thm:RationalSurgeries} is based on
generalization of the the knot filtration of~\cite{Knots}
and~\cite{RasmussenThesis} to the case of rationally null-homologous
knots, together with a generalization of the integer surgeries
description from~\cite{IntSurg} (where in fact
Theorem~\ref{thm:RationalSurgeries} is proved in the case where
$q=1$). Note that in the more general case, the knot filtration is
naturally a filtration by relative $\SpinC$ structures on the knot
complement (rather than integers).  Rational surgeries on $K\subset Y$
can be realized as Morse surgeries on a knot in the connected sum
of $Y$ with a lens space.  The resulting knot is gotten by forming the
connected sum of $K$ with a model (homologically non-trivial) knot in
the lens space.  Theorem~\ref{thm:RationalSurgeries} is then realized
as a combination of a straighforward calculation involving this model knot,
combined with a K\"unneth principle for connected sums, followed by
the general Morse surgeries description.

We turn now to various applications of Theorem~\ref{thm:RationalSurgeries}.

\subsection{Applications to knots with $L$-space surgeries}

In Section~\ref{sec:LSpaceKnots}, we give applications of
Theorem~\ref{thm:RationalSurgeries} to knots in $S^3$ which 
admit $L$-space surgeries.

Recall that an $L$-space is a rational homology three-sphere $Y$ whose
Floer homology $\HFp$ in each $\SpinC$ structure is isomorphic (as a
relatively-graded $\Z[U]$-module) to $\HFp(S^3)$. This is equivalent
to the condition that $\HFa(Y,\spinc)\cong \Z$ for each
$\spinc\in\SpinC(Y)$.  Recall also~\cite{HolDiskFour}, \cite{AbsGraded}
that if $Y$ is a rational homology three-sphere, then $\HFa(Y,\spinc)$
is a $\Q$-graded group. Thus, the Heegaard Floer homology of an
$L$-space is determined by the ``correction term''
function $$d\colon
\SpinC(Y)\longrightarrow \Q$$ which associates to each
$\spinc\in\SpinC(Y)$ the degree in which $\HFa(Y,\spinc)$ is
supported, compare also~\cite{Froyshov}. 

The set of $L$-spaces includes all lens spaces and indeed all
three-manifolds with elliptic geometry, c.f.~\cite{NoteLens}; for more
examples, see also~\cite{Nemethi}. Another interesting family is given
by the branched double-covers of alternating knots in $S^3$,
c.f.~\cite{BrDCov}.

Let $K\subset S^3$ be a knot in the three-sphere. Write its
symmetrized Alexander polynomial as
$$\Delta_K(T)=a_0 + \sum_{i>0}a_i(T^i+T^{-i}),$$
and let
\begin{equation}
\label{eq:DefTorsion}
t_i(K)=\sum_{j=1}^\infty j a_{|i|+j}.
\end{equation}

Note that for any knot $C$ in $S^3$, there is a canonical affine map
$\Zmod{p}\cong \SpinC(S^3_{p/q}(C))$.

\begin{theorem}
\label{thm:RatSurgeryLSpace}
Let $K\subset S^3$ be a knot which admits an $L$-space surgery, for
some $r=\frac{p}{q}\in \Q$ with $r\geq 0$.  Then, 
for all integers $i$ with  $|i|\leq
\frac{p}{2q}$
we have that
\begin{equation}
\label{eq:DDifferences}
d(S^3_{p/q}(K),i)-d(S^3_{p/q}(O),i) =
-2t_{\lfloor\frac{|i|}{q}\rfloor}(K),
\end{equation}
 while for all $|i|>
\frac{p}{2q}$, we have that $t_i(K)=0$.  
\end{theorem}

In the case where $q=1$, a version of the above theorem is established
in~\cite{AbsGraded} (c.f. Theorem~\ref{AbsGraded:thm:PSurgeryLens}
of~\cite{AbsGraded}). Theorem~\ref{thm:RatSurgeryLSpace} (in the case
where $q=2$) also gives the symmetry used in~\cite{UnknotOne} to find
an obstruction to a knot having unknotting number equal to one.

The following is a quick consequence of this result, together with the
fact that knot Floer homology distinguishes the unknot
(c.f.~\cite{GenusBounds}) (though alternative proofs could be given
which model the proof in~\cite{KMOSz} more closely):

\begin{cor}
\label{cor:GordonConjecture}
If $S^3_{p/q}(K)\cong S^3_{p/q}(O)$ as oriented three-manifolds, then
$K=O$.
\end{cor}

The above result, which was conjectured by Gordon
in~\cite{GordonConjecture}, was first established using Seiberg-Witten
monopoles in~\cite{KMOSz}. Thanks to a theorem of Eliashberg and
Etnyre (c.f.~\cite{Eliashberg} and~\cite{Etnyre}), it is now possible
to prove results of this type purely in the context of Heegaard Floer
homology, see also~\cite{GenusBounds}.

\subsection{On cosmetic surgeries}
\label{subsec:CosmeticSurgeries}

Let $Y$ be a closed, oriented three-manifold, and $K\subset Y$ be a framed
knot. Given a rational number $r$, let $Y_r(K)$ denote the
three-manifold obtained by Dehn surgery along $K$ with slope $r$ (with
respect to the initial framing). If there are two distinct rational
numbers $r$ and $s$ with the property that $Y_r(K)$ and $Y_s(K)$ are
homeomorphic (but the homeomorphism is not
required to preserve the orientation inherited from $Y$),
then the surgeries are called {\em cosmetic}.
A pair of surgeries on $K$ with $r\neq s$ is called
{\em truly cosmetic} if $Y_r(K)\cong Y_s(K)$
as oriented manifolds.

Amphicheiral knots have cosmetic surgeries; specifically, if $K$ is an
amphicheiral knot, then $S^3_r(K)\cong -S^3_{-r}(K)$. The unknot $O$ admits
infinitely many truly cosmetic surgeries: $S^3_{p/q}(O)=S^3_{p/{p+q}}(O)$.
Lackenby~\cite{Lackenby} has shown that under general conditions on a
knot $K\subset Y$, there are at most finitely many cosmetic surgeries,
see also~\cite{BleilerHodgsonWeeks}.
It is conjectured~\cite{BleilerHodgsonWeeks}
that if $Y_r(K)\cong Y_s(K)$, then $Y-\nbd{K}$ admits
an automorphism which carries the slope $r$ to the slope $s$.

The present state of Heegaard Floer homology -- and specifically the
surgery formulas given here -- work best for excluding cosmetic
surgeries on knots in $S^3$.  For example, we have the following
result:

\begin{theorem}
\label{intro:GenusOne}
If $K$ is a knot with Seifert genus equal to one,
and $S^3_r(K)\cong S^3_s(K)$ with $r\neq s$, then
$S^3_r(K)$ is
an $L$-space.
\end{theorem}

The conclusion of the above theorem places severe restrictions on
$K$.   In particular, according to results of~\cite{NoteLens}, it
follows that $K$ must have the same knot Floer homology (and in
particular the same Alexander poylnomial) as the trefoil knot $T$, and
thus according to Theorem~\ref{thm:RationalSurgeries}, $S^3_r(K)$ and
$S^3_r(T)$ have the same (graded) Floer homology groups. For
particular integers $p$, the existence
of a truly cosmetic surgery on such a knot $K$ with specified
numerator $p$ can be ruled out by an explicit, finite search.

\begin{theorem}
\label{intro:OppositeSigns}
Let $K\subset S^3$ and suppose that
$S^3_r(K)\cong \pm S^3_s(K)$, then either $S^3_r(K)$ is an $L$-space
or $r$ and $s$ have opposite signs.
\end{theorem}

For the above theorem, both possible conclusions can hold.  The
simplest example is the unknot
which admits cosmetic surgeries with
positive slopes. A more interesting example of cosmetic
surgeries with positive slopes is provided by 
the trefoil knot $K$,
which has the property that $S^3_9(K)\cong -S^3_{9/2}(K)$,
c.f.~\cite{Mathieu}. Examples where $r$ and $s$ have opposite signs
are given by amphicheiral knots.

Our methods can be refined to exclude cosmetic surgeries for certain
numerators $p$. We study here the case where $p=3$.

\begin{theorem}
\label{intro:PEqualsThree}
Suppose that $K\subset S^3$ is a knot
with the property that 
$S^3_{p/q}(K)\cong S^3_{p/q'}(K)$ as oriented manifolds.
In the case where $p=3$, we can conclude that $q=q'$.
\end{theorem}

\subsection{Heegaard Floer homology of Seifert fibered spaces}

We give some other applications of the general surgeries description
along a rationally null-homologous knot. In
Section~\ref{sec:Seiferts}, we use it to describe the Heegaard Floer
homology of any Seifert fibered space whose first Betti number is even.

\subsection{Organization}

This paper is organized as follows. In Section~\ref{sec:TopPre}, we
review some of the topological preliminaries required by the knot
filtrations, including the notion of {\em relative $\SpinC$
  structures} for three-manifolds with torus boundary. In
Section~\ref{sec:Construction}, we give the construction of the knot
filtration. In the next two sections, we turn to some properties of
the knot Floer homology which are rather straightforward adaptations
of the corresponding results for null-homologous knots (\cite{Knots},
\cite{RasmussenThesis}): the relationship between knot Floer homology
and ``large'' surgeries on a rationally null-homologous knot
(Section~\ref{sec:Large}), and the K\"unneth principle for connected
sums of knots (Section~\ref{sec:Kunneth}). The first result is an
ingredient in the Morse surgery formula from
Section~\ref{sec:Surgery}. In Section~\ref{sec:RatSurg}, we show how
the K\"unneth principle, together with the Morse surgery formula, give
the rational surgery formula described in this introduction. In
Section~\ref{sec:LSpaceKnots}, we turn to knots which admit $L$-space
surgeries.  In Section~\ref{sec:CosmeticSurgeries}, we give the
applications of the rational surgery formula to the problem of
cosmetic surgeries on a knot in $S^3$.  In Section~\ref{sec:Seiferts},
we turn to the Heegaard Floer homology groups of Seifert fibered
spaces with even first Betti number.

\subsection{Acknowledgements} The
authors would like to thank Yi Ni, Andr{\'a}s Stipsicz, and Jacob
Rasmussen for interesting conversations during the course of this
work.  We are also very grateful to Walter Neumann for sharing with us
his expertise in hyperbolic geometry.

\section{Preliminaries}
\label{sec:TopPre}

We recall some of the background material and notation which will be
required for the construction of the knot filtration.  The bulk of
this material is about the constructions relating (doubly-pointed)
Heegaard diagrams for knots, and relative $\SpinC$ structures for
three-manifolds with boundary. We include also a few key properties of
Heegaard Floer homology which will be used repeatedly throughout,
see also~\cite{HolDisk}.

\subsection{Surgeries}
\label{subsec:Surgeries}

Let $K\subset Y$ be a knot. the boundary of a tubular neighborhood of
the knot $K$ is a torus $T$ equipped with a canonical choice of
(isotopy class of) embedded curve $\mu$, a {\em meridian} for $K$.  A
{\em longitude} for $K$ is any embedded curve $\lambda$ in $T$ which
meets a meridian transversally in a single point.

Given a homologically non-trivial, embedded curve $\gamma$ in $T$, we
can form the three-manifold $Y_\gamma(K)$ which is gotten by attaching
a solid torus to $Y-\nbd{K}$ so that $\gamma$ bounds a disk in the
attached solid torus. We say that $Y_\gamma(K)$ is obtained from $Y$
by {\em Dehn surgery along $K$}. In the special case where $\gamma$ is
a longitude for $K$, there is a canonical two-handle cobordism from
$Y$ to $Y_\gamma(K)$, and we say that $Y_\gamma(K)$ is obtained from
$Y$ by {\em Morse surgery on $K$}. A choice of longitude for $K$ is
also called a {\em framing} of $K$.

Let $K\subset Y$ be a rationally null-homologous knot in a closed,
oriented three-manifold, equipped with a framing $\lambda$. 
Since $K\in Y$ has finite order, there is a 
pair of integers
$n$ and $d$ with minimal absolute value which satisfy the property that
\begin{equation}
\label{eq:OrderD}
d \cm \lambda = n\cm \mu \in H_1(Y-K;\Z),
\end{equation}
In particular, it follows that the induced homology class
$[K_\lambda]\in H_1(Y;\Z)$ has order $|d|$. Here, and in the future,
$K_\lambda$ denotes a copy of $K$ displaced into $Y-\nbd{K}$ using the
framing $\lambda$.

\subsection{Relative $\SpinC$ structures for three-manifolds with boundary}
\label{subsec:RelSpinC}

Following Turaev~\cite{Turaev}, we say that two vector fields $v_1$ and
$v_2$ on a closed, oriented three-manifold are {\em homologous} if
they are homotopic in the complement of a ball in $Y$. The set of
homology classes of vector fields can be given the structure of an
affine space for $H^2(Y;\Z)$, and indeed, it can be identified with
the space of $\SpinC$ structures over $Y$, $\SpinC(Y)$.

This construction can be readily generalized to the case of a
three-manifold with torus boundary, c.f. Chapter I.4 of~\cite{Turaev}.
Specifically, on an oriented three-manifold $M$ with torus boundary,
two vector fields $v_1$ and $v_2$ on $M$ which vanish nowhere and
point outwards at $\partial M$ are said to be {\em homologous} if they
are homotopic in the complement of a ball in the interior of $M$. The
set of homology classes of nowhere vanishing vector fields can be
naturally given the structure of an affine space for the relative
cohomology $H^2(M,\partial M;\Z)$. The homology classes of nowhere
vanishing vector fields are called {\em relative $\SpinC$ structures
on $M$}, and are denoted $\RelSpinC(M,\partial M)$.

Let $K\subset Y$ be a knot in a closed, oriented three-manifold, we
can construct the three-manifold with torus boundary $M=Y-\nbd{K}$.
In this case, we denote the relative $\SpinC$ structures over $M$ by
$\RelSpinC(Y,K)$.

Orienting the core of the solid torus, we obtain a canonical isotopy
class of nowhere vanishing vector field which points inward at the
boundary, and which has a closed orbit which is the core of the solid
torus (with its given orientation). More explicitly, the isotopy class
of this vector field on $D^2\times S^1$ is characterized by the property
that in the interior of the solid torus, $w$ is everywhere transverse to the
tangent planes to $D^2$. 

Giving $K$ an orientation (and denoting this
oriented knot by $\OrK$), we glue in this vector field $w$ to obtain a
natural map $$\Fill{Y,\OrK}\colon \RelSpinC(Y,K) \longrightarrow
\SpinC(Y)$$
which is equivariant with respect to the action by $H^2(Y,K;\Z)$; i.e.
letting 
$$\iota \colon H^2(Y,K;\Z)\longrightarrow H^2(Y;\Z)$$ be the
natural map, we have for each $k\in H^2(Y,K;\Z)$,
$$\Fill{Y,\OrK}(\xi+k)=\Fill{Y,\OrK}(\xi)+\iota(k).$$ As the notation
suggests, $\Fill{Y,\OrK}$ depends on the orientation for $K$; indeed,
for any $\xi\in\SpinC(Y,K)$,
\begin{equation}
\label{eq:DependenceOnOrientation}
\Fill{Y,\OrK}(\xi)-\Fill{Y,-\OrK}(\xi)=\PD[\OrK].
\end{equation}

Note that if $Y$ is a three-manifold, and $K\subset Y$ is a knot with
meridian $\mu$, then the fibers of the map
$\Fill{Y,\OrK}$ are the orbits of $\RelSpinC(Y,K)$ under the action
by $\Z\cdot \PD[\mu]$, where we think of $\PD[\mu]\in H^2(Y,K;\Z)$;
i.e. $\Fill{Y,\OrK}$ realizes an identification
$$\SpinC(Y)\cong \frac{\RelSpinC(Y,K)}{\Z\cdot \PD[\mu]}.$$

\subsection{Doubly-pointed Heegaard diagrams}
\label{subsec:DoublyPointedDiagrams}

A {\em doubly-pointed Heegaard diagram} is a collection of data
$(\Sigma,\alphas,\betas,w,z)$ where $\Sigma$ is an oriented surface of
genus $g$, $\alphas=\{\alpha_1,...,\alpha_g\}$ is a $g$-tuple of
homologically linearly independent, pairwise disjoint, embedded curves
in $\Sigma$ (a $g$-tuple of attaching circles), and
$\betas=\{\beta_1,...,\beta_g\}$ is another $g$-tuple of attaching
circles, and $w$ and $z$ are two distinct pair of points in
$\Sigma-\alpha_1-...-\alpha_g-\beta_1-...-\beta_g$.  

A doubly-pointed Heegaard diagram gives rise to an oriented
three-manifold $Y$ together with an oriented knot $\OrK\subset Y$. The
orientation on $\Sigma$ induces an orientation on $U_\alpha$
(so that $\partial U_\alpha=\Sigma$), which can then be uniquely extended
to an orientation over $Y$. The knot $K$ is obtained
as a union of two arcs, $\eta_\alpha$ and $\eta_\beta$.  The arc
$\eta_\alpha$ is gotten by connecting $w$ and $z$ in
$\Sigma-\alpha_1-...-\alpha_g$, and orienting it as a path from $w$ to
$z$. This arc is then pushed into $U_\alpha$ so that only its boundary
meets $\Sigma$ (at $w$ and $z$). The arc $\eta_\beta$ is obtained in
an analogous manner, only reversing the roles of the circles in $\alphas$
and $\betas$. The oriented knot $\OrK$ is gotten by the difference
$\eta_\alpha-\eta_\beta$. 

\subsection{Relative $\SpinC$ structures associated to intersection points}
\label{subsec:RelSpinCPoints}

Let $(\Sigma,\alphas,\betas,w,z)$ be a doubly-pointed Heegaard
diagram.  Fix $\x,\y\in \Ta\cap\Tb$. There are paths 
$$a\colon
[0,1]\longrightarrow \Ta, b\colon [0,1]\longrightarrow \Tb$$ with
$\partial a = \partial b = \x-\y$. These paths can be viewed as arcs
in $\Sigma$ (supported inside the $\alphas\cup\betas$). The difference
$a-b$ is a closed one-cycle in $\Sigma$ which is disjoint from $w$ and
$z$.  Indeed, since $\Sigma-w-z$ is a subset of $Y-K$, this one-cycle
can be viewed as an element ${\underline\epsilon}(\x,\y)\in
H_1(Y-K;\Z)$.

We construct a map $$\spincrel_{w,z}\colon \Ta\cap\Tb \longrightarrow
\RelSpinC(Y,K),$$ as follows (compare the analogous construction
from~\cite{HolDisk}).  A Heegaard diagram
$(\Sigma,\alphas,\betas,w,z)$ can be realized by a self-indexing Morse
function $f\colon Y \longrightarrow [0,3]$, with a single index zero
and three critical point (and $g$ index one and two critical points),
together with a Riemannian metric $g$, for which $\Sigma$ is the
mid-level $f^{-1}(3/2)$, $\alpha_i$ is the locus of points flowing out
of the $i^{th}$ index one critical point (via gradient flow), and
$\beta_j$ is the locus of points flowing into the $j^{th}$ index two
critical point.  Thus, $\x\in\Ta\cap\Tb$ can be thought of as a
$g$-tuple of gradient flow-lines $\gamma_\x$ containing all the index
one and two critical points.

Now, the knot $K$ is realized as a union of two flow-lines $\gamma_w$
and $\gamma_z$ which connect the index zero and index three critical
points, meeting $\Sigma$ in the points $w$ and $z$ respectively. 
The oriented knot $\OrK$ is gotten by $\gamma_z-\gamma_w$. 

We construct a vector field representing the relative $\SpinC$
structure $\spincrel_{w,z}(\xi)$ as follows. Modify the gradient
vector field $\grad f$ in a neighborhood of the flows $\gamma_\x$ so
that it has no zeros at any of the index one or two critical points.
This modification involves a choice of nowhere vanishing vector field
in a regular neighborhood of $\gamma_\x$, but it will follow easily
from the construction that this choice will not affect the relative
$\SpinC$ structure of the induced vector field.  Next, we modify the
vector field in a neighborhood of $\gamma_w$ to obtain a new vector
field $v$ which has no zeros at either the index zero or
three-critical points. In fact, this can be achieved so that the knot
$K$ is a closed orbit of the resulting vector field, whose orientation
as a flow-line agrees with the orientation induced from $\OrK$. This
$v$ modification involves a choice $X$ of nowhere vanishing vector
field on the neighborhood of $\gamma_w$, with fixed behaviour on
$K$. When calling attention to this choice, we write $v=v(X)$.  The
resulting vector field $v$ over $Y$ is the vector field representing
the $\SpinC$ structure $\spinc_w(\x)\in \SpinC(Y)$ associated to the
intersection point $\x$ and reference point $w$, as described
in~\cite{HolDisk} (and this fact is independent of the choice of $X$).

Our representative $v$ has been constructed so that there is a
neighborhood $D_2\times S^1$ of the closed flow-line $\{0\}\times
S^1\cong \OrK$ (where here $D_2$ is a disk of radius two centered at
the origin) with the property that the disks $D_2\times \theta$ for
any $\theta\in S^1$ are transverse to the vector field $v$. Consider a
concentric disk $D_1\subset D_2$. We can continuously extend
$v_2=v|_{Y-D_2\times S^1}$ to a new vector field $v_1=v_1(\x,X)$
over $Y-D_1\times S^1$, by using a vector field
over $(D_2-D_1)\times S^1$ which is everywhere transverse to the
annuli $(D_2-D_1)\times \theta$, and which point towards the origin at
$\partial D_1\times S^1$.  Thus, the vector field $v_1$
over $Y-(D_1\times S^1)$ inherits the
vector field $v'$ which points outwards at the boundary.  It is easy
to see that the isotopy class of $v_1$ is uniquely determined by the
isotopy class of our initial vector field $v$.

The induced relative $\SpinC$ structure $v_1$ over
$Y-\nbd(K)=Y-(D_1\times S^1)$ depends, of course, on $w$, $z$, $\x$,
and our choice $X$, and we write it correspondingly as
$v_1(w,z,\x,X)$. It is easy to see that
$$\Fill{Y,\OrK}(v_1(w,z,\x,X))=\spinc_w(\x).$$ 
It is easy to see that the induced relative $\SpinC$ structures over
$Y-\nbd{K}$ depends on the choice of $X$ by 
$$v_1(w,z,\x,X)-v_1(w,z,\x,X')=a\cm \PD[\mu],$$
where here $a\in\Z$ is a universal constant (depending on only $X$ and $X'$;
in fact, it is even independent of the ambient three-manifold $Y$).
We choose $X$ now to satisfy a normalization condition as follows.

Consider the unknot $K\subset S^3$. An orientation $\OrK$ specifies a
canonical relative $\SpinC$ structure
$\spincrel_0\in\RelSpinC(S^3,\OrK)$.  This is the relative $\SpinC$
structure represented by a vector field $v$ over $S^3-\nbd{\OrK}\cong
D^2\times S^1$ which is everywhere transverse to the disks $D^2$.  The
direction is specified by the condition that $v$ can be represented by
a vector field with closed orbits which have linking number one with
our original knot $\OrK$.  Our condition on $X$ now is that for the
standard genus one doubly-pointed diagram for the oriented unknot with
a single intersection point $\x$, the relative $\SpinC$ structure
induced by $v_1(w,z,\x,X)$ is the canonical $\SpinC$ structure for the
oriented unknot.

With this choice for $X$, it is easy now to see that the relative
$\SpinC$ structure underlying $v_1(w,z,\x,X)$ depends only on $w$,
$z$, and $\x$, inducing the required map $\relspinc_{w,z}\colon
\Ta\cap\Tb \longrightarrow \RelSpinC(Y,K)$.

We investigate its dependence on $w$, $z$, and $\x$ in the following lemma.
Continuing notation from~\cite{HolDisk}, letting $\x,\y\in\Ta\cap\Tb$,
let $\pi_2(\x,\y)$ denote the space of homotopy classes of Whitney
disks connecting $\x$ and $\y$, and for fixed 
$$p\in\Sigma-\alpha_1-...-\alpha_g-\beta_1-...-\beta_g$$
and $\phi\in\pi_2(\x,\y)$,
let $n_p(\phi)$ denote the intersection number of $\phi$ with 
the submanifold
$\{p\}\times\Sym^{g-1}(\Sigma)\subset \Sym^g(\Sigma)$.

\begin{lemma}
\label{lemma:SpinCRelChange}
Given $\x,\y\in\Ta\cap\Tb$, we have that
\begin{equation}
\label{eq:SpinCRelChangeGen}
\spincrel_{w,z}(\y)-\spincrel_{w,z}(\x)=
\PD[{\underline\epsilon}(\x,\y)].
\end{equation}
In particular, if there is some $\phi\in\pi_2(\x,\y)$, then
\begin{equation}
\label{eq:SpinCRelChange}
\spincrel_{w,z}(\x)-\spincrel_{w,z}(\y)=(n_z(\phi)-n_w(\phi))\cm
\PD[\mu].
\end{equation}
\end{lemma}

\begin{proof}
  We begin by establishing Equation~\eqref{eq:SpinCRelChangeGen}.
  Vector fields representing $\spincrel_{w,z}(\x)$ and
  $\spincrel_{w,z}(\y)$ can be chosen so that they agree everywhere
  except in a regular neighborhood of ${\underline\epsilon}(\x,\y)$.
  It follows that $\spincrel_{w,z}(\x)$ and $\spincrel_{w,z}(\y)$
  differ by some multiple of the Poincar\'e dual of this curve (in the
  case where the curve is one). The fact that this multiple is one
  follows from the analgous property of $\spinc_w(\x)$ established in
  Lemma~\ref{HolDiskOne:lemma:VarySpinC} of~\cite{HolDisk}.
  
  We turn to Equation~\eqref{eq:SpinCRelChange}. The homotopy class
  $\phi\in\pi_2(\x,\y)$ gives rise to a null-homotopy of
  ${\underline\epsilon}(\x,\y)$ inside $Y$. This null-homology meets
  the knot $K$ with intersection number $n_z(\phi)-n_w(\phi)$. Thus,
  it can be modified to give a homology of
  ${\underline\epsilon}(\x,\y)$ with $(n_z(\psi)-n_w(\psi))\cm
  \PD[\mu]$ in $Y-K$. Equation~\eqref{eq:SpinCRelChange} now follows
  from Equation~\eqref{eq:SpinCRelChangeGen}.
\end{proof}

\subsection{Heegaard triples and relative $\SpinC$ structures}
\label{subsec:HeegaardTriples}

A {\em Heegaard triple} is a closed, oriented two-manifold $\Sigma$,
equipped with three $g$-tuples of attaching circles,
$(\Sigma,\alphas,\betas,\gammas)$, c.f.~\cite{HolDisk}. This gives
rise to a four-manifold $X_{\alpha\beta\gamma}$ which has three
boundary components $-Y_{\alpha\beta}$, $Y_{\beta\gamma}$, and
$Y_{\alpha\gamma}$. (Of course, $Y_{\alpha\beta}$ here denotes the
three-manifold described by the Heegaard diagram
$(\Sigma,\alphas,\betas)$; $Y_{\alpha\gamma}$ and $Y_{\beta\gamma}$
are defined analogously.)

Suppose $K\subset Y$ is a framed knot. We can construct a Heegaard
triple $(\Sigma,\alphas,\betas,\gammas)$, where $Y_{\alpha\beta}$
represents $Y$, $Y_{\alpha\gamma}$ represents framed surgery along
$K$, and $Y_{\beta\gamma}$ is a connected sum of $g-1$ copies of
$S^2\times S^1$, and filling in $Y_{\beta\gamma}$ by a boundary
connected sum of $g-1$ copies of $B^3\times S^1$, we obtain a
four-manifold which is diffeomorphic to $W_\lambda(K)$. This is
obtained by first writing down a Heegaard triple with the property
that $K$ is supported entirely inside $U_\beta$, as the core of one of
the handles in the handlebody, and $\beta_g$ represents its meridian,
and realizing the framing $\lambda$ as a curve $\gamma_g$ which is
disjoint from the $\beta_i$ for $i=1,...,g-1$.  We then let $\gamma_i$
for $i=1,...,g-1$ be small, isotopic translates of the corresponding
$\beta_i$. We would like to choose convenient reference points. These
are points $w$ and $z$ are chosen so that there is an arc from $z$ to
$w$ which crosses none of the $\alpha_i$ or $\gamma_i$ for
$i=1,...,g$, or $\beta_j$ for $j=1,...,g-1$, so that an arc from $z$
to $w$ crosses $\beta_g$ transverally once. An ordering on $w$ and $z$
is specified by an orientation on $K$ We call such a doubly-pointed
Heegaard triple $(\Sigma,\alphas,\betas,\gammas,w,z)$ to be a {\em
  doubly-pointed Heegaard triple subordinate to the framed, oriented
  knot $\OrK\subset Y$ with framing $\lambda$}.

By constructing the diagram carefully, we can arrange for there to be
a unique intersection point $\Theta\in\Tb\cap\Tc$ representing the
generator $\Theta$ of the top-most non-trivial Floer homology group of
$\HFa(\#^{g-1}(S^2\times S^1))$. 
We denote this intersection point, and its corresponding homology class,
by $\Theta$.

As in~\cite{HolDisk}, we let $\pi_2(\x,\Theta,\y)$ denote the space of
homotopy classes of Whitney triangles, i.e. continuous maps of the
triangle into $\Sym^g(\Sigma)$ which map the three edges to $\Ta$,
$\Tb$, and $\Tc$, and the vertices to the given points. 
In Section~\ref{HolDisk:sec:HolTriangles} of~\cite{HolDisk}, it is
shown that homotopy classes of triangles give rise to $\SpinC$ structures
over $X_{\alpha,\beta,\gamma}$.

\begin{prop}
	\label{prop:FillW}
  Let $(\Sigma,\alphas,\betas,\gammas,w,z)$ be a doubly-pointed
  Heegaard triple subordinate to an oriented knot $\OrK\subset Y$,
  equipped with framing $\lambda$. The map
$$\psi\in\pi_2(\x,\Theta,\y) \mapsto \spincrel_{w,z}(\x)
+(n_w(\psi)-n_z(\psi))\mu$$ 
descends to give a well-defined map
$$\FillW{Y,\lambda,\OrK}\colon \SpinC(W_{\lambda}(K))\longrightarrow
\RelSpinC(Y,\OrK)$$
(i.e. which is independent of the choice of
Heegaard triple).  Moreover, letting $\CapSurf\in
H_2(W_\lambda(K),Y)\cong \Z$ denote a generator, we have that
\begin{equation}
\label{eq:PDF}
\FillW{Y,\lambda,\OrK}(\spinc+\PD[F_\lambda]) = \FillW{Y,\lambda,\OrK}(\spinc)+ \PD[K_\lambda].
\end{equation}
\end{prop}

\begin{proof}
  Recall~\cite{HolDisk} that two triangles
  $\psi\in\pi_2(\x,\Theta,\y)$ and $\psi'\in\pi_2(\x',\Theta,\y')$
  induce the same $\SpinC$ structure on $W_\lambda(K)$ is and only if
  there are $\phi_1\in\pi_2(\x,\x')$ (relative to the subspaces $\Ta$
  and $\Tb$) and $\phi_2\in\pi_2(\y',\y)$ (relative to the subspaces
  $\Ta$ and $\Tc$) with the property that
  $$\cald(\psi)-\cald(\psi')-\cald(\phi_1)-\cald(\phi_2)=n\cm [\Sigma]$$
  for some integer $n$.
  (Here, as in~\cite{HolDisk}, $\cald(\psi)$ denotes the two-chain in $\Sigma$
  induced by the homotopy class $\psi$, whose multiplicity at some point $p$
  is the intersection number of $\psi$ with $p\times \Sym^{g-1}(\Sigma)$.)
  But now $n_w(\phi_2)=n_z(\phi_2)$ (since $w$ and $z$ lie in the same
  component of $\Sigma-\alpha_1-...-\alpha_g-\gamma_1-...-\gamma_g$);
  combining this with the fact that
  $$\spincrel_{w,z}(\x')=\spincrel_{w,z}(\x)+(n_w(\phi_1)-n_z(\phi_1))\cm \mu$$
  (Lemma~\ref{lemma:SpinCRelChange}), it follows that in this case
  $$\spincrel_{w,z}(\x)+(n_w(\psi)-n_z(\psi))\mu =
  \spincrel_{w,z}(\x')+(n_w(\psi')-n_z(\psi'))\mu.$$
  
  Independence of the Heegaard triple is a routine consequence that
  any two Heegaard triples can be connected by a sequence of
  stabilizations and handleslides among the $\alpha_i$, $\beta_j$ for
  $j=1,...,g-1$, and the distinguished meridian $\beta_g$
  handlesliding over the other $\beta_j$.
  
  We now verify Equation~\eqref{eq:PDF}. Suppose that we have
  $\psi\in\pi_2(\x,\Theta,\y)$ and $\psi'\in\pi_2(\x',\Theta,\y)$.
  The difference $\cald(\psi)-\cald(\psi')$ gives a two-chain $C$ in
  $\Sigma$, which has no corners on $\gamma_i$ for $i=1,...,g$. The
  corners of $C$ occur at $x_i'\in\x'$ and $x_i\in\x$. Indeed, the
  boundary of the two-chain consists of ${\underline\epsilon}(\x,\x')$ 
  and some number of copies of the $\gamma_i$.
  As such, it can be thought of as furnishing a homology 
  $${\underline\epsilon}(\x,\x')=\ell\cm \gamma_g +
  (n_w(\psi-\psi')-n_z(\psi-\psi'))\mu$$
  in the knot complement (we
  are thinking of $\lambda$ as a curve in the Heegaard surface; of
  course, it is identified with $K_\lambda$ thought of as a curve in
  the knot complement).  Indeed, $\psi-\psi'$ corresponds to a
  surface-with-boundary $Z$ in the four-manifold underlying the
  Heegaard triple, whose boundary is supported entirely $Y$. The
  integer $\ell$ corresponds to the multiplicity with which $Z$ meets
  the cocore of the attaching two-handle. Thus,
  $\spinc_w(\psi)-\spinc_w(\psi')=\ell\cm \CapSurf$, and the equation
  follows.
\end{proof}

For a pointed Heegaard triple $(\Sigma,\alphas,\betas,\gammas,z)$, the
group of two-dimensional homology classes is identified with the group
of relations $a+b+c=0$ in $H_1(\Sigma;\Z)$ where where $a$ resp. $b$
resp. $c$ lies in the span of $\{[\alpha_i]\}_{i=1}^g$
resp. $\{[\beta_i]\}_{i=1}^g$ resp.  $\{[\gamma_i]\}_{i=1}^g$. These
relations correspond to two-chains $P$ in $\Sigma$ with boundary a
formal linear combination of the attaching circles, and with
$n_z(P)=0$. Given a Whitney triangle $\psi\in\pi_2(\x,\y,\w)$, there
is a formula for the evaluation $\langle
c_1(\spinc_z(\psi)),[P]\rangle$, where $[P]$ denotes the second
homology class corresponding to $P$, in terms of the Heegaard triple,
c.f. Proposition~\ref{HolDiskFour:prop:COneFormula}
of~\cite{HolDiskFour}. To describe this, we need two notions, the {\em
Euler measure} of the periodic domain and the {\em dual spider number}
of the triangle with respect to the triply-periodic domain.

Let $\Phi\colon F \longrightarrow \Sigma$ be a representative
for $P$, where here $\Phi$ is an immersion in a neighborhood of $\partial F$.
The line bundle
$\Phi^*(T\Sigma)$ has a canonical
trivialization over $\partial F$, since $\Phi$ induces an isomorphism
$$TF \cong \Phi^*(T\Sigma),$$ and $TF$ is canonically trivialized near
$\partial F$ (using the outward normal orientation on $F$). We define
$\EulerMeasure({\mathcal P})$ to be the Euler number of $\Phi^*(T\Sigma)$,
relative to this trivialization at $\partial F$,
$$\EulerMeasure({\mathcal P})=\langle c_1(\Phi^*(T\Sigma),\partial F), F\rangle.$$

Note first that the orientation of $\Sigma$, and the orientations on
all the attaching circles $\alphas$, $\betas$, and $\gammas$ naturally
induce ``inward'' normal vector fields to the attaching circles
(i.e. if $\gamma\colon S^1\longrightarrow \Sigma$ is a unit speed
immersed curve, this inward normal vector is given by
$J\frac{d\gamma}{dt}$).  Let $\alpha_i'$, $\beta_i'$, and $\gamma_i'$
denote copies of the corresponding attaching circles $\alpha_i'$,
$\beta_i'$, and $\gamma_i'$, translated slightly in these normal
directions. Let $\alphaprs$, $\betaprs$, and $\gammaprs$ denote the
corresponding $g$-tuples, and $\TaPr$, $\TbPr$, and $\TcPr$ be the
corresponding tori in $\Sym^g(\Sigma)$. By construction, then,
$u(e_{\alpha})$ misses $\TaPr$, $u(e_{\beta})$ misses $\TbPr$, and
$u(e_{\gamma})$ misses $\TcPr$.

Let $x\in\Delta$ be a point in the interior, chosen in general
position, so that the $g$-tuple $u(x)$ misses all of $\alphaprs$,
$\betaprs$, and $\gammaprs$. Choose three paths $a$, $b$, and $c$ from
$x$ to $e_{0}$, $e_{1}$, and $e_{2}$ respectively. The central
point $x$ and the three paths $a$, $b$, and $c$ is called a {\em dual
spider}. We can think of the paths $a$, $b$, and $c$ as one-chains
in $\Sigma$.
Recall that $\partial\PerDom$ has three types of boundaries: the
$\alpha$, $\beta$, and $\gamma$ boundaries, which we denote
$\partial_\alpha \PerDom$, $\partial_\beta \PerDom$, and
$\partial_\gamma\PerDom$.
Let 
$\partial'_\alpha\PerDom$, $\partial'_\beta\PerDom$, and
$\partial'_\gamma\PerDom$ respectively denote the one-chains obtained
by translating the corresponding boundary components using the induced
normal vector fields.
The dual spider number of $u$ and $\PerDom$ is defined by
$$\Spider(u,\PerDom)=n_{u(x)}(\PerDom)+\#(a\cap\partial'_\alpha\PerDom)+\#(b\cap\partial'_\beta\PerDom)+\#(c\cap
\partial'_\gamma\PerDom).$$

According to Proposition~\ref{HolDiskFour:prop:COneFormula}
of~\cite{HolDiskFour},
\begin{equation}
\label{eq:COneFormula}
\langle c_1(\spinc_z(u)), H(\PerDom)\rangle
=
\EulerMeasure(\PerDom) + \#(\partial\PerDom)-2n_z(\PerDom)+2 \Spider(u,\PerDom),
\end{equation}
where here $\#(\partial\PerDom)$ represents the number of boundary components
of $\PerDom$, with multiplicity.

\subsection{Filtered complexes}
\label{subsec:FiltCx}

A $\Z\oplus \Z$-filtered complex is a chain complex
$C$ whose underlying Abelian group decomposes as 
$$C=\bigoplus_{(i,j)\in\Z\oplus\Z} C\{i,j\},$$
and whose boundary operator
$\partial$ carries $C\{i_0,j_0\}$ into the subset
$$C\{i\leq i_0~\text{and}~ j\leq j_0\} = 
\bigoplus_{\{i',j'\big| i'\leq i~\text{and}~ j'\leq j\}}C\{i',j'\}\subset C.$$
A map $f\colon C\longrightarrow C'$ between filtered complexes is
called a {\em filtered map} if it carries $C\{i_0,j_0\}$ into
$C'\{i\leq i_0~\text{and}~ j\leq j_0\}$.  Two filtered chain maps
$$f_0, f_1\colon C\longrightarrow C'$$ are called {\em filtered
homotopic} if there is a filtered map $H\colon C \longrightarrow C'$
with $$f_0-f_1 = \partial' \circ H + H\circ \partial.$$ Two filtered
chain maps $C$ and $C'$ are called {\em filtered chain homotopy
equivalent} if there are filtered chain maps $f\colon C\longrightarrow
C'$ and $g\colon C'\longrightarrow C$ with the property that $f\circ
g$ and $g\circ f$ are filtered homotopic to the corresponding identity
maps.

A $\Z\oplus \Z$-filtered $\Z[U]$-complex $C$ is a filtered chain complex
equipped with an endomorphism $U\colon C\longrightarrow C$ whose restriction
to $C\{i,j\}$ maps to $C\{i-1,j-1\}$. 

If $C_1$ and $C_2$ are $\Z\oplus \Z$-filtered $\Z[U]$-chain complexes,
then we can form their tensor product $C_1\otimes_{\Z[U]} C_2$.  
This can be given a $\Z\oplus \Z$-filtration by 
$$(C\otimes C')\{i,j\} =
\frac{
\bigoplus_{\{(i_1,j_1), (i_2,j_2)\big| 
(i_1,j_1)+(i_2,j_2)=(i,j)\}}
C\{i_1,j_1\}\otimes C\{i_2,j_2\}}
{U_1\cm \xi_1 \otimes \xi_2\sim
\xi_1 \otimes U_2\cm \xi_2}.$$

If $C$ is a filtered complex, and $(a,b)\in \Z\oplus \Z$, let
$C[(a,b)]$ denote the filtered complex whose underlying chain complex
is isomorphic, but whose filtration is shifted by $(a,b)$; i.e.
$$C[(a,b)]\{i,j\} = C\{a+i,b+j\}.$$

A {\em relatively filtered map} $f\colon C \longrightarrow C'$ is a
chain map which respects the filtration on $C$ and the filtration
$C'[(a,b)]$ for some $(a,b)\in\Z\oplus \Z$.

\subsection{Absolute gradings}
\label{subsec:AbsGradings}

Heegaard Floer homology is natural under cobordisms.  Indeed, if $W$
is a smooth, connected, oriented cobordism with $\partial W = -Y_1\cup
Y_2$ which is equipped with a $\SpinC$ structure $\spinc$ whose
restriction $\spinct_i=\spinc|_{Y_i}$ for $i=1,2$ has torsion first
Chern class, then there is an induced chain map $$\fp_{W,\spinc}\colon
\CFp(Y_1,\spinct_1)
\longrightarrow \CFp(Y_2,\spinct_2)$$ which is homogeneous
of degree 
\begin{equation}
\label{eq:DimensionShift}
\frac{c_1(\spinc)^2-2\chi(W)-3\sigma(W)}{4}
\end{equation}
(c.f. Theorem~\ref{HolDiskFour:thm:AbsGrade} of~\cite{HolDiskFour}).

Let $\InjMod{}$ denote the module $\Z[U,U^{-1}]/U\cm
\Z[U]$. Recall~\cite{HolDisk} that $\HFp(S^3)\cong \InjMod{}$.
The $\Q$-grading on Floer homology is characterized by 
Equation~\eqref{eq:DimensionShift}, together with 
the normalization that $\HFp_d(S^3)$ is trivial for
all $d<0$, non-trivial in degree $d=0$.

Following our usual notational conventions, we write $\InjMod{(d)}$
for the module $\InjMod{}$, thought of as a graded $\Z[U]$-module,
where multiplication by $U$ lowers absolute degree by $2$, and the
non-zero homogeneous elements of lowest degree have degree $d$.
In this notation, then, $\HFp(S^3)\cong \InjMod{(0)}$.

\subsection{Approximating $\CFp$}
\label{subsec:CFd}

Following~\cite{IntSurg}, it is useful to have the following
approximation to $\CFp(Y)$: fix an integer $\delta\geq 0$, let
$\CFd(Y)\subset \CFp(Y)$ denote the subcomplex which is annihilated
by multiplication by $U^{\delta+1}$, and let $\HFd(Y)$ denote the 
homology of $\CFd(Y)$. In particular, for $\delta=0$,
this construction gives $\HFa(Y)$. Note that for all $\delta\geq 0$,
$\HFd(Y)$ is a three-manifold invariant.

\section{Construction of the knot filtration}
\label{sec:Construction}

The aim of this section is to construct the knot filtration: a
filtration of $\CFinf(Y)$ induced by a rationally null-homologous knot
$K\subset Y$. Most of this discussion is a straightforward generalization
of the corresponding constructions for null-homologous knots $K\subset Y$
described in~\cite{Knots} and~\cite{RasmussenThesis}.

Let $K\subset Y$ be a rationally null-homologous knot, and let
$(\Sigma,\alphas,\betas,w,z)$ be a corresponding doubly-pointed
Heegaard diagram. 

For fixed $\xi\in\SpinC(Y,K)$, let $\EquivClass\subset
(\Ta\cap\Tb)\times \Z\times \Z$ be the subset of elements
$[\x,i,j]$ satisfying
\begin{equation}
\label{eq:GenCFinfYK}
\spincrel_{w,z}(\x)+(i-j)\cm \PD[\mu]=\xi.
\end{equation}
According to Lemma~\ref{lemma:SpinCRelChange}, if 
$[\x,i,j]\in\EquivClass$, and $\phi\in\pi_2(\x,\y)$,
then $[\y,i-n_w(\phi),j-n_z(\phi)]\in\EquivClass$.

Let 
$\CFKinf(\Sigma,\alphas,\betas,w,z,\xi)$ be the chain complex generated
by $[\x,i,j]\in \EquivClass$, endowed with the differential
$$\partial^\infty[\x,i,j]= \sum_{\y\in\Ta\cap\Tb}
\sum_{\{\phi\in\pi_2(\x,\y)\big|\Mas(\phi)=1\}}
\#\left(\frac{\ModFlow(\phi)}{\R}\right)\cm [\y,i-n_w(\phi),j-n_z(\phi)],$$
where as usual, $\ModFlow(\phi)$ is the moduli
space of pseudo-holomorphic representatives of $\phi$, and
$\Mas(\phi)$
denotes its expected dimension.

The map $$\Filt\colon \EquivClass\longrightarrow \Z\oplus\Z$$ given by
$\Filt[\x,i,j]=(i,j)$ induces a $\Z\oplus \Z$ filtration on
$\CFKinf(\Sigma,\alphas,\betas,w,z)$.

\begin{theorem}
The filtered chain homotopy type of
$\CFKinf(\Sigma,\alphas,\betas,w,z,\xi)$ is an invariant
of the underlying oriented knot $\OrK\subset Y$ and choice of 
relative $\SpinC$ structure $\xi\in\SpinC(Y,K)$. 
\end{theorem}

\begin{proof}
This is a routine adaptation of the corresponding statement
in Theorem~\ref{Knots:thm:KnotInvariant} of~\cite{Knots}.
\end{proof}

We denote this complex by $\CFKinf(Y,\OrK,\xi)$.

\begin{prop}
        \label{prop:Filtrations}
        Let $C_{\xi}=\CFKinf(Y,\OrK,\xi)$.
        We have that
        $$C_{\xi+\PD[\mu]}=C_{\xi}[(0,-1)].$$
        The short exact sequences of $\Z[U]$-chain complexes
        $$
        \begin{CD}
        0@>>>C_{\xi}\{i<0\}@>>>C_{\xi}@>>>C_{\xi}\{i\geq 0\}@>>> 0 
        \end{CD}
        $$
        and
        $$
        \begin{CD}
        0@>>> C_{\xi}\{i= 0\} @>>>
        C_{\xi}\{i\geq 0\}
        @>{U}>>
        C_{\xi}\{i\geq 0\}
        @>>> 0
        \end{CD}
        $$
        are isomorphic to the 
        short exact sequences
        $$
        \begin{CD}
        0@>>>\CFm(Y,\spinc) @>>>
        \CFinf(Y,\spinc) @>>>
        \CFp(Y,\spinc) @>>>0,
        \end{CD}
        $$
        and
        $$
        \begin{CD}
        0@>>>\CFa(Y,\spinc)@>>>\CFp(Y,\spinc)@>{U}>> \CFp(Y,\spinc)@>>>0
        \end{CD}
        $$      
        respectively, where here $\spinc=\Fill{Y,\OrK}(\xi)$.
        Similarly, 
        the short exact sequences of $\Z[U]$-chain complexes
        $$
        \begin{CD}
        0@>>>C_{\xi}\{j<0\}@>>>C_{\xi}@>>>C_{\xi}\{j\geq 0\}@>>> 0 
        \end{CD}
        $$
        and
        $$
        \begin{CD}
        0@>>> C_{\xi}\{j= 0\} @>>>
        C_{\xi}\{j\geq 0\}
        @>{U}>>
        C_{\xi}\{j\geq 0\}
        @>>> 0
        \end{CD}
        $$
        are isomorphic to the 
        short exact sequences
        $$
        \begin{CD}
        0@>>>\CFm(Y,\spinc') @>>>
        \CFinf(Y,\spinc') @>>>
        \CFp(Y,\spinc') @>>>0,
        \end{CD}
        $$      
        and
        $$
        \begin{CD}
        0@>>>\CFa(Y,\spinc') @>>>
        \CFp(Y,\spinc') @>{U}>>
        \CFp(Y,\spinc') @>>>0,
        \end{CD}
        $$      
        where here $\spinc'=\Fill{Y,-\OrK}(\xi)$.
\end{prop}

The set of relative $\SpinC$ structures $\Fill{Y,\OrK}^{-1}(\spinc)$
inducing a fixed $\SpinC$ structure
over $Y$ is a well-ordered set, under the rule that
$\xi_1\leq \xi_2$ if $\xi_2=\xi_1 + j \cm \PD[\mu]$ for some $j\geq 0$.

This ordering on relative $\SpinC$ structures gives rise to an
ordering of the generators of $\CFa(Y)$ by $\xi\mapsto
\spincrel_{w,z}(\xi)$. It is easy to see that this ordering gives
rise to a filtration of 
the complex $\CFa(Y)=\bigoplus_{\x\in\Ta\cap\Tb} \Z\x$ endowed
with the usual differential $$\partial \x =
\sum_{\{\y\in\Ta\cap\Tb\big| n_w(\phi)=0, \Mas(\phi)=1\}}
\#\left(\frac{\ModFlow(\phi)}{\R}\right)\cm \y.$$
The homology of the associated graded object is the knot Floer
homology 
$$\HFKa(Y,\OrK)=\bigoplus_{\xi\in\RelSpinC(Y,\OrK)} \HFKa(Y,\OrK,\xi),$$
where $\HFKa(Y,\OrK,\xi)$ is generated by $\x$ with $\spincrel_{w,z}(\x)=\xi$,
and whose differential counts holomorphic disks with
$n_w(\phi)=n_z(\phi)=0$.

\subsection{Relationship with knot Floer homology for null-homologous knots}
\label{subsec:OlderNotation}

If $K\subset Y$ is a null-homologous knot, then a choice of Seifert
surface for $Y$ gives an identification $\SpinC(Y,K)\cong \Z\oplus
\SpinC(Y)$. Thus, with this additional choice, we identify the
$\RelSpinC(Y)$-filtration of $\CFa(Y,\spinc)$
with a $\Z$-filtration. This identification is used, for example,
in~\cite{Knots}, where the knot filtration is described as a filtration by
$\Z$, rather than by relative $\SpinC$ structures.

For example, if $\OrK$ is a knot in $S^3$ and $s$ is an integer, then
$\HFKa(S^3,\OrK,\xi)$ in the present notation corresponds to
$\HFKa(Y,\OrK,s)$ for $s\in\Z$ in the notation from the introduction
or ~\cite{Knots}, where here $s$ and $\xi$ are related by
$c_1(\xi)=2s\PD[\mu]$.

\section{Knot Floer homology and large surgeries}
\label{sec:Large}

We describe here the result of forming ``large surgeries'' on a
rationally null-homologous, framed knot $K\subset Y$. This result
generalizes a corresponding result in the null-homologous case,
c.f.~\cite{Knots} and \cite{RasmussenThesis}, and it is a special case of
Theorem~\ref{thm:SurgeryFormula}

Let $Y$ be a knot given with framing $\lambda$. We can form a
three-manifold $Y_{m\mu+\lambda}(K)$ obtained by filling the curve
$m\cm \mu+\lambda$ in $Y-\nbd{K}$.  Let $$W_m'(K) \colon
Y_{m\mu+\lambda}(K) \longrightarrow Y$$
denote the two-handle
cobordism obtained by turning around the two-handle cobordism from
$-Y$ to $-Y_{m\cm\mu+\lambda}(K)$. Note that $$H_2(W_m'(K),Y;\Z)\cong
\Z,$$
and let $\CapSurf\subset W_m'(K)$ be a surface-with-boundary representing a
generator.  Clearly, $\PD[\CapSurf]|_{Y}=\PD[K]\in H^2(Y;\Z)$.  Note
that for sufficiently large $m$, the self-intersection number of
$\CapSurf$ is negative.

Let $\OrK\subset Y$ be an oriented knot in a three-manifold $Y$,
whose induced homology class is trivial in rational homology.
We fix also a framing $\lambda$ of $K$.
For a fixed relative $\SpinC$ structure $\xi\in\RelSpinC(Y,K)$,
let $C_\xi$ be the $\Z\oplus\Z$-filtered chain complex
$\CFKinf(Y,\OrK,\xi)$. 
There are two projection maps
\begin{eqnarray}
C_\xi\{\max(i,j)\geq 0\} \longrightarrow C_\xi\{i\geq 0\}
&{\text{and}}&
C_\xi\{\max(i,j)\geq 0\} \longrightarrow C_\xi\{j\geq 0\}.
\label{eq:CanonProj}
\end{eqnarray}
Denoting
\begin{eqnarray}
\Ap_\xi(Y,\OrK)=C_\xi\{\max(i,j)\geq 0\}&{\text{and}}&
\Bp_\xi(Y,\OrK)=\CFp(Y,\Fill{Y,\OrK}(\xi))
\label{eq:DefAp}
\end{eqnarray}
and using the identifications
(from Proposition~\ref{prop:Filtrations})
\begin{eqnarray*}
C_\xi\{i\geq 0\}\cong \CFp(Y,\Fill{Y,\OrK}(\xi))
&{\text{and}}&
C_\xi\{j\geq 0\}\cong \CFp(Y,\Fill{Y,-\OrK}(\xi))\cong 
\CFp(Y,\Fill{Y,\OrK}(\xi+\PD[K_\lambda]),
\end{eqnarray*}
where here $K_\lambda$ is the push-off of $K$ inside $Y-K$ using
the framing $\lambda$,
we can view the canonical projection maps of Equation~\eqref{eq:CanonProj}
as maps
\begin{eqnarray}
\label{eq:DefOfVH}
\vertp_\xi\colon \Ap_\xi(Y,\OrK) \longrightarrow \Bp_\xi(Y,\OrK)
&{\text{and}}&
\horp_\xi\colon \Ap_\xi(Y,\OrK) \longrightarrow \Bp_{\xi+\PD[K_\lambda]}(Y,\OrK).
\end{eqnarray}

\begin{theorem}
\label{thm:LargeNSurgeries}
Let $K\subset Y$ be a rationally null-homologous knot in a closed,
oriented three-manifold, equipped with a framing $\lambda$.  Then, for
all sufficiently large $m$, there is a map $$\Xi \colon \SpinC(Y_{m\cm\mu + \lambda}(K))
\longrightarrow
\RelSpinC(Y,K)$$ 
with the property that for all
$\spinct\in\SpinC(Y_{m\mu+\lambda}(K))$, the group
$\CFp(Y_{m\cm\mu+\lambda}(K),\spinct)$ is represented by the chain
complex 
$$\Ap_{\Xi(\spinct)}=C_{\xi}\{\max(i,j)\geq 0\},$$ in the
sense that there are isomorphisms (of relatively $\Z$-graded
$\Z[U]$-complexes) $$\Psi^+_{\spinct,m}\colon
\CFp(Y_{m\cm\mu+\lambda}(K),\spinct)
\longrightarrow \Ap_{\Xi(\spinct)}(Y,\OrK).$$ 

Furthermore, fix $\spinct\in\SpinC(Y_{m\mu+\lambda}(K))$, 
and let $\xi=\Xi(\spinct)$. There are $\SpinC$ structures
$\spincx,\spincy\in W'_{m}(K)$ with 
$\FillW{Y,m\mu+\lambda,\OrK}(\spincx)=\xi$,
and $\spincy=\spincx+\PD[\CapSurf]$ with the property that 
the maps 
$\vertp_\xi$ and $\horp_\xi$ correspond to the maps
induced by the cobordism $W'_m(K)$ equipped with $\spincx$ and $\spincy$
respectively. More precisely,
the following squares commute:
\[\begin{CD}
\CFp(Y_{m\cm\mu+\lambda}(K),\spinct) @>{\fp_{W'_m(K),\spincx_s}}>> \CFp(Y,\Fill{Y,\OrK}(\xi))) \\
@V{\Psi^+_{\spinct,m}}VV @VV{=}V \\
\Ap_\xi(Y,\OrK) @>{\vertp_\xi}>> \Bp_\xi(Y,\OrK)
\end{CD}
\]
and
\[\begin{CD}
\CFp(Y_{m\mu+\lambda}(K),\spinct) @>{\fp_{W'_m(K),\spincy_s}}>> \CFp(Y,\Fill{Y,-\OrK}(\xi)) \\
@V{\Psi^+_{\spinct,m}}VV @VV{=}V \\
\Ap_\xi(Y,\OrK) @>{\horp_\xi}>> \Bp_{\xi+\PD[\OrK]}(Y,\OrK).
\end{CD}
\]
\end{theorem}

The following result is also easy consequences of the proof:

\begin{prop}
  \label{prop:OnlyTwo}
  Let $K\subset Y$ be a rationally null-homologous knot in a closed,
  oriented three-manifold, equipped with framing $\lambda$. For any
  $\delta>0$, there is an integer $N$ with the property that for all
  $m\geq N$ and all $\spinct\in \SpinC(Y_{m\cm\mu+\lambda}(K))$, there
  are at most two $\SpinC$ structures in
  $\SpinC(W_{m\cm\mu+\lambda}(K))$ with restriction to
  $Y_{m\cm\mu+\lambda}(K)$ equal to $\spinct$ for which the induced map
  $$ F^{\delta}_{W'_{m\cm\mu+\lambda}(K),\spinc}\colon \HFd(Y_{m\cm\mu+\lambda}(K),\spinct)\longrightarrow
  \HFd(Y)$$
  is non-trivial. These are the $\SpinC$ structures $\spincx$ and $\spincy$
  associated to $\spinct$
  from Theorem~\ref{thm:LargeNSurgeries} above.
\end{prop}

We return to the proofs of Theorem~\ref{thm:LargeNSurgeries} 
and Proposition~\ref{prop:OnlyTwo}
after some preliminary discussion and lemmas.

\begin{figure}
\mbox{\vbox{\epsfbox{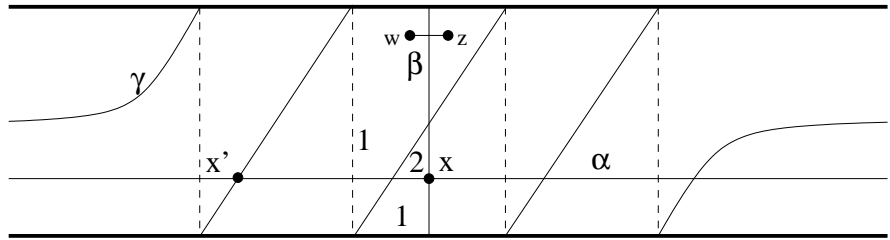}}}
\caption{\label{fig:LargeSurgeries}
{\bf{Illustration of the Heegaard triple}} The integers denote
(non-zero) local multiplicities of the ``small triangle'' 
connecting $x$ and $x'$. This picture is taking place in a cylindrical
region in $\Sigma$.}
\end{figure}

We work with a family of doubly-pointed Heegaard triples for the
framed knot $\OrK$ $(\Sigma,\alphas,\betas,\gammas,w,z)$, so that
there are identifications $Y_{\alpha,\gamma}\cong
Y_{m\mu+\lambda}(K)$, $Y_{\beta,\gamma}\cong \#^{g-1}(S^2\times S^1)$,
$Y_{\alpha,\beta}\cong Y$. We can give Heegaard triples for all of the
$W'_m$ which differ only in $\gamma_g$, which winds along the meridian
$\mu$. We call this region the ``winding region'' (c.f.
Figure~\ref{fig:LargeSurgeries} for an illustration of a winding
region; in this picture, the subscript for $\gamma_g$ is dropped).

\begin{defn}
\label{def:SmallTriangles}
An intersection point $\x\in\Ta\cap\Tc$ is said to be {\em supported
in the winding region} if its component along $\gamma_g$ lies in the
winding region. Given $\x,\y\in\Ta\cap\Tc$ (both supported in the
winding region), we say that $\phi\in\pi_2(\x,\y)$ is supported in the
if segment of $\partial(\cald(\phi))$ in $\gamma_g$ is a subset of
the winding region. (Note that if $\phi\in\pi_2(\x,\y)$ is supported in
the winding region, then so is any $\phi'\in\pi_2(\x,\y)$.) An
equivalence class of intersection point for $Y_{m\cm\mu+\lambda}$ is
said to be supported in the winding region if every intersection point
in the equivalence class is supported in the winding region, and any
two intersection points can be connected by Whitney disks supported in
the winding region.  Finally, $\psi\in\pi_2(\x,\Theta,\y)$ is said to
be a {\em small triangle} if the $\gamma_g$ arc in $\partial \cald(\psi)$ 
is supported in the winding
region. Note that in this case, $\y$ is the ``closest point'' in
$\Ta\cap\Tb$ to $\x$.
\end{defn}

\begin{defn}
  The {\em depth} of a basepoint $z$ in the winding region is the
  minimum absolute value of the algebraic intersection number of
  $\gamma_g$ with any arc $A\subset \Sigma-\alpha_1-...-\alpha_g$
  which connects $z$ with some basepoint $z_0$ outside the winding
  region with $z$.  Given an integer $\epsilon>0$, we say that $z$ is
  $\epsilon$-centered if its depth is less than or equal to
  $(m-\epsilon)/2$. Similarly, a choice of meridian $\mu$ is called
  $\epsilon$-centered if each point $z\in\mu$ is $\epsilon$-centered.
\end{defn}

For a  Heegaard diagram with an $\epsilon$-centered meridian, with
basepoints $w$ and $z$ on either side of $\mu$. Then,
for each small triangle $\psi$, we have that at least one of $n_w(\psi)$ or $n_z(\psi)=0$, and also
\begin{equation}
\label{eq:BoundSmallTriangle}
\max(|n_w(\psi)|, |n_z(\psi)|)\leq \frac{m-\epsilon}{2}.
\end{equation}
This is an immediate consequence of the definitions.

\begin{lemma}
\label{lemma:SupportedInRegion}
There is an integer $\epsilon>0$ with the property that for all
sufficiently large $m$, each $\spinct\in\SpinC(Y_{m\mu+\lambda}(K))$
can be represented by an equivalence class of intersection points
supported in the winding region and an $\epsilon$-centered choice
of meridian.
\end{lemma}

\begin{proof}
If an equivalence class of intersection points in $\Ta\cap\Tc$ is not
supported in the winding region, we say it is {\em bad}.  It is clear
that for all large $m$, the set of bad equivalence classes is bounded.
We choose $\epsilon$ so that $2\epsilon$ is greater than this number.
Now move the basepoint as needed.
\end{proof}

Given $\xi\in\RelSpinC(Y,K)$, define
$$\Psi^\infty_{\xi}\colon \CFinf(Y_{m \cm \mu + \lambda} (K),\spinct)
\longrightarrow \CFinf(Y,\OrK,\xi)$$ by 
$$\Psi^\infty_{\xi}[\x,i] =
\sum_{\y\in\Ta\cap\Tb}
\sum_{\{\psi\in\pi_2(\x,\Theta,\y)\big|
\FillW{W'_m(K)}(\spinc_w(\psi))=\xi)\}}
\left(\#\ModFlow(\psi)\right)\cm [\y,i-n_w(\psi),i-n_z(\psi)],$$
where here $\x$ represents $\spinc_w(\psi)|_{Y_{m\mu+\lambda}(K)}$. 
We claim that this is a chain map.

Given $\spinct\in\SpinC(Y_{m\mu+\lambda}(K))$, we realize $\spinct$ by an
equivalence class of intersection points supported in the winding
region. Define $\xi=\Xi(\spinct)$ by
$$\xi=\spincrel_{w,z}(\x)+ (n_w(\psi)-n_z(\psi))\cdot \mu,$$
where $\psi\in\pi_2(\x',\Theta,\x)$ is any small
triangle connecting an intersection point $\x'\in\Ta\cap\Tc$
representing $\spinct$ with its ``nearest point'' $\x\in\Ta\cap\Tb$.
The restriction of $\Psi^\infty_\xi$ to $\CFp(Y_{m\mu+\lambda}(K),\spinct)$ induces
a chain map
$$\Psi^+_{\xi}\colon 
\CFp(Y_{m\mu+\lambda}(K),\spinct)\longrightarrow
\Ap_\xi(Y,\OrK).$$

\vskip.2cm
\noindent{\bf{Proof of Theorem~\ref{thm:LargeNSurgeries}.}}
The argument from~\cite{Knots} shows
that $\Psi^+$ induces an isomorphism. Post-composing $\Psi^+_{}$
with the projection $C_\xi\{\max(i,j)\geq 0\}$ to $C\{i\geq 0\}$
(i.e. $\vertp$), we obtain the map induced by the cobordism
$W'_{m}(K)$ equipped with the $\SpinC$ structure $\spincx_s$; i.e.
the first square in the statement of the theorem
commutes. Commutativity of the second  square follows similarly.
\qed
\vskip.2cm

\begin{lemma}
  \label{lemma:CalcEvaluation}
        There is a constant $c$ with the property that for
        all sufficiently large $m$, for any of the $\SpinC$ structures
        $\{\spincx(\spinct),\spincy(\spinct)\}_{\spinct\in\SpinC(Y_{m\mu+\lambda}(K))}$
        appearing in Theorem~\ref{thm:LargeNSurgeries}, we have that
        \[\begin{array}{rl}
        -c&\leq \langle c_1(\spincx),[\CapSurf]\rangle \leq 2m+c \\
        -2m-c&\leq \langle c_1(\spincy),[\CapSurf]\rangle \leq c. \\
        \end{array}\]
\end{lemma}

\begin{proof}
There is a constant $C$ (depending on $\x\in\Ta\cap\Tb$, but independent
of $m$) with the
property that for any small triangle $\psi\in\pi_2(\x',\Theta,\x)$,
\begin{equation}
\label{eq:ChernFormula}
\langle c_1(\spinc_w(\psi)),[\CapSurf]\rangle + [\CapSurf]\cdot[\CapSurf]
= C + 2(n_w(\psi)-n_z(\psi)). \end{equation} 

This can be seen as follows. Consider the function $$f(\x)=\langle
c_1(\spinc_w(\psi)),[\CapSurf]\rangle-2(n_w(\psi)-n_z(\psi)),$$
where
here $\psi\in\pi_2(\x',\Theta,\x)$ is a small triangle.  We claim that
this is independent of the choice of $\x'$.  This can be seen by
varying $\x'$, and appealing to Equation~\eqref{eq:COneFormula}.
Similarly, to verify its independence of $m$, it suffices to consider
a fixed small triangle $\psi\in\pi_2(\x,\Theta,\x')$, and observe that
$\langle c_1(\spinc_w(\psi)),[\CapSurf]\rangle $ increases by one as
the number $m$ is increased by one, and furthermore
$[\CapSurf]\cdot[\CapSurf]$ also decreases by one.
Note here that by
``fixed small triangle'', we mean that $\x$ is fixed, as is
$n_w(\psi)-n_z(\psi)$, though, of course, $\gamma_g$ is varied. This
assertion is an easy consequence of Equation~\eqref{eq:COneFormula}.
More specifically, with the integer $d$ chosen as in
Equation~\eqref{eq:OrderD}, we see that $d\CapSurf$ can be represents
a generator of $H_2(W'_m(K);\Z)$. Specifically, 
let $\PerDom_m$ be a generator for the space of triply-periodic domains, and write
$\partial [\PerDom_m]=A+B+C$,
where $A$, $B$ and $C$ are first
homology classes in the spans of $\alphas$, $\betas$, and $\gammas$
respectively. In view of Equation~\eqref{eq:OrderD}, $\partial
\PerDom_m$ has $C$ component given by $d(m\mu+\lambda)$ and $B$
component $-(dm+n)\mu$, modulo the other $\beta_i$ for $i<g$. Thus,
$\# (\partial \PerDom_m)=dm+c$ for some constant $c$.  It is easy to
see that the other quantities in Equation~\eqref{eq:COneFormula} are
independent of $m$.

Thus, we have seen that for any small triangle $\psi\in\pi_2(\x',\Theta,\x)$,
$$\langle
c_1(\spinc_w(\psi)),[\CapSurf]\rangle + [\CapSurf]\cdot[\CapSurf]
-2(n_w(\psi)-n_z(\psi))$$ depends only on 
$\x\in\Ta\cap\Tb$.  Since there are only finitely many 
intersection points $\Ta\cap\Tb$, 
it follows that there is some constant
$C'$ with the property that
\begin{equation}
\label{eq:BoundCOne}
\Big|\langle c_1(\spinc_w(\psi),[\CapSurf]\rangle +
[\CapSurf]\cdot[\CapSurf] -2(n_w(\psi)-n_z(\psi))\Big| \leq C'.
\end{equation}
By restricting to $\epsilon$-centered base points for some
integer $\epsilon$ independent of $m$ (which we can do in light of
Lemma~\ref{lemma:SupportedInRegion}), we can arrange for
$2|n_w(\psi)-n_z(\psi)|\leq m-\epsilon$ for all small triangles
(c.f. Equation~\eqref{eq:BoundSmallTriangle}).
The result now follows, bearing in mind that
$\spincx=s_w(\psi)$, and 
$\spincy=\spinc_z(\psi)=\spincx(\psi)-\PD[\CapSurf]$.
\end{proof}

\vskip.2cm
\noindent{\bf Proof of Proposition~\ref{prop:OnlyTwo}.}
We choose $m$ large enough that $W_m'(K)$ has negative-definite
intersection form.  Given $\spinct\in\SpinC(Y)$, let ${\mathfrak
  S}(\spinct)\in\SpinC(W_m'(K))$ denote the set of $\SpinC$ structures
whose restriction to $Y$ is $\spinct$.  This set ${\mathfrak S}(t)$
has the form $\{\spinc_0+n \cm \PD[\CapSurf]\}_{n\in\Z}$ for some
fixed $\spinc_0\in\SpinC(W_{m}'(K))$. The function $n\mapsto
c_1(\spinc_0+n\cm \PD[\CapSurf])^2$ is a quadratic function of $n$
which is bounded above.

Choosing $m$ larger than the constant $c$ from
Lemma~\ref{lemma:CalcEvaluation}, it follows easily that
$c_1(\spincy)^2 \geq c_1(\spincy+\PD[\CapSurf])^2$,
while 
$c_1(\spincx)^2\geq c_1(\spincx-\PD[\CapSurf])^2$.
Since $\spincy=\spincx+\PD[\CapSurf]$, it follows readily that 
at least
one of $c_1(\spincx)^2$ or $c_1(\spincy)^2$ is a maximum of $c_1(\spinc)^2$
for $\spinc\in{\mathfrak S}(\spinct)$.

As in the proof of Lemma~\ref{lemma:CalcEvaluation} (cf.
Inequality~\eqref{eq:BoundCOne}), there is a constant $B$ 
independent of $m$
with
the property that if
$\langle c_1(\spincx),[\CapSurf]\rangle + [\CapSurf]\cm[\CapSurf] \geq B$, 
then for any corresponding
small triangle $\psi$ representing $\spincx$, we have $n_w(\psi)>0$.
Moreover, in this case, $\horp$, which corresponds to the $\SpinC$
structure $\spincy$, induces an isomorphism, while $\spincy$ maximizes
$c_1(\spinc)^2$ among all $\spinc\in{\mathfrak S}(\spinct)$. In the
same way, when $\langle c_1(\spincx),[\CapSurf]\rangle +[\CapSurf]\cm[\CapSurf]
\leq -B$, then $\spincx$ induces an isomorphism,
and it maximizes $c_1(\spinc)^2$ among all $\spinc\in{\mathfrak
  S}(\spinct)$.  

In either case, in view of Equation~\eqref{eq:DimensionShift}, we see
that the degree of any element of $\HFa(Y_{m\mu+\lambda}(K),\spinct)$
lies within a bounded distance (independent from $m$ and $\spinct$)
from $-{c_1(\spinc)^2}/{4}$, where $\spinc\in{\mathfrak S}(\spinct)$
minimizes $c_1(\spinc)^2$. The claimed result follows.  \vskip.2cm

Proposition~\ref{prop:OnlyTwo} has the following consequence:

\begin{cor}
  \label{cor:LargeNFloerHomology}
  Let $K\subset Y$ be a null-homologous knot in an integer homology
  three-sphere, and fix an integer $\delta\geq 0$. There is a constant $C$
  with  the property that for all sufficiently large $m$
  and any $\spinct\in\SpinC(Y_{m\mu+\lambda}(K))$, there is a chain complex
  for $\CFd(Y_{m\mu+\lambda}(K),\spinct)$ with the property that
  if $\CFd_d(Y_{m\mu+\lambda}(K),\spinct)$ is non-trivial, then 
  $$|d-\frac{m}{4}|\leq C.$$
\end{cor}

\begin{proof}
  It is easy to see that for fixed $\spinct\in\SpinC(Y)$, if we
  consider ${\mathfrak S}(\spinct)$, the element $\spinc_0$ which
  maximizes $c_1(\spinc)^2$ has $c_1(\spinc_0)^2=-m+c$, where here $c$
  is some constant (independent of $m$). Thus, according to
  Equation~\eqref{eq:DimensionShift}, the map
  $\fd_{W_\lambda(K),\spinc_0}$ carries an element of degree $d$ to an
  element of $d-\frac{m}{4}+c$. Now, according to
  Theorem~\ref{thm:LargeNSurgeries}, we have a chain complex
  representing $\CFd(Y_{m\mu+\lambda}(K),\spinct)$ (for any
  choice of $\spinct)$ whose breadth is constant, indepedent of $m$
  and the choice of $\spinct$. According to the proof of
  Proposition~\ref{prop:OnlyTwo}, $\fd_{W_\lambda(K),\spinc_0}$,
  the map of degree $d-\frac{m}{4}+c$,
  carries some element of this complex non-trivially to $\CFd(Y)$.
\end{proof}

\section{K{\"u}nneth principle}
\label{sec:Kunneth}

Let $\OrK_1\subset Y_1$ and $\OrK_2\subset Y_2$ be a pair of
oriented three-manifolds equipped with oriented knots. Then, 
we can form the connected sum to obtain an oriented knot
$\OrK_1\# \OrK_2\subset Y_1\# Y_2$.
Indeed, given $\xi_i\in\RelSpinC(Y_i,K_i)$,
we can form their connected sum $\xi_1\#\xi_2$. 
This induces a gluing map $$\RelSpinC(Y_1,K_1)\times
\RelSpinC(Y_2,K_2)\longrightarrow \RelSpinC(Y_1\# Y_2,K_1\#K_2),$$
written $\xi_1, \xi_2\mapsto \xi_1\#\xi_2$.
which is equivariant with respect to the natural map
$$H^2(Y_1,K_1;\Z)\oplus H^2(Y_2,K_2;\Z) \longrightarrow H^2(Y_1\# Y_2,
K_1\# K_2;\Z).$$

More explicitly, we can think of $\xi_i$ ($i=1,2$) as specified by a
non-vanishing vector field which contains $\OrK_i$ as a closed orbit.
We realize the connected sum as attaching a one-handle to $Y_1\coprod
Y_2$ along a pair of three-balls $B_i$ supported on $\OrK_i$. We can
assume that on each sphere $S_i=\partial B_i$, there are eactly two
points where the vector field is normal to $S_i$, the two points where
$S_i$ meets $\OrK_i$: the ``in-going'' and ``out-going points''. Here,
$p\in S_i$ is ``in-going'' if $\xi_i$ points into $Y_i-B_i$. We can
then match the in-going point on $Y_1-B_1$ with the out-going point on
$Y_2-B_2$ and vice versa, to construct a nowhere vanishing vector
field on $Y_1\#Y_2$. This new vector field has $K_1\# K_2$ as a closed
orbit, and it has prescribed from in the connected sum region. This
vector field gives rise to a relative $\SpinC$ structure
$\SpinC(Y_1\#Y_2,\OrK_1\#\OrK_2)$, inducing the gluing map above.

The gluing map can be described in terms of Heegaard diagrams as
follows.  Let $(\Sigma_1,\alphas_1,\betas_1,w_1,z_1)$ and
$(\Sigma_2,\alphas_2,\betas_2,w_2,z_2)$ be the doubly-pointed Heegaard
diagrams compatible with the oriented knots $\OrK_1\subset Y_1$ and
$\OrK_2\subset Y_2$. Consider the oriented two-manifold $\Sigma_3$
obtained as the connected sum of $\Sigma_1$ with $\Sigma_2$,
identifying neighborhoods of $z_1$ with $w_2$. Then, the
doubly-pointed Heegaard diagram
$(\Sigma_3,\alphas_1\cup\alphas_2,\betas_1\cup\betas_2,w_1,z_2)$ is
compatible with $\OrK_1\#\OrK_2\subset Y_1\# Y_2$. Now, given
$\x_1\in{\mathbb T}_{\alpha_1}\cap{\mathbb T}_{\beta_1}$ and 
$\x_2\in{\mathbb T}_{\alpha_2}\cap{\mathbb T}_{\beta_2}$, 
we can think of $\x_1\times\x_2$ as an intersection point
${\mathbb T}_{\alpha_1\cup\alpha_2}\cap {\mathbb T}_{\beta_1\cup\beta_2}$.
Then, 
$$\spincrel_{w_1,z_2}(\x_1\times\x_2)=\spincrel_{w_1,z_1}(\x_1)\#\spincrel_{w_2,z_2}(\x_2).$$

The following is a straightforward generalization of the connected sum
principle for knot homology in the case of null-homologous knots,
c.f.~\cite{Knots}, \cite{RasmussenThesis}. A proof can be found, for
example, in Theorem~\ref{Knots:thm:ConnectedSumsOfKnots}
of~\cite{Knots}; we omit it here.

\begin{theorem}
\label{thm:Kunneth}
Fix $\xi_i\in\RelSpinC(Y_i,K_i)$ for $i=1,2$. 
There is a filtered chain homotopy equivalence 
$$\CFKinf(Y_1,K_1,\xi_1)\otimes_{\Z[U]} \CFKinf(Y_2,K_2,\xi_2)\
\longrightarrow \CFKinf(Y_1\#Y_2,K_1\#K_2,\xi_1\#\xi_2).$$
\end{theorem}

\begin{defn}
A {\em $U$-knot} is a knot in a three-manifold $Y$ with the property
that the induced filtration of $\CFinf(Y,\spinc)$ is trivial. More
precisely, given any $\xi\in\RelSpinC(Y,K)$, 
$\CFKinf(Y,K,\xi)$ is chain
homotopy equivalent via a relatively $\Z\oplus \Z$-filtered
chain homotopy to the chain complex $R$ which is a free, rank one
$\Z[U,U^{-1}]$ module with the trivial differential.
\end{defn}
 
The unknot in $S^3$ is a $U$-knot. Indeed, the results
of~\cite{GenusBounds} can be interpreted as saying that the only
$U$-knot in $S^3$ is the unknot. Other $U$-knots will be 
described in Lemma~\ref{lemma:LensSpaceUKnots}.

\begin{cor}
\label{cor:ConnSumUKnot}
If $K_2\subset Y_2$ is a $U$-knot, then for each
$\xi_1\in\RelSpinC(Y_1,\OrK_1)$ and $\spinc_2\in\SpinC(Y_2)$, 
there is some $\xi_2\in\RelSpinC(Y_2)$ representing $\spinc_2$, 
with the property that
$$\CFKinf(Y_1,\OrK_1,\xi_1)\cong 
\CFKinf(Y_1\# Y_2,\OrK_1\#\OrK_2,\xi_1\#\xi_2)$$
as
$\Z\oplus\Z$-filtered chain complexes.
\end{cor}

\begin{proof}
  Let $\xi\in\RelSpinC(Y_2)$ be any relative $\SpinC$ structure
  representing $\spinc_2$. Then, $\CFKinf(Y_2,\OrK_2,\xi)$ is
  quasi-isomorphic to the $\Z\oplus\Z$-filtered chain complex
  $\Z[U,U^{-1}]$ which contains a non-trivial element in filtration
  $(0,n)$ for some $n$. Letting 
  $$\xi_2=\xi+n\PD[\mu]$$ (where $\mu$ of
  course is a meridian for $K_2$), we see that
  $\CFKinf(Y_2,\OrK_2,\xi_2)$ is quasi-isomorphic to $\Z[U,U^{-1}]$
  with a generator in filtration level $(0,0)$ (c.f.
  Proposition~\ref{prop:Filtrations}). For this choice of $\xi_2$, the
  corollary is an immediate consequence of Theorem~\ref{thm:Kunneth}.
\end{proof}

\section{The Morse surgery formula}
\label{sec:Surgery}

We describe the Morse surgery formula, which expresses the Heegaard
Floer homologies of Morse surgery along a rationally null-homologous
knot in terms of its knot Floer homology.

Given $\spinc\in\SpinC(Y_\lambda(K))$, let
\begin{eqnarray*}
\BigAp_\spinc(Y,\OrK)&=&\bigoplus_{\{\xi\in\RelSpinC(Y_\lambda(K),\OrK)\big|
\Fill{Y_\lambda(K),\OrK}(\xi)=\spinc\}}\Ap_\xi(Y,\OrK) \\
\BigBp_\spinc(Y,\OrK)&=&\bigoplus_{\{\xi\in\RelSpinC(Y_\lambda(K),\OrK)\big|
\Fill{Y_\lambda(K),\OrK}(\xi)=\spinc\}}\Bp_\xi(Y,\OrK),
\end{eqnarray*}
where here $\Ap_\xi(Y,\OrK)$ and $\Bp_\xi(Y,\OrK)$ are defined as in
Equation~\eqref{eq:DefAp}.

In the above definition, we view $K$ also as a knot in $Y_\lambda(K)$.
An orientation on $K\subset Y$ naturally induces also an orientation
on the induced knot in $Y_\lambda(K)$ in a natural way: an orientation
on $K$ corresponds to an orientation on a meridian for $K$, which can
be thought of as supported in $Y-\nbd{K}$, which in turn can be
thought of as a subset of $Y_\lambda(K)$, where it in turn induces an
orientation of the induced knot in $Y_\lambda(K)$.

\begin{theorem}
\label{thm:SurgeryFormula}
Let $K\subset Y$ be a knot in a rational homology three-sphere, and let
$\lambda$ be a framing on $K$ with the property that $Y_\lambda(K)$
is also a rational homology three-sphere.
The Heegaard Floer homology $\HFp(Y_\lambda(K),\spinc)$ is calculated by 
the homology of the mapping cone of the chain map 
$$\Dp_{\spinc}\colon \BigAp_{\spinc}(Y,\OrK)\longrightarrow 
\BigBp_{\spinc}(Y,\OrK)$$
defined by 
$$\Dp_{\spinc}\left(\{a_\xi\}_{\xi\in\Fill{Y_\lambda(K),\OrK}^{-1}(\spinc)}\right)=
\{b_\xi\}_{\xi\in\Fill{Y_\lambda(K),\OrK}^{-1}(\spinc)},$$
where
$$b_{\xi}=\horp_{\xi-\PD[K_\lambda]}(a_{\xi-\PD[K_\lambda]})-\vertp_{\xi}(a_\xi).$$
\end{theorem}

The proof is modeled on the proof of the main result
from~\cite{IntSurg}.  Specifically, it is based on two key
ingredients: the relationship between the knot Floer homology of
$K\subset Y$ and the Heegaard Floer homology of three-manifolds
obtained as sufficiently large surgeries on $K$ from
Section~\ref{sec:Large}, and also an exact exact triangle which
relates, for all integers $m$, the Heegaard Floer homologies of
$Y_\lambda(K)$, $Y_{m\cm \mu+\lambda}(K)$, and $Y$ (the latter taken
with twisted coefficients).

We have not strived for maximal generality in
Theorem~\ref{thm:SurgeryFormula}.  As a technical tool, we make heavy
use of the rational grading on $\HFp$ which is defined on the Heegaard
Floer homology of any three-manifold equipped with a $\SpinC$
structure whose first Chern class is torsion,
cf. Section~\ref{HolDiskFour:sec:AbsGrade}
of~\cite{HolDiskFour}. Thus, Theorem~\ref{thm:SurgeryFormula} holds --
and indeed the proof we give below applies -- whenever we consider
$\SpinC$ structures over $Y_\lambda(K)$ whose first Chern class is
torsion.

\subsection{A surgery exact sequence}
\label{subsec:ExactSeq}

The following is essentially a restatement of
Theorem~\ref{IntSurg:thm:ExactSeq} of~\cite{IntSurg}:

\begin{theorem}
\label{thm:ExactSeq}
Let $Y$ be a closed, oriented three-manifold, and $K\subset Y$ be 
a knot with framing $\lambda$. Then, for all positive integers $m$,
there is a long exact sequence
\[
\begin{CD}
...@>>>\HFp(Y_\lambda(K)) @>{\Fp{1}}>>\HFp(Y_{m\cm\mu+\lambda}(K)) @>{\Fp{2}}>> \bigoplus^m \HFp(Y)
@>{\Fp{3}}>>...,
\end{CD}
\]
where here $\mu$ denotes the meridian of $K$.
\end{theorem}

We describe the maps in the above theorem below, after which
we make a 
few comments on the proof.

We place a basepoint $p$ on $\beta_g$, and consider twisted homology
with coefficients in $\Zmod{m}$; i.e. write
$\Z[\Zmod{m}]=\Z[T]/(T^m-1)$, and consider the chain complex
$\CFp(Y)\otimes_\Z \Z[\Zmod{m}]$ endowed with the differential
$${\underline \partial}^+ [\x,i] =
\sum_{\y\in\Ta\cap\Tb}\sum_{\{\phi\in\pi_2(\x,\y)\big|\Mas(\phi)=1
  \}}\#\left(\frac{\ModFlow(\phi)}{\R}\right) \cm T^{m_{p}(\phi)}\cm
[\y,i-n_{z}(\phi)]$$
where as usual here $\x\in\Ta\cap\Tb$, $i\geq 0$,
$\pi_2(\x,\y)$ denotes the space of homotopy classes of Whitney disks
connecting $\x$ and $\y$, $\Mas(\phi)$ denotes the Maslov index of
$\phi$, and terms in the above equation for which $i-n_z(\phi)<0$ are
to be dropped.  Moreover, $m_p(\phi)$ denotes the multiplicity of the
basepoint $p$ in the boundary of $\phi$; i.e. $p$ determines a
codimension one submanifold $\beta_1\times...\times \beta_{g-1}
\times\{p\}\subset \Tb$, and $m_p(\phi)$ denotes the intersection
number with the restriction of the boundary of $\phi$ with this
subset. We denote the complex by $\CFp(Y;\Z[\Zmod{m}])$. (In the
terminology of~\cite{HolDisk}, this is the chain complex for $Y$ with
twisted coefficients in $\Z[\Zmod{m}]$, where it is denoted
$\uCFp(Y;\Z[\Zmod{m}])$, however, as in~\cite{IntSurg}, we drop the
underline here in the interest of notational simplicity.) Recall that there
is an
isomorphism of chain complexes of modules over $\Z[\Zmod{m}]$,
\begin{equation}
\label{eq:DefChangeP}
\theta\colon \CFp(Y;\Z[\Zmod{m}])\stackrel{\cong}{\longrightarrow}
\CFp(Y)\otimes_\Z \Z[\Zmod{m}],
\end{equation}
where here the right-hand-side is
endowed with the differential which is the original differential on
$\CFp(Y)$ tensored with the identity map on $\Z[\Zmod{m}]$.  
There is
a corresponding identification $$\HFp(Y;\Z[\Zmod{m}])\cong
\HFp(Y)\otimes_\Z\Z[\Zmod{m}]\cong \bigoplus^m\HFp(Y).$$

The quantity $m_p(\phi)$ has an expression more in the spirit of
the constructions from Section~\ref{sec:Construction}.
Fix two basepoints, $w$ and $z$ on either side of $\beta_g$,
so that there is an arc in $\Sigma$ connecting them, but disjoint from
all the attaching circles except for $\beta_g$, which it meets transversally
in the single intersection point $p$. Then, 
\begin{equation}
\label{eq:DefMp}
m_p(\phi)=n_w(\phi)-n_z(\phi).
\end{equation}

We now define the chain maps $\fp_1$, $\fp_2$, and $\fp_3$ inducing
maps on homology $\Fp{1}$, $\Fp{2}$, and $\Fp{3}$ appearing in
Theorem~\ref{thm:SurgeryFormula}.  

The map $\fp_1$ is defined by
counting pseudo-holomorphic triangles between $\Ta$, $\Tc$, and $\Td$.
More precisely, note that the Heegaard triple
$(\Sigma,\alphas,\gammas,\deltas,z)$ determines a four-manifold
$X_{\alpha,\gamma,\delta}$ with three boundary components
\begin{equation}
\label{eq:DefX}
Y_{\alpha,\gamma}\cong Y_{\lambda}(K), \hskip1cm
Y_{\alpha,\delta}\cong Y_{m\cm\mu+\lambda}(K), \hskip1cm \text{and}
\hskip1cm Y_{\gamma,\delta}\cong \#^{g-1}(S^2\times S^1)\# L(m,1).
\end{equation}

\begin{defn}
  \label{def:DefCan}
Let $\spinccan \in\SpinC(\#^{g-1}(S^2\times S^1)\# L(m,1))$ denote the
{\em canonical $\SpinC$ structure}, i.e. this is the one which extends over
the tubular neighborhood $N$ of a sphere with self-intersection number $m$
(after attaching a two-handle and a collection of $g-1$ three-handles)
in such a way that its first Chern class evaluates as $\pm m$ on the
two-sphere generator of the two-dimenisonal homology of $N$.
\end{defn}
 
Let $\Theta_{\gamma\delta}$ denote the Floer homology class
corresponding to the generator (over $\Wedge^*
H_1(Y_{\gamma,\delta})\otimes \Z[U]$) of
$$\HFleq(Y_{\gamma,\delta},\spinccan) \cong \Wedge^*
H^1(Y_{\gamma,\delta})\otimes \Z[U]$$ in its canonical $\SpinC$
structure $\spinccan$.  (As usual, we  arrange for
the homology class
$\Theta_{\gamma\delta}$ to be represented by a single intersection
point in $\Tc\cap\Td$, which we also denote by
$\Theta_{\gamma\delta}$.)

We then define
\begin{equation}
\label{eq:DefF1}
f^+_1([\x,i]) = \sum_{\y\in\Ta\cap\Td}\sum_{\{\psi\in\pi_2(\x,\Theta_{\gamma\delta},\y)\big|\Mas(\psi)=0\}}
\#\ModFlow(\psi)\cm [\y,i-n_z(\psi)].
\end{equation}
Similarly, we define $\fp_2\colon \CFp(Y_{m\cm\mu+\lambda}(K))\longrightarrow \CFp(Y;\Z[\Zmod{m}])$ by
\begin{equation}
\label{eq:DefF2}
f^+_2([\y,i]) = \sum_{\w\in\Ta\cap\Tb}\sum_{\{\psi\in\pi_2(\y,\Theta_{\delta\beta},\w)
\big|\Mas(\psi)=0\}} \#\ModFlow(\psi)\cm [\w,i-n_z(\psi)]\cm T^{m_p(\psi)},
\end{equation}
where $m_p(\psi)$ is the natural extension of $m_p$ to triangles.
(In particular, in the case where we consider doubly-pointed Heegaard diagrams,
$m_p(\psi)=n_w(\psi)-n_z(\psi)$ as in Equation~\eqref{eq:DefMp}.)

To define $\fp_3$, we proceed as follows.
Fix
$\psi\in\pi_2(\Theta_{\gamma\beta},
\Theta_{\beta\delta},\Theta_{\gamma\delta})$.
The congruence class $c$ of $m_p(\psi)$ modulo $m$ is independent of the
choice of $\psi$. The map
$$\fp_3\colon \CFp(Y;\Z[\Zmod{m}])\longrightarrow \CFp(Y_{\lambda}(K))$$
is given by the formula
\begin{equation}
\label{eq:DefF3}
\fp_3(T^s\cm [\x,i])= 
\sum_{\y\in\Ta\cap\Tc} \sum_{\{\psi\in\pi_2(\x,\Theta_{\beta\gamma},\y)\big|
\stackrel{\Mas(\psi)=0,}{s+m_p(\psi})\equiv c\pmod{m}
\}}
\#\ModFlow(\psi)\cm [\y,i-n_z(\psi)].
\end{equation}

\vskip.2cm
\noindent{\bf{Proof of Theorem~\ref{thm:ExactSeq}.}}
Let $\Hp{1}$ and $\Hp{2}$ denote the null-homotopies of 
$\fp_2\circ\fp_1$ and $\fp_3\circ\fp_2$ respectively.
With this notation,
the proof of Theorem~\ref{IntSurg:thm:ExactSeq} applies to show that
the chain maps
\begin{eqnarray*}
\phi^+\colon \CFp(Y_\lambda(K))\longrightarrow \MCone(\fp_2)
&{\text{and}}&
\psi^+\colon \MCone(\fp_2) \longrightarrow \CFp(Y_{\lambda}(K))
\end{eqnarray*}
defined by 
\begin{equation}
\label{eq:DefPhi}
\phi^+(\xi)=(\fp_1(\xi),\Hp{1}(\xi)).
\end{equation}
and 
\begin{equation}
\label{eq:DefPsi}
\psi^+(x,y)=\Hp{2}(x)+\fp_3(y)
\end{equation}
respectively are quasi-isomorphisms, and hence establishing
Theorem~\ref{thm:ExactSeq} as stated above.
\qed
\vskip.2cm

\subsection{The analogue of Theorem~\ref{thm:SurgeryFormula}
  for the case of $\CFd$} Given an integer $\delta\geq 0$, let $\CFd$
denote the approximation of $\CFp$ described in
Subsection~\ref{subsec:CFd}. We state and prove analogue of
Theorem~\ref{thm:SurgeryFormula} for this group, following the pattern
of proof from~\cite{IntSurg}.

Given $\delta\geq 0$, let
\begin{eqnarray*}
\Ad_\xi(Y,\OrK)=C_\xi\{0\leq \max(i,j) \leq \delta\}&{\text{and}}&
\Bd_\xi(Y,\OrK)=\CFd(Y,\Fill{Y,\OrK}(\xi)).
\end{eqnarray*}
Given $\spinc\in\SpinC(Y_\lambda(K))$, let
\begin{eqnarray*}
\BigAd_\spinc(Y,\OrK)&=&\bigoplus_{\{\xi\in\RelSpinC(Y_\lambda(K),K)\big|
\Fill{Y_\lambda(K),\OrK}(\xi)=\spinc\}}\Ad_\xi(Y,\OrK) \\
\BigBd_\spinc(Y,\OrK)&=&\bigoplus_{\{\xi\in\RelSpinC(Y_\lambda(K),K)\big|
\Fill{Y_\lambda(K),\OrK}(\xi)=\spinc\}}\Bd_\xi(Y,\OrK).
\end{eqnarray*}
Consider the map
$$\Dd_{\spinc}\colon \BigAd_{\spinc}(Y,\OrK)\longrightarrow 
\BigBd_{\spinc}(Y,\OrK)$$
defined by 
$$\Dd_{\spinc}\left(\{a_\xi\}_{\xi\in\Fill{Y_\lambda(K),\OrK}^{-1}(\spinc)}\right)=
\{b_\xi\}_{\xi\in\Fill{Y_\lambda(K),\OrK}^{-1}(\spinc)},$$
where
$$b_{\xi}=\hord_{\xi-\PD[K_\lambda]}(a_{\xi-\PD[K_\lambda]})-\vertd_{\xi}(a_\xi).$$

Our aim in the present subsection is to prove the following analogue
of Theorem~\ref{thm:SurgeryFormula} for $\CFd$.

\begin{theorem}
  \label{thm:CFdSurgeryFormula}
  The group $\CFd(Y_\lambda(K),\spinc)$ is quasi-ismorphic to the
  mapping cone $\Xd_\spinc(\lambda)$ of
  $$\Dd_{\spinc}\colon \BigAd_\spinc(Y,\OrK)\longrightarrow \BigBd_\spinc(Y,\OrK).$$
\end{theorem}

We deduce Theorem~\ref{thm:CFdSurgeryFormula} from
a combination of
Theorems~\ref{thm:ExactSeq} and~\ref{thm:LargeNSurgeries}.
Theorem~\ref{thm:SurgeryFormula} will follow from
Theorem~\ref{thm:CFdSurgeryFormula}, together with some further
observations about gradings, as explained in
Subsection~\ref{subsec:Gradings}.

We turn our attention now towards proving Theorem~\ref{thm:CFdSurgeryFormula}.
We assume that both $d$ and $n$ have the same sign, and hence without
additional loss of generality that both are positive, returning to the
case where their signs are opposite in
Subsection~\ref{subsec:OppSigns}. 

It will be useful to have the following lemma.  In the statement,
``sufficiently large'' is meant with respect to the ordering on
$\SpinC(Y,K)$ with fixed filling $\spinct$ over $Y$ described earlier.
Specifically, recall that $\xi_1\leq \xi_2$ if $\xi_2=\xi_1 + j \cm
\PD[\mu]$ for some $j\geq 0$.
In particular, we say that a condition $P$ holds for all all
sufficiently large relative $\SpinC$ structures $\RelSpinC(Y,\OrK)$ if
there is a finite collection $\Xi\subset \RelSpinC(Y,\OrK)$ with the
properties that each $\spinct\in\SpinC(Y)$ can be represented by some
$\eta\in\RelSpinC(Y,\OrK)$ with $\eta\in\Xi$, and also
for any $\xi\in\RelSpinC(Y,\OrK)$ which has the property that
$\xi\geq \eta$ for some $\eta\in\Xi$ then $P$ holds for $\xi$.

\begin{lemma}
        \label{lemma:LargeBehaviour}
        If $\xi$ is sufficiently large, then
        $$\vertd_\xi\colon \Ad_\xi(Y,\OrK)
        \longrightarrow \Bd_\xi(Y,\OrK)$$
        is an isomorphism, while 
        $$\hord_\xi\colon \Ad_\xi(Y,\OrK)
        \longrightarrow \Bd_{\xi+\PD[K_\lambda]}(Y,\OrK)$$
        is trivial; moreover, if $\xi$ is sufficiently small, then 
        $\hord_\xi$ is an isomorphism, while $\vertd_\xi$ is trivial.
\end{lemma}

\begin{proof}
Since there are finitely many $\x\in\Ta\cap\Tb$,
there is a maximal $\xi_0$ among all $\spincrel_{w,z}(\x)$
with $\x\in\Ta\cap\Tb$. Clearly, for any $\xi\geq \xi_0$,
there are no generators $[\x,i,j]$ of $\CFKinf(Y,K,\xi)$
with $i<0$ and $j\geq 0$; thus, in fact, 
$\Ap_{\xi}(Y,\OrK)=C_\xi\{i\geq 0\}$,
i.e. $\vertp_\xi$ is the identity map, and the corresponding 
statement for  $\vertd_\xi$ follows at once.

The other assertions follow from similar reasoning.
\end{proof}

\vskip.3cm
\noindent{\bf Proof of Theorem~\ref{thm:CFdSurgeryFormula} when
$d$ and $n$ have the same sign.}  Choose $m=nk$, where $n$ is as in
Equation~\eqref{eq:OrderD}.  It is easy to see that $\PD[K_\lambda]$
has order $dk+1$ in $H^2(Y_{m\cm\mu+\lambda})$ whereas $\PD[K_\lambda]$
has order $nk$ in the quotient group $H^2(Y,Y-K;\Z)/m\PD[\mu]$.

Theorem~\ref{thm:ExactSeq} expresses $\CFd(Y_\lambda(K))$ as the
mapping cone of a map from $\CFd(Y_{m\cm \mu+\lambda}(K))$ to
$\CFd(Y,\Z[\Zmod{m}])$. 
We can think of $\CFd(Y,\Z[\Zmod{m}])$ more
invariantly, as a sum $$
\bigoplus_{[\xi]\in \RelSpinC(Y,K)/\Z
m\cm \mu}([\xi],\CFp(Y,\Fill{Y,\OrK}(\xi))),$$ where, as the notation
suggests, the index set consists of $m\cm \mu$-orbits in $\RelSpinC(Y,K)$
(each of which induces the same $\SpinC$ structure over $Y$, of
course). The identification is induced by the map $$T^i\otimes \xi
\mapsto (\relspinc_{w,z}(\x)+i\cdot \PD[\mu], \x).$$ 
Now,  the map $$\fd_3\colon
\CFd(Y_{m\mu+\lambda}(K))\longrightarrow \CFd(Y;\Z[\Zmod{m}])$$ can be written as
$$a\mapsto \sum_{\spinc\in\SpinC(W'_m(K))}
([\FillW{Y,\OrK}(\spinc)],\fd_{\spinc}(a)),$$
where here $W'_m(K)$ denotes the natural cobordism from $Y_{m\mu+\lambda}(K)$,
and $\FillW{Y,\OrK}$ is the map from $\SpinC$ structures over $W'_m(K)$ to relative
$\SpinC(Y,\OrK)$ as in Proposition~\ref{prop:FillW}.

Choosing $k$ sufficiently large, (where $m=nk$) according to
Theorem~\ref{thm:LargeNSurgeries}, we have that the summand
$\CFd(Y_{m\mu+\lambda},\spinct)\subset \CFd(Y_{m\mu+\lambda})$ is
identified with $\Ad_{\xi}(Y,\OrK)$ for $\xi=\Xi(\spinct)$.
In fact, there are only two homotopically non-trivial components of
$\fd_2|_{\CFd(Y_{m\mu+\lambda}(K),\spinct)}$, according
to Proposition~\ref{prop:OnlyTwo}: these are the
components belonging to $\spincx$ and
$\spincy$, whose corresponding maps are
identified with $\vertd_\xi$ and $\hord_\xi$ respectively. In view of
Equation~\eqref{eq:PDF},
$$\FillW{Y,\OrK}(\spincx+\PD[{\CapSurf}])=
\FillW{Y,\OrK}(\spincx)+m\mu+\lambda;$$
and hence, if $\Xi(\spinct)=\xi$ and
$\Xi(\spinct+\PD[K_\lambda])=\xi+\PD(K_\lambda)$, then the ranges of
$\vertd_\xi$ and $\horp_{\xi+\PD[K_\lambda]}$ coincide.

In effect, we have shown that $\CFd(Y_\lambda(K))$ is quasi-isomorphic
to the mapping cone of a map which splits according to $K_\lambda$
orbits in $\SpinC(Y,K)$ or equivalently, $\SpinC$ structures over $Y_\lambda(K)$.
In each such orbit, the map has the form
$$(\fd_2)' \colon \bigoplus_{s\in\Zmod{(dk+1)}} \Ad_{\xi + s\cm \PD[K_\lambda]}(Y,\OrK)
\longrightarrow \bigoplus_{s\in\Zmod{dk}} 
\Bd_{\xi+s\cm \PD[K_\lambda]}(Y,\OrK),$$
defined by
adding  
$$\vertd_s\colon \Ad_{\xi+s\cm \PD[K_\lambda]}(Y,\OrK) \longrightarrow
\Bd_{\xi+s\cm \PD[K_\lambda]}(Y,\OrK) $$
and
$$\hord_s\colon \Ad_{\xi+s\cm \PD[K_\lambda]}(Y,\OrK) \longrightarrow
\Bd_{\xi+(s+1)\cm \PD[K_{\lambda}]}(Y,\OrK). $$
According to Lemma~\ref{lemma:LargeBehaviour} and our hypothesis on
the sign of $d$, we see that for any $\xi$,
if $s$ is sufficiently large, then letting
$\xi'=\xi+s\cm \PD[K_\lambda]$, we have that $\hord_{\xi'}$
is null-homotopic, and $\vertd_{\xi'}$ is an isomorphism. It follows
from this (together with the analogous statement for $s$
sufficiently small) that the mapping cone of $(\fd_2)'$ is identified
with the mapping cone of $$ (\fd_2)''\colon \bigoplus_{s\in\Z}
\Ad_{\xi+(sn)\cm\PD[\mu]} (Y,\OrK)
\longrightarrow \bigoplus_{s\in\Z} \Bd_{\xi+s\cm \PD[K_\lambda]}(Y,\OrK),$$
gotten by adding 
all the maps $\vertd_s$ and $\hord_s$. This establishes the theorem.
\qed

The case where $d$ and $n$ have opposite signs is handled in
Subsection~\ref{subsec:OppSigns}.

\subsection{Gradings, and the proof of Theorem~\ref{thm:SurgeryFormula}
(in the case where $d$ and $n$ have the same sign)}
\label{subsec:Gradings}

In the proof of  Theorem~\ref{thm:CFdSurgeryFormula}, for each 
$\delta\geq 0$, we establish quasi-isomorphisms
$$\phi^\delta\colon \CFd(Y_\lambda(K),\spinc)
\longrightarrow \Xd_\spinc(\lambda).$$
We wish to conclude that
$$\HFp(Y_\lambda(K),\spinc)\cong H_*(\Xp_\spinc(\lambda)).$$
To this end, we must  establish that $\Xp_\spinc(\lambda)$ 
is a relative $\Z$-graded complex, which can be given an
absolute grading so that the map $\phi^\delta$
is homogeneous of degree zero. 

Clearly, the groups $\Ap_\xi(Y,\OrK)$ and $\Bp_\xi(Y,\OrK)$ are relatively
$\Z$-graded, and the maps
$\vertp_\xi$ and $\horp_\xi$ are all relatively $\Z$-graded maps. It
is now a formal consequence of the shape of 
$$\Dp_\spinc\colon \Ap_\spinc(Y,\OrK)\longrightarrow
\Bp_\spinc(Y,\OrK)$$
that $\Ap_\spinc(Y,\OrK)$ and
$\Bp_\spinc(Y,\OrK)$ can be given absolute $\Z$-gradings so that
$\Dp_\spinc$ drops this grading by one. This naturally endows
the mapping cone $\Xp_\spinc(\lambda)$ with an absolute $\Z$-grading.
In fact, this grading is uniquely determined up to an overall shift.

In the proof of Theorem~\ref{thm:ExactSeq}, we considered maps
$$\fd_1\colon \CFd(Y_\lambda(K)) \longrightarrow
\CFd(Y_{nk\cm\mu+\lambda}(K))$$ and 
$$\fd_2\colon \CFd(Y_{nk\cm \mu + \lambda}(K))
\longrightarrow
\CFd(Y;\Z[\Zmod{m}]).$$

Fix $\xi\in\RelSpinC(Y,K)$.  For $s\in\Z$, consider the projection
$$\Pi^A_s \colon \CFp(Y_{nk\cm \mu+\lambda}(K))
\longrightarrow\CFp(Y_{nk\cm
  \mu+\lambda}(K),\Fill{Y_{nk\cm\mu+\lambda}(K),\OrK}(\xi+s\cm
\PD[K_\lambda]).$$
We prove the following analogue of
Proposition~\ref{IntSurg:prop:RelativelyGraded} of~\cite{IntSurg}.

\begin{prop}
  \label{prop:DegreesOfMaps} Fix an absolute lift of the relative
  $\Z$-grading on $\CFp(Y_\lambda(K),\Fill{Y_\lambda(K),\OrK}(\xi))$, and an integer
  $\delta\geq 0$. Then, there is a constant $b$ so that, for all
  sufficiently large $k$, there are absolute lifts of the relative
  $\Z$-gradings on both 
  $$\Ad_{\xi+s\cm\PD[K_\lambda]}(Y,\OrK)\cong\CFd(Y_{nk\cm \mu +
  \lambda}(K),\spinct)
\subset \CFd(Y_{nk\cm
  \mu+\lambda}(K))$$ and $\CFd(Y,\Fill{Y,\OrK}(\xi+s\cm \PD[K_\lambda]))$ for all
  $|s|\leq b$, with the property that $\Pi^A_s\circ \fd_1$ and also
  the restriction of $\fd_2$ to 
  $$([s],\CFd(Y,\Fill{Y,\OrK}(\xi+s\cm
  \PD[K_\lambda]))\subset \CFd(Y,\Z[\Zmod{m}])$$ have degree zero.
\end{prop}

We return to a proof of Proposition~\ref{prop:DegreesOfMaps} after
introducing a lemma.

Recall (c.f. Equation~\eqref{eq:DefX}) that the map $\fp_1$ was
defined by counting holomorphic triangles for a Heegaard triple
describing a four-manifold which we denoted $X_{\alpha\gamma\delta}$,
whose boundaries consist of $Y_\lambda(K)$, $Y_{(nk)\mu+\lambda}(K)$
and $\#^{g-1}(S^2\times S^1)\# L(nk,1)$. We denote this four-manifold
here by $X(k)$, to call attention to its dependence on $k$.

The following lemma is an adaptation of the proof of
Lemma~\ref{IntSurg:lemma:GradingsOnX} from~\cite{IntSurg}.

\begin{lemma}
  \label{lemma:OnlyOne}
  Fix a constant $C_0$. Then, for all sufficiently large $k$, the
  following statment holds. Each $\SpinC$ structure over $Y_{nk\cm \mu
    + \lambda}(K)$ has at most one extension $\spinc$ over $X(k)$
  whose restriction to $\#^{g-1}(S^2\times S^1)\# L(nk,1)$ is the
  canonical $\SpinC$ structure (in the sense of
  Definition~\ref{def:DefCan}, and for which
  \begin{equation}
    \label{eq:MaximalSpinCInX}
    C_0\leq c_1(\spinc)^2 + nk.
  \end{equation}
\end{lemma}

\begin{proof}
  Observe that $H_2(X(k);\Z)\cong\Z$ is generated by a surface $\Sigma$
  with 
  $$\Sigma^2 = -nkd^2(dk+1).$$
  This can be seen, for example, by observing that the group of
  triply-periodic
  domains in $X(k)$ is generated by a relation of the form
  $$-(dk+1)d\cm \lambda + d\cm (nk\cm\mu+\lambda) + a + b,$$
  where $a$
  is a sum of curves among the $\alpha_i$ with $i=1,...,g$, and $b$ a
  sum of curves among the $\beta_j$ with $j=1,...,g-1$. The
  intersection number of the first two curves -- which gives the
  self-intersection number of the homology class corresponding to the
  triply-periodic domain -- is $-nkd^2(dk+1)$.
  
  Now, for any other $\SpinC$ structure over $X(k)$ which interpolates
  between the same two $\SpinC$ structures on $Y_{\lambda}(K)$ and
  $Y_{nk\cm\mu + \lambda}(K)$ (and whose restriction to the remaining
  boundary boundary component $\#^{g-1}(S^2\times S^1)\# L(nk,1)$ is
  the canonical $\SpinC$ structure) has the form $\spinc+j\PD[\Sigma]$,
  for some integer $j\neq 0$; thus,
  \begin{eqnarray*}
    c_1(\spinc+j \cm \PD[\Sigma])^2-c_1(\spinc)^2 &=&
    4(j^2 \Sigma\cm \Sigma + j \langle c_1(\spinc),[\Sigma]\rangle) \\
    &\leq& -4j^2nkd^2(dk+1)\left(1-\frac{|\alpha|}{j}\right) \\
    &\leq& -2 nkd^2(dk+1).
  \end{eqnarray*}
  Of course, if $k$ is sufficiently large, then 
  Inequality~\eqref{eq:MaximalSpinCInX} is violated.
\end{proof}

\vskip.2cm
\noindent{\bf{Proof of Proposition~\ref{prop:DegreesOfMaps}.}}
Fix an absolute lift of the relative $\Z$-grading on
$\CFp(Y_\lambda(K),\spinc)$, and fix some $\delta\geq 0$. Now, $\fd_1$
decomposes as a sum of homogeneous terms, indexed by those
$\spinc\in\SpinC(X(k))$ whose restriction to the boundary component
$\#^{g-1}(S^2\times S^1)\# L(nk,1)$ is the canonical $\SpinC$
structure. By a suitable adapation of
Equation~\eqref{eq:DimensionShift}, we see that each term is
homogeneous of degree $(c_1(\spinc)^2+nk)/4$ (in this application,
note that $X(k)$ has three, rather than two, boundary components, and
we are considering the pairing with a fixed homology class on the
third term).  It follows readily from Lemma~\ref{lemma:OnlyOne} that
there is at most one such $\SpinC$ structure which can induce a
non-trivial map from $\CFd(Y_\lambda(K),\Fill{Y_\lambda(K)}(\xi))$ to
$\CFd(Y_{nk\cm \mu+\lambda}(K))$.

Thus, we can grade $\CFd(Y_{nk\cm\mu+\lambda}(K))$
so that each component of $\fd_1$ has degree zero.

Similarly, we can endow $([s],\CFd(Y,\xi+s\cm
\PD[K_\lambda])\subset \CFd(Y,\Z[\Zmod{m}])$ with the grading for which
$$\vertd_s\colon \Ad_{\xi+s\cm \PD[K_\lambda]}(Y,\OrK)\cong
\CFd(Y_{nk\cm\mu+\lambda}(K),\xi+s\cm\PD[K_\lambda])\longrightarrow
([s], \CFd(Y,\xi+s \cm\PD[K_\lambda]))$$
has degree zero. We
must check that this is compatible with the grading for which
$$\hord_s\colon \Ad_{\xi+s\cm \PD[K_\lambda]}(Y,\OrK)\cong
\CFd(Y_{nk\cm\mu+\lambda}(K),\xi+s\cm\PD[K_\lambda])\longrightarrow
([s], \CFd(Y,\xi+(s+1) \cm\PD[K_\lambda]))$$
has degree zero.

To this end, it suffices to prove the following.  Fix $K\subset Y$
with framing $\lambda$. Fix also $\xi\in\RelSpinC(Y,K)$ and $t\in\Z$.
Then, for sufficiently large $\delta$, there is a homogeneous
element
\begin{equation}
\label{eq:RequiredClass}
a=\{a_s\}_{-b\leq s \leq b} \in \bigoplus_{-b\leq s \leq b} H_*(\Ad_{\xi+s\cm\PD[K_\lambda]}(Y,\OrK))
\end{equation}
with $a_t, a_{t+1}\neq 0$, and $H_*(\Dd_\spinc(a))=0$, i.e. where here
$H_*(\Dd_\spinc(a))$ denotes the map on homology induced by
$$\Dd_\spinc\colon
\bigoplus_{-b\leq s \leq b} \Ad_{\xi+s\cm \PD[K_\lambda]}(Y,\OrK)
\longrightarrow
\bigoplus_{-b+1\leq s\leq b} \Bd_{\xi+s\cm\PD[K_\lambda]}(Y,\OrK).$$

To see this, we proceed as follows. Abbreviate $\bigoplus_{-b\leq s
  \leq b} \Ad_{\xi+s\cm\PD[K_\lambda]}(Y,\OrK)$ by $\Ad$ (or $\Ad(b)$, when
we wish to refer to its dependence on $b$).  We have also $A^\infty$
and $\Ap$ which are obtained analogously.  Similarly, write $\Bd$ for
$\bigoplus_{-b+1\leq s\leq b}\Bd_{\xi+s\cm \PD[K_\lambda]}(Y,\OrK)$.
We have maps $\Dd\colon \Ad \longrightarrow \Bd$, $\Dp\colon
\Ap\longrightarrow \Bp$ and $\Dinf\colon \Ainf\longrightarrow \Binf$.
Note that $H_*(M(\Dinf\colon \Ainf\longrightarrow \Binf))\cong
\Z[U,U^{-1}]$.  Moreover, there is a map
$$\Pi^\infty\colon M(f^{\infty}\colon \Ainf\longrightarrow
\Binf)\longrightarrow \Binf$$
induced by $\Pi^\infty(\{a_s\}_{-b\leq
  s\leq b}) = \{\vertd_s(a_s)\}_{-b+1\leq s\leq b}$.
(Here, $f^{\infty}$ is the map obtained by adding the 
$v^\infty_s$ and $h^\infty_s$ to give a map from $A^\infty(b)$ to $B^\infty(b)$.) It is not
difficult to see  that for a generator
$$\alpha^\infty\in H_*(M(\Dinf\colon \Ainf\longrightarrow \Binf)),$$
each component of $\Pi^\infty(\alpha^\infty)$ is non-trivial. We can
also arrange for $a^\infty$ to be a homogeneous element. By
multiplying the generator through by $U^{-i}$ if necessary, we can
arrange for each component of $\Pi^\infty(\{a_s\}_{-b\leq s\leq
  b})$ to inject into $H_*(\Bp_{\xi+s\cm\PD[K_\lambda]}(Y,\OrK))$.
Thus, we obtain a homogeneous element $\alpha^+\in H_*(M(\Dp\colon
\Ap\longrightarrow \Bp))$ (the image of $\alpha^\infty$ under the natural map),
whose image $\Pi^+(\alpha^+)$ has non-trivial
components in $H_*(\Bp_{\xi+s\cm\PD[K_\lambda]}(Y,\OrK))$
(for all $-b+1\leq s\leq b$). Now,
for sufficiently large $\delta$,  $U^\delta
(a^+)=0$. Thus, we can find a homogeneous element $a^\delta\in
H_*(M(\fd\colon \Ad\longrightarrow\Bd))$ with the property that
$D^\delta (a^\delta)=0$, but $\Pi^\delta(a^\delta)$ has non-trivial
image in each component $H_*(\Bd_{\xi+s\cm\PD[K_\lambda]}(Y,
\OrK))$. 

To repeat: given $t$, we obtain a $\delta$ so that
$H_*(\MCone(f^\delta\colon \Ad(b) \longrightarrow \Bd(b)))$ has a
homology class of the required shape as in
Equation~\eqref{eq:RequiredClass}. Here $b$ is any cut-off, chosen
sufficiently large (note that $H_*(\MCone(f^\delta_{b}))$ stabilizes
when $b$ is sufficiently large).  Given this $\delta$, we then find an
appropriately large $k$ so that
$\CFd(Y_{nk\cm\mu+\lambda}(K),\spinct)$ is for each
$\spinct\in\SpinC(Y_{nk\cm\mu+\lambda}(K))$ is represented by
$\Ad_{\xi+s\cm\PD[K_\lambda]}(Y,\OrK)$ for some $s$
(Theorem~\ref{thm:LargeNSurgeries}). Thus, we have found a homology
class in $H_*(\MCone(f^\delta))$ of the required shape, establishing
the compatibility of the various gradings.

\qed
\vskip.2cm

With the above proposition in place, the proof of
Theorem~\ref{thm:SurgeryFormula} follows from
Theorem~\ref{thm:CFdSurgeryFormula} by following the pattern
from~\cite{IntSurg}.

Specifically, according to Proposition~\ref{prop:DegreesOfMaps}, 
an absolute grading on $\CFd(Y_\lambda(K))$ induces absolute gradings on both
\[
\bigoplus_{-\Width\leq s \leq \Width}
\Ad_{\xi+s\cm\PD[K_\lambda]}
\]
and
\[
\bigoplus_{
-\Width+n\leq s \leq \Width} 
([s], \CFd(Y,\Fill{Y,\OrK}(\xi+s\cm\PD[K_\lambda]))) \subset
\CFd(Y;\Z[\Zmod{nk}]),
\]
so that $\fd_1$ and $\fd_2$ both are graded maps with degree zero.
Let 
$$\Pi^B_s\colon \CFd(Y;\Z[\Zmod{nk}])\longrightarrow
 \CFd(Y,\Fill{Y,\OrK}(\xi+s\cm \PD[K_\lambda]))$$
denote the projection onto 
the summand 
$([s], \CFd(Y,\xi+s\cm\PD[K_\lambda]))$.

\begin{lemma}
  \label{lemma:RespectGradings}
  With respect to the above gradings, given $\delta\geq 0$, we have
  that for any sufficiently large $k$ and $|s|\leq \Width$, the map
  $\Pi^B_s\circ \Hd_1$ is homogeneous of degree $+1$.
\end{lemma}

\begin{proof}
This follows exactly the proof of Lemma~\ref{IntSurg:lemma:GradedHomotopy}
of~\cite{IntSurg}
\end{proof}

\noindent{\bf{Proof of Theorem~\ref{thm:SurgeryFormula}, in the case where
    $d$ and $n$ have the same sign.}}  This now follows exactly as
    in~\cite{IntSurg}. Specifically,
    Theorem~\ref{thm:CFdSurgeryFormula} provides, for any $\delta>0$,
    a quasi-isomorphism between $\CFd(Y_\lambda(K),\spinc)$ and
    $\MCone(\Dd_{\spinc})$. 
    
    There is an integer $\Width$ with the property that for all $s\geq
    \Width$, $\vertp_{s}$ and $\horp_{-s}$ are isomorphisms. Thus, we
    can truncate the mapping cone at this level to obtain a
    quasi-isomorphic complex.  According to
    Proposition~\ref{prop:DegreesOfMaps}, the truncated mapping cone
    inherits a grading, and according to
    Lemma~\ref{lemma:RespectGradings}, the quasi-isomorphism from
    Theorem~\ref{thm:CFdSurgeryFormula} respects these gradings. In
    effect, we have shown that for any $\delta\geq 0$, there is a
    graded isomorphism of $\HFd(Y_\lambda(K))$ with
    $\HFd(\MCone(\Dd_{\spinc}))$. It is now a formal consequence (see
    Lemma~\ref{IntSurg:lemma:HFdLemma} of~\cite{IntSurg} for details) that
    $\HFp(Y_\lambda(K))\cong\HFp(\MCone(\Dp_{\spinc}))$, as well.
    \qed

\subsection{The signs of $d$ and $n$ are opposite}
\label{subsec:OppSigns}

In the case where the signs of $d$ and $n$ in
Equation~\eqref{eq:OrderD} do not coincide, we write instead
\begin{equation}
\label{eq:OrderDMinus}
 d\cm \mu = -n \cm \mu \in H_1(Y-L;\Z)
\end{equation}
where in this new notation now, both $d$ and $n$ are positive.
Letting $m=nk$, we see that in this case,
$\PD[K_\lambda]$ has order $dk-1$ in $H^2(Y_{m\cm\mu+\lambda})$.

Theorem~\ref{thm:CFdSurgeryFormula} is obtained as before, with minor
notational changes. Specifically, applying Theorem~\ref{thm:ExactSeq}
as before, we have a description of $\HFd(Y_\lambda(K),\spinc)$ as the
homology of a mapping cone of of a map from $\CFd(Y_{m\cm
  \mu+\lambda}(K))$ to $\CFd(Y,\Z[\Zmod{m}])$, which decomposes as
$$(\fd_2)' \colon \bigoplus_{s\in\Zmod{(dk-1)}} \Ad_{\xi + s\cm \PD[K_\lambda]}(Y,\OrK)
\longrightarrow \bigoplus_{s\in\Zmod{dk}} 
\Bd_{\xi+s\cm \PD[K_\lambda]}(Y,\OrK),$$
defined by
adding  
$$\vertd_s\colon \Ad_{\xi+s\cm \PD[K_\lambda]}(Y,\OrK) \longrightarrow
\Bd_{\xi+s\cm \PD[K_\lambda]}(Y,\OrK) $$
and
$$\hord_s\colon \Ad_{\xi+s\cm \PD[K_\lambda]}(Y,\OrK) \longrightarrow
\Bd_{\xi+(s+1)\cm \PD[K_{\lambda}]}(Y,\OrK). $$
Combining Lemma~\ref{lemma:LargeBehaviour} with the
positivity of $d$ and $n$ in Equation~\eqref{eq:OrderDMinus},
we see that for any $\xi$, if $s$ is sufficiently large 
and $\xi'=\xi+s\cm \PD[K_\lambda]$, we have that
$\hord_{\xi'}$ is 
null-homotopic, and $\vertd_{\xi'}$ is an
isomorphism. It follows from this (together with the analogous
statement for $s$ sufficiently small) that the mapping cone of
$(\fd_2)'$ is identified with the mapping cone of $$
(\fd_2)''\colon
\bigoplus_{s\in\Z} \Ad_{\xi+(sn)\cm\PD[\mu]} (Y,\OrK) \longrightarrow
\bigoplus_{s\in\Z} \Bd_{\xi+s\cm \PD[K_\lambda]}(Y,\OrK),$$
gotten by
adding all the maps $\vertd_s$ and $\hord_s$. This establishes the
Theorem~\ref{thm:CFdSurgeryFormula} in the present case.

For fixed $\delta\geq 0$ and $k$ sufficiently large, $\fd_2$ remains a
homogeneous map as before. However, in the present case,
we study $\fp_3$, the map associated to the cobordism $W_\lambda$ from
$Y$ to $Y_\lambda$, in place of $\fp_1$.
To this end, we have the following:

\begin{lemma}
  \label{lemma:DistinguishThree}
  Fix an integer $\delta\geq 0$ and constant $C_0$. For all
  sufficiently large $k$, the following condition holds. For each
  $\spinc_0\in\SpinC(W_\lambda(K))$, there is at most one $\SpinC$
  structure of the form $\spinc=\spinc_0+ jk \PD[\CapSurf]$ with $j\in\Z$
  for which 
  $$C_0\leq c_1(\spinc)^2+nk. $$
  Here, $\CapSurf$ denotes a generator of $H^2(W_\lambda(K),Y)$.
\end{lemma}

\begin{proof}
  This follows from the fact that $W_\lambda$ is a negative-definite
  cobordism.
\end{proof}

The relevance of the lemma is the following: $\fd_3(\xi\otimes T^i)$
is a sum of the maps associated to $\SpinC$ structures $\spinc$ over
the cobordism $W_\lambda(K)$ which differ by addition of
$nk\PD[\mu]=dk\PD[\lambda]$, under the identification
$H^2(Y,K;\Z)\cong H^2(W_\lambda(K))$. By Equation~\eqref{eq:PDF},
the latter in turn correspond to $dk \Z\cdot\CapSurf$ orbits.

It follows from Lemma~\ref{lemma:DistinguishThree} that if $k$ is
sufficiently large, then given $\spinct\in\SpinC(Y)$ and
$i\in\Zmod{nk}$, there is a unique $\SpinC$ structure which
contributes non-trivially to $\fd_3(\xi\cdot \otimes T^i)$ for all
$\xi\in\HFd(Y,\spinct)$. Specifically, according to
Equation~\eqref{eq:DimensionShift}, the non-triviality of the map
places a lower bound (independent of $k$) on the square of the first
Chern class of any such $\SpinC$ structure
(c.f. Corollary~\ref{cor:LargeNFloerHomology}).

With these remarks in place, it now follows quickly that the map
$$\psi^\delta\colon \MCone(\fd_2) \longrightarrow \CFd(Y_\lambda(K))$$
given by $\psi^\delta(x,y)=H^\delta_2(x)+\fd_3(y)$ is a relatively
graded map.

The proof of Theorem~\ref{thm:SurgeryFormula} is then completed as before.
\section{The proof of Theorem~\ref{thm:RationalSurgeries}}
\label{sec:RatSurg}

\begin{figure}
\mbox{\vbox{\epsfbox{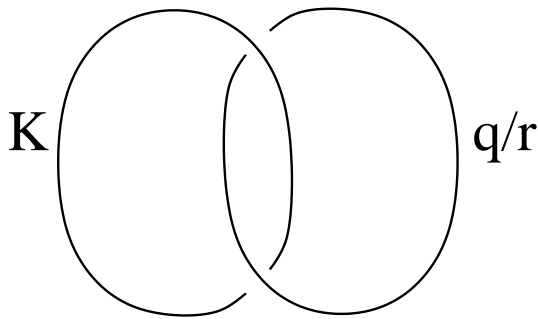}}}
\caption{\label{fig:Oqr}
  {\bf The knot $O_{q/r}$.}  Thinking of $K$ as a knot in the lens
  space obtained by performing $q/r$ on the other component of this
  link, we obtain the knot $O_{q/r}$.}
\end{figure}

If $K\subset Y$ is a null-homologous knot, then $Y_{p/q}(K)$ can be
realized by surgery with coefficient $a$ inside the knot $K\#
O_{q/r}\subset Y\#L(q,r)$, where $a$
is the greatest  integer smaller than or equal to $p/q$,
$a=\lfloor \frac{p}{q}\rfloor$, and
$$\frac{p}{q}=a+\frac{r}{q},$$
and $O_{q/r}\subset L(q,r)$ is the knot which is obtained by
viewing one component of the Hopf link as a knot
inside the lens space $L(q,r)$,
thought of as $q/r$ surgery on the other component of the Hopf link,
c.f. Figure~\ref{fig:Oqr}.

In view of these remarks, Theorem~\ref{thm:RationalSurgeries} is
proved by combining a model calculation of the knot Floer homology of
$O_{q/r}$, the K\"unneth principle for connected sums, together with
the surgery formula of Theorem~\ref{thm:SurgeryFormula}.

Of course, $L(q,r)-O_{q/r}$ is a solid torus, and consequently, 
there is an affine identification
$\RelSpinC(L(q,r),O_{q/r})\cong \Z$.

\begin{lemma}
\label{lemma:LensSpaceUKnots}
The knot $O_{q/r}\subset L(q,r)$ is a $U$-knot. Indeed,
there is an affine identification $\phi$ fitting into a commutative diagram
$$
\begin{CD}
\Z@>{\phi}>> \RelSpinC(L(q,r),O_{q/r}) \\
@VVV @VV{\Fill{L(q,r),O_{q/r}}}V \\
\Zmod{q} @>{\cong}>> \SpinC(L(q,r))
\end{CD}
$$
(where the left vertical arrow is the natural quotient map)
with the
property that
$$\HFKa(L(q,r),O_{q/r},\phi(i))\cong 
\left\{\begin{array}{ll}
        \Z      &       {\text{if $0\leq i \leq q-1$}} \\
        0       &       {\text{otherwise.}}
\end{array}
\right.
$$
\end{lemma}

\begin{proof}
Consider the standard genus one Heegaard decomposition of the lens
space $L(q,r)$, where $\alpha$ is a curve of slope $0/1$ and $\beta$
is a curve of slope $q/r$.  Placing two basepoints $w$ and $z$,
separated by an arc which is disjoint from $\alpha$ and meets $\beta$
exactly once, we obtain a doubly-pointed Heegaard diagram for
$O_{q/r}\subset L(q,r)$. Of course $\alpha\cap\beta$ consists of
exactly $q$ points $\{x_0,...,x_{q-1}\}$, each one representing a
different $\SpinC$ structure over $L(q,r)$. See Figure~\ref{fig:Uknot}
for an illustration.

\begin{figure}
\mbox{\vbox{\epsfbox{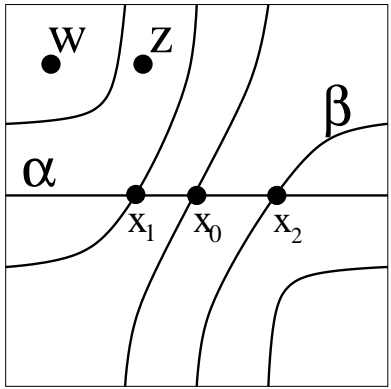}}}
\caption{\label{fig:Uknot}
{\bf Heegaard diagram for $O_{3/2}$.}
The intersection points $x_0$, $x_1$, $x_2$ are ordered so that
$\epsilon(x_0,x_1)=\epsilon(x_1,x_2)$ is a generator of the homology
of $L(3,2)-O_{3/2}$.}
\end{figure}

Now, 
\begin{equation}
\label{eq:Isomorphism}
H_1(L(q,r)-O_{q/r})\cong \frac{H_1(T^2-w-z)}{\Z \cdot \alpha + \Z
  \cdot \beta} \cong \frac{H_1(T^2)}{\Z \cdot \alpha} 
\cong \Z,
\end{equation}
which is generated by the homology class represented by a curve
$\gamma$ in $T^2$ with slope $1/0$. We specify an ordering on the
$\{x_i\}_{i=0}^{q-1}$ as follows. Let $m$ be the intersection point of
$\beta$ with the arc connecting $w$ and $z$. Our ordering is
specified by the conditions that the arc in $\beta$ from $x_i$ to $x_{i+1}$
is disjoint from all the other $x_k$ and $m$.
Clearly, with this ordering, 
$\epsilon(x_i,x_{i+1})=\gamma$. 
Let $\phi(i)$ correspond
to $\spincrel_{w,z}(x_i)$. The lemma follows.
\end{proof}

Suppose that $K$ is a null-homologous knot in an integral homology
three-sphere $Y$.
In this case, $H^2(Y,K)\cong \Z$, and also $H^2(Y\#L(q,r),K\#
O_{q/r})\cong \Z$.

\begin{lemma}
\label{lemma:ConnSumLensSpaceUKnots}
Suppose that $K\subset Y$ is a null-homologous knot in an integer homology
three-sphere.
Under the connected sum 
$$(Y,K)\# (L(q,r),O_{q/r}) \longrightarrow (Y\# L(q,r), K\# O_{q/r}),$$
the following diagram commutes:
$$
\begin{CD}
\Z\oplus \Z @>{f}>> \Z \\
@V{\cong}VV @VV{\cong}V \\
H^2(Y,K)\oplus H^2(L(q,r),O_{q/r})@>>> H^2(Y\# L(q,r), K\# O_{q/r}), \\
\end{CD}
$$
where $f(x,y)=qx+y$. 
Moreover, under this correspondence, if $K_\lambda$ is the push-off of $K$
with respect to the integral framing $a$, then $\PD[K_\lambda]$ represents
the element $p\in \Z\cong H^2(Y\# L(q,r),K\# O_{q/r})$,
where $p/q=a+q/r$.
\end{lemma}

\begin{proof}
Let $m$ and $\ell$ denote the meridian and the longitude
of $O_{q/r}$. We have an isomorphism
$H^2(L(q,r),O_{q/r})\cong \Z$
under which $m$ and $\ell$ are mapped to $r$ and $q$
respectively. Now in $Y\# L(q,r)-K\# O_{q/r}$,
clearly the meridian of $K$ is homologous to a
longitude for $O_{q/r}$, and hence it is mapped to $q$
under the isomorphism. 

Moreover, the push-off of $K_\lambda$ is homologous 
$m+a\ell$, which in turn is mapped to $r+a q=p$, under this above map.
\end{proof}

\vskip.2cm
\noindent{\bf{Proof of Theorem~\ref{thm:RationalSurgeries}.}}
Fix an identification $\RelSpinC(Y\# L(q,r),K\#O_{q/r})\cong \Z$, and
correspondingly think of $\Ap_s(Y\# L(q,r),K\#O_{q/r})$ as indexed by
$s\in\Z$.  Since $O_{q/r}$ is a $U$-knot
(c.f. Lemma~\ref{lemma:LensSpaceUKnots}), the K\"unneth principle
for connected sums, in the form of Corollary~\ref{cor:ConnSumUKnot}
applies to show that 
$$\Ap_s(Y\# L(q,r),K\#O_{q/r})\simeq \Ap_{g(s)}(Y)$$ for some $g\colon
\Z\longrightarrow \Z$
(where here $\simeq$ means filtered quasi-isomorphic).  
Indeed, according to Lemma~\ref{lemma:ConnSumLensSpaceUKnots},
the formula $g(s)$ satisfies $s=q g(s)+j$ where, according
to Lemma~\ref{lemma:LensSpaceUKnots}, $0\leq j < q$; i.e.
$$g(s)=\lfloor
\frac{s}{q} \rfloor.$$ Recall also that $\PD[K_\lambda]$ represents
$p$ times a generator of $H^2(Y\# L(q,r),K\# O_{q/r})$.

With these remarks in place, Theorem~\ref{thm:RationalSurgeries} is
obtained as a direct application of Theorem~\ref{thm:SurgeryFormula}.
\qed
\vskip.2cm

\subsection{The case of $\HFa$}

Some of the algebra is simpler when one considers $\HFp$,  rather
than $\HFa$. To this end, we let
$$\Da_{i,p/q}\colon \BigAa_i \longrightarrow \BigBa_i;$$
i.e.
$$\Da_{i,p/q} \{a_{\lfloor \frac{i+ps}{q}\rfloor}\}_{s\in\Z}
=\{b_{\lfloor \frac{i+ps}{q}\rfloor}\}_{s\in\Z},$$
where here
$$b_{\lfloor \frac{i+ps}{q}\rfloor}
=\verta_{\lfloor \frac{i+ps}{q}\rfloor}(a_{\lfloor \frac{i+ps}{q}\rfloor})
+\hora_{\lfloor \frac{i+p(s-1)}{q}\rfloor}(a_{\lfloor \frac{i+p(s-1)}{q}\rfloor}).$$
The proof of Theorem~\ref{thm:RationalSurgeries} adapts readily
to this context to give the following:

\begin{theorem}
\label{thm:RationalSurgeriesa}
Let $K\subset Y$ be a null-homologous knot.
There is a relatively graded isomorphism of groups
$$H_*(\Xa_{i,p/q})\cong \HFa(Y_{p/q}(K),i)$$
for each $i\in\Zmod{p}$.
\end{theorem}

\subsection{Absolute gradings}
\label{subsec:AbsGrade}

In fact, Theorem~\ref{thm:RationalSurgeries} (and indeed
Theorem~\ref{thm:SurgeryFormula}) actually determines
$\HFp(Y_{p/q}(K),i)$ as an absolutely graded group.

An absolute grading on $\Xp_{i,p/q}(K)$ compatible with the relative
grading is specified by fixing an absolute grading on $\Bp_{\lfloor
  \frac{i+ps}{q}\rfloor}$ (for any $s$), thought of as a subcomplex of
$\Xp_{i,p/q}(K)$ (note that it is independent of the choice of $K$).
This absolute grading in turn is in turn determined by a grading on
its homology $H_*(\Bp_{\lfloor \frac{i+ps}{q}
  \rfloor})\cong\InjMod{}$.  Finally, that datum is fixed by the
requirement $H_*(\Xp_{i,p/q}(O))\cong \InjMod{}$ in such a way that
its bottom-most non-trivial element is supported in dimension
$d(S^3_{p/q}(O),i)=d(L(p,q),i)$. (Explicit, recursive formulas for
these rational numbers can be found in
Proposition~\ref{AbsGraded:prop:dLens} of~\cite{AbsGraded}.)

This assertion follows easily from the proof of
Theorem~\ref{thm:RationalSurgeries}. Specifically, as in
Subsection~\ref{subsec:Gradings}, the composite of
$$\fd_1\colon \CFd(Y_\lambda(K)) \longrightarrow
\CFd(Y_{nk\cm\mu+\lambda}(K))$$
with the projection onto
$\CFd(Y_{nk\cm\mu+\lambda}(K),\spinct)\subset
\CFd(Y_{nk\cm\mu+\lambda}(K))$ is a homogeneous map (provided that $n$
is sufficiently large).  Moreover, its degree depends on the
intersection form of the cobordism $X(k)$ and the first Chern class of
$\spinct$.  Moreover, the proof of Theorem~\ref{thm:CFdSurgeryFormula}
shows that this map is also identified with the projection
$\Xd_{\spinc}(\lambda)\longrightarrow \Ad_\xi(Y,\OrK)$. Since the
intersection form of $X(k)$ does
not depend on the particular knot $\OrK$, the claim follows.

\section{Knots which admit $L$-space surgeries}
\label{sec:LSpaceKnots}

Suppose that $K\subset S^3$ is a knot which admits an $L$-space
surgery with positive slope $r$.  Examples of such knots include all
torus knots, and also knots from Berge's list, c.f.~\cite{Berge}.
Moreover, an alternating knot $K$ with unknotting number equal to one
gives rise to another knot $C\subset S^3$ which admits an $L$-space
surgery.  This new knot $C$ is obtained by performing the unknotting
operation, but connecting the two strands which were crossed by an arc
$\gamma$; $C$ then is the branched double-cover $\gamma$,
compare~\cite{UnknotOne}.

Theorem~\ref{thm:RatSurgeryLSpace} says that for such a knot, the
Alexander polynomial is determined from the correction terms. In this
application, we find it convenient to use the group $\HFp$; to do
this, recall the characterization of $L$-spaces in terms of this
group~\cite{NoteLens}:
a rational homology three-sphere
$Y$ is an $L$-space if and only if for each $\spinct\in\SpinC(Y)$,
$$\HFp(Y,\spinct)\cong \InjMod{}.$$ The correction terms $d(Y,\spinct)$ 
of an $L$-space is the minimal grading of any homogeneous
element of $\HFp(Y,\spinct)$.

The symmetrized Alexander polynomial $\Delta_K(T)$ plays a role, since
it is the Euler characteristic of the knot Floer homology; 
i.e.
$$\sum_s \chi(\HFKa_*(S^3,K,s)) \cdot T^s
= \Delta_K(T),$$ 
c.f. Equation~\eqref{Knots:eq:EulerChar} of~\cite{Knots}
or~\cite{RasmussenThesis}.

\vskip.2cm
\noindent{\bf{Proof of Theorem~\ref{thm:RatSurgeryLSpace}.}}
It is easy to see that both maps $\vertp_s\colon H_*(\Ap_s)
\longrightarrow H_*(\Bp)$ and $\horp_s\colon H_*(\Ap_s)\longrightarrow
H_*(\Bp)$ are surjective, since $\Bp\cong \HFp(S^3)$, and both maps are isomorphisms
in all sufficiently large degrees. 

Fix $i\in\Zmod{p}$, and let $s\in\Z$ be an arbitrary representative.
Split
$\BigAp_i=\Ap_s \oplus J$ 
where $J=\bigoplus_{t\neq s}
\Ap_{\lfloor \frac{i+sp}{q}\rfloor}$.
From the above remarks, it is clear that the restriction of
$\Dp_{i,p/q}$ to $H_*(J)$ surjects onto $H_*(\BigBp_i)$. It follows at
once that there is a surjection
$$\varphi \colon \Ker \left(H_*(\Dp_{i,p/q})\colon
  H_*(\BigAp_i)\longrightarrow H_*(\BigBp_i)\right) \longrightarrow
H_*(\Ap_s).$$
Moreover, it follows from this surjectivity, combined with
Theorem~\ref{thm:RationalSurgeries}, that
$$\Ker \left(H_*(\Dp_{i,p/q})\colon H_*(\BigAp_i)\longrightarrow H_*(\BigBp_i)\right)
\cong \HFp(S^3_{p/q}(K),i)\cong \InjMod{}.$$
The surjectivity of $\varphi$, together with the fact that
$H_*(\Ap_s)\cong H_*(\Bp)$ in all sufficiently large degree
combine to show that $H_*(\Ap_s)\otimes_\Z \Q \cong \InjMod{}\otimes_\Z\Q$.
Applying the same argument, only taking coefficients in $\Zmod{p}$ for arbitrary
$p$ shows that $H_*(\Ap_s)\cong\InjMod{}$.

Having established that $H_*(\Ap_s)\cong \InjMod{}$, and that
$\vertp_s$ and $\horp_s$ are isomorphisms in all sufficiently large
degrees, it follows that $\vertp_s$ is modeled on multiplication by
$U^{V_s}$, and $\horp_s$ is multiplication by $U^{H_s}$, where all
$V_s, H_s \geq 0$.

The condition that $\Ker \left(H_*(\Dp_{i,p/q})\colon
H_*(\BigAp_i)\longrightarrow H_*(\BigBp_i)\right)\cong\InjMod{}$
ensures readily that in each $i\in\Zmod{p}$, there can be at most one
integer $s\in\Z$ for which both $V_{\lfloor \frac{i+ps}{q}
\rfloor}$ and $H_{\lfloor \frac{i+ps}{q}\rfloor} \geq 0$.
Our aim now is to determine which value of $s$ has this property.

Let 
$m(s)=\min(V_s,H_s)$. We claim that 
\begin{eqnarray}
\label{eq:MonotoneT}
m(s_1)\leq
m(s_2)&{\text{if}}& s_1\leq s_2\leq 0~\text{or}~ s_1\geq s_2\geq 0. 
\end{eqnarray}

To this end, we have an exact sequence $$\begin{CD}
0@>>>C\{i<0~\text{and}~j\geq s\} @>>> C\{\max(i,j-s)\} @>{\vertp}>>
C\{i\geq 0\}@>>> 0.
\end{CD}$$
Since $\chi(\HFKa(K,s))=a_s$, it follows at once that for $s\geq 0$,
$$\chi(C\{i<0~\text{and}~j\geq s\})=t_s(K),$$
which in turn is the
same as $V_s$. Similarly, $\chi(C\{i\geq
0~\text{and}~j<s\})=t_s(K)+2s$, which is $H_s$. In particular, this
(together with a symmetric argument for $s\leq 0$) shows that
\begin{eqnarray}
\label{eq:AvsBIneq}
V_s \leq H_s~\text{for}~s\geq 0 &{\text{and}}& V_s \geq H_s~\text{for}~s\leq 0.
\end{eqnarray}
Incidentally, we have just established that
\begin{equation}
  \label{eq:MisT}
  t_s=m(s)\geq 0.
\end{equation}

Note that there is a natural quotient map $\Ap_{s}\longrightarrow
\Ap_{s+1}$, and indeed, the projection of $\Ap_s \longrightarrow
\Bp_s$ factors through this quotient. It follows at once that
$\{V_s\}_{s\in\Z}$ is a non-increasing sequence in $s$. Dually,
$\{H_s\}_{s\in\Z}$ is an non-decreasing sequence in $s$.  Combining
this with Equation~\eqref{eq:AvsBIneq},
Equation~\eqref{eq:MonotoneT} follows.

Having established that for each $i\in\Zmod{p}$, there is exactly one
non-zero integer among the $\{m(\lfloor\frac{i+ps}{q}\rfloor)\}_{s\in\Z}$,
Equation~\eqref{eq:MonotoneT} together with the fact that
$m(s)=m(-s)$ (an easy consequence of Equation~\eqref{eq:MisT}) shows
that $m(\lfloor\frac{i}{q}\rfloor)\neq 0$ implies that 
$$
\Big|\lfloor \frac{i}{q}\rfloor \Big|\leq 
\Big|\lfloor \frac{p}{2q} \rfloor \Big|
$$
for all $s\in\Z$.

Finally, note that the bottom-most generator of
$\Ap_{\lfloor\frac{i}{q}\rfloor}$ is an element whose degree is
$2m(\lfloor\frac{i}{q}\rfloor)$ less than the corresponding generator
for the unknot. In view of the remarks from
Subsection~\ref{subsec:AbsGrade}, the theorem now follows.  \qed
\vskip.2cm

\begin{lemma}
        \label{lemma:UKnotLemma}
        Let $K\subset S^3$ be a knot.
        The map
        $\vertp_s\colon \Ap_s\longrightarrow \Bp$
        is an isomorphism on homology for all $s\geq 0$ if and only if
        $\HFKa(K,s)=0$ for all $s>0$. 
\end{lemma}

\begin{proof}
        It is easy to see by descending induction on $s$ that
        $\vertp_s$ is an isomorphism on homology for all $s\geq d$ if
        and only if $\HFKa(K,s)=0$ for all $s>d$. 
\end{proof}

\vskip.2cm
\noindent{\bf{Proof of Corollary~\ref{cor:GordonConjecture}.}}
Suppose that $S^3_{p/q}(K)\cong S^3_{p/q}(O)$. Then, for some
permutation $\sigma\colon \Zmod{p}\longrightarrow \Zmod{p}$, we have
that $\HFp(S^3_{p/q}(O),i)\cong \HFp(S^3_{p/q}(K),\sigma(i))$ for all
$i$.  It follows from the proof of Theorem~\ref{thm:RatSurgeryLSpace}
(Equation~\eqref{eq:MisT})
that $t_{i}(K)\geq 0$ for all $i$. Summing
Equation~\ref{eq:DDifferences} over all $i$ and using the hypothesis
that $S^3_{p/q}(K)\cong S^3_{p/q}(O)$, we conclude that $t_i(K)\equiv
0$. It follows that $V_s=0$ for $s\geq 0$, i.e.
$\vertp_s\colon H_*(\Ap_s)\longrightarrow\HFp(S^3)$ is an isomorphism
for all $s\geq 0$.  According to Lemma~\ref{lemma:UKnotLemma},
$\HFKa(K,s)=0$ for all $s\neq 0$, and hence, according
to~\cite{GenusBounds}, $K$ is the unknot.  \qed \vskip.2cm

\section{On cosmetic surgeries}
\label{sec:CosmeticSurgeries}

Under favorable circumstances, the existence of an
orientation-preserving homeomorphism $\HFa(S^3_r(K))\cong
\HFa(S^3_s(K))$ for distinct $r$ and $s$ forces the knot Floer
homology of $K$ to agree with that for the unknot, and hence for the
knot to be unknotted.  This is not always the case, though. For
example, the knot $K=9_{44}$ has the property that
$\HFp(S^3_{+1}(K))\cong \HFp(S^3_{-1}(K))$, although
$S^3_{+1}(K)\not\cong S^3_{-1}(K)$.  Specifically, according to
Theorem~\ref{KT:thm:SmallKnots} of~\cite{calcKT}, the knot $9_{44}$
has
$$\HFKa(K,i)=\Z^{|a_i|}_{(i)}$$
where the subscript indicates the degree in which the summand is
supported, and $a_i$ is the $i^{th}$ coefficient of the Alexander 
polynomial
$$\Delta_K(T)=\sum_{i\in\Z} a_i\cm T^i = T^{-2}-4T^{-1}+7-4T +T^2.$$
It is now a straightforward application of the surgery formula,
either in the form given in the present paper or from~\cite{IntSurg}, that
$$\HFp(S^3_{+1}(K))\cong \InjMod{(0)}\oplus \Z_{(0)}^2\oplus
\Z_{(-1)}^2\cong \HFp(S^3_{-1}(K)).$$
On the other hand, Walter
Neumann~\cite{Neumann} informs us that the manifolds $S^3_{+1}(K)$ and
$S^3_{-1}(K)$ can be distinguished (using the computer program {\em
  Snap}~\cite{Snap}, see also~\cite{SnapPea}) by their hyperbolic volume; one
has hyperbolic volume roughly $5.52$, the other has hyperbolic volume
roughly $5.27$.

Our results here can be proved by restricting attention to Floer homology
with coefficients in any field $\Field$, which we suppress from the notation.
For concreteness, we restrict to $\Field=\Zmod{2}$, so that 
e.g. $\HFa(S^3)$ denotes the homology $H_*(\CFa(S^3)\otimes_{\Z}\Field)$.

\begin{defn}
Define $\nu(K)$ by
$$\nu(K)=\min\{s\in\Z\big| 
\verta_s\colon \Aa_s \longrightarrow \CFa(S^3)~\text{induces a non-trivial
map in homology}\}.$$
\end{defn}

\begin{lemma}
\label{lemma:Reflect}
If $K\subset S^3$ and $\Reflect{K}$ denotes its reflection, 
then either $\nu(K)$ or $\nu(\Reflect{K})$ is non-negative.
\end{lemma}

\begin{proof}
Let $\{\Filt_s\}_{s\in\Z}$ be the knot filtration of $\CFa(S^3)$
(i.e. in the notation of Section~\ref{sec:Construction}
$\Filt_s=C_{\xi_0}\{j<s\}$, where $\xi_0\in\RelSpinC(S^3,K)$
is the relative $\SpinC$ structure with trivial first Chern class)
and
let 
$$\tau(K)=\min\{s\in\Z\big| H_*(\Filt_s)\longrightarrow
\HFa(S^3)~\text{is non-trivial}\}.$$  
Recall that $\tau(K)=-\tau(\Reflect{K})$ and
also that $\nu(K)=\tau(K)$ or $\tau(K)+1$
(c.f. Lemma~\ref{4BallGenus:lemma:Reflection} of~\cite{4BallGenus} and 
Proposition~\ref{4BallGenus:prop:FourDInterp} from the same refence
respectively; or alternatively, see~\cite{RasmussenThesis}).
\end{proof}

The integer $\nu(K)$ dictates the maps on homology induced by all the
$\hora_s$ and $\verta_s$, according to the following two lemmas.

\begin{lemma}
\label{lemma:BehaviourOfNu}
For all $s\geq \nu(K)$, $\verta_s$ induces a non-trivial map in homology.
\end{lemma}

\begin{proof}
This follows at once from the fact that 
the image of $\verta_s$ is $\Filt_s$, and $\Filt_s \subseteq \Filt_t$
if $s\leq t$.
\end{proof}

\begin{lemma}
\label{lemma:BehaviourOfNuTwo}
If $\nu(K)\geq 0$, then for all $s>0$, $(\hora_s)_*=0$.
\end{lemma}

\begin{proof}
The image of $\hora_s\colon \Aa_s\longrightarrow \Ba$ is 
identified with $\Filt_{-s}\subset \Ba$, thus, if $s>0$ and $\tau(K)\geq 0$,
then the map on homology induced by $\hora_s$ is trivial.
\end{proof}

\begin{prop}
\label{prop:CalcRanks}
Let $K\subset S^3$ be a knot, and fix a pair of relatively prime
integers $p$ and $q$.
Then
\begin{equation}
\label{eq:CalcRank}
\Rk \HFa(S^3_{p/q}(K)) =
|p| + 2\max(0,(2\nu(K)-1)|q|-|p|) +|q| \left(\sum_{s} \left(\Rk H_*(\Aa_s)-1\right)\right) 
\end{equation}
\end{prop}

\begin{proof}
  Recall that if $K\subset S^3$, and
  $\Reflect{K}$ denotes its reflection, then $S^3_r(K)\cong
  -S^3_{-r}(\Reflect{K})$. Thus, in view of Lemma~\ref{lemma:Reflect}
  and Equation~\eqref{eq:RanksAreTheSame}, we can assume without loss
  of generality that $\nu(K)\geq 0$. Also, of course, we can assume
  that $q\geq 0$. Using Theorem~\ref{thm:RationalSurgeriesa} to
  express $\HFa(S^3_{p/q}(K))$ in terms of $H_*(\Xa_{p/q})$,
  calculating the rank of $H_*(\Xa_{p/q})$ is now straightforward
  application of Lemmas~\ref{lemma:BehaviourOfNu} and
  \ref{lemma:BehaviourOfNuTwo}.  In fact, the cokernel $\Da$ has rank
  $\max(0,(2\nu(K)-1)q-p)$, while its kernel has rank
  $$p+\max(0,(2\nu(K)-1)q-p)
  +q \left(\sum_{s} \left(\Rk H_*(\Aa_s)-1\right)\right).$$
\end{proof}

(The above proposition holds even in the case where $p$ and $q$ are not
relatively prime, with the understanding that if $(p,q)=a$, then
$\HFa(S^3_{p/q}(K))$ denotes the direct sum of $a$ many copies of
$\HFa(S^3_{p'/q'}(K))$, where here $p'=p/a$, $q'=q/a$.)

The following result is analogous to Lemma~\ref{lemma:UKnotLemma}.

\begin{prop}
        \label{prop:CharacterizeUnknot}
        Let $K$ be a knot with $\nu(K)=0$ and $\rk H_*(\Aa_s(K))=1$ for all     
        $s\in \Z$. Then, $K$ is the unknot.
        More generally, 
        $$g(K)=\max(\nu(K), \{s\in \Z\big| \rk H_*(\Aa_{s-1})>1\}).$$
\end{prop}

\begin{proof}
        Recall~\cite{GenusBounds} 
        that 
        $$g(K)=\max \{t\in\Z \big| \HFKa(K,t)\neq 0\}.$$ 
        This readily
        implies that the map $\verta\colon \Aa_s \longrightarrow
        \CFa(S^3)$ induces an isomorphism in homology for all $s\geq
        g(K)$: $\Aa_s$ has a subcomplex $C\{i<0,j=s\}$, the kernel of
        $\verta_s$, which is easily seen to be acyclic (as it is
        filtered by subcomplexes whose associated graded has homology
        isomorphic to $\bigoplus_{t>s}\HFKa(K,t)$), and moreover
        its image in $C\{i=0\}$ consists of the subcomplex
        $C\{i=0,j\leq s\}$, whose quotient is acyclic.

        Thus, we have proved that 
        $$g(K) \geq \max(\nu(K),\{s\in\Z\big| \rk H_*(\Aa_s)>1\}.$$
        Now consider 
        $$\begin{CD} 0@>>>
        C\{i<0,j\geq g-1\} @>>> C\{\max(i,j-g+1)\geq 0\} @>{\vertp_{g-1}}>>
        C\{i\geq 0\} @>>> 0, \end{CD}$$ where the midddle term here
        is, of course $\Ap_{g-1}$. Since
        $\HFK(K,t)=H_*(C\{(0,t)\})=0$ for all $t> g$, it follows
        from the natural filtration that $H_*(C\{i<0j\geq g-1\})\cong
        H_*(\{(-1,g-1)\})\cong H_*(C\{0,g\})=\HFKa(K,g)\neq 0$.  
        It follows at once that $\verta_{g-1}$ is not an isomorphism.
        It follows at once that either the map is not surjective, in which
        case $\nu(K)=g-1$, or it has kernel, in which case $\rk H_*(\Aa_{g-1})>1$.
\end{proof}

Before stating the next result, we recall that the total rank of $\HFa(Y)$
is independent of the orientation of $Y$:
\begin{equation}
\label{eq:RanksAreTheSame}
\rk \HFa(-Y)=\rk \HFa(Y).
\end{equation}
(Indeed, in Proposition~\ref{HolDiskTwo:prop:Duality} of~\cite{HolDiskTwo},
it is shown that $\HFa^*(-Y)\cong \HFa_*(Y)$.)

\begin{theorem}
        \label{thm:OppositeSigns}
        Let $K\subset S^3$ be a knot, and suppose that $r$ and $s$
        are distinct rational numbers with the property that
        $S^3_r(K)\cong \pm S^3_s(K)$. Then, either $S^3_r(K)$ is
        an $L$-space, or $r$ and $s$
        have opposite signs.
\end{theorem}

\begin{proof}
  Since $H_1(S^3_{\pm p/q}(K);\Z)\cong \Zmod{p}$, we can fix $p$
  throughout.  Now,
  according to Proposition~\ref{prop:CalcRanks}, for fixed $p$, and
  positive integral $q$, $\rk \HFa(S^3_{p/q}(K))$ is a monotone
  non-decreasing function of $q$. In fact, the function is
  strictly monotone except possibly for sufficiently small $q$, for
  which the rank is $p$.  But this ensures that $S^3_{p/q}(K)$ is an
  $L$-space.  The same remarks hold for the function $\rk
  \HFa(S^3_{-p/q}(K))$ for fixed $p$ and positive, integral $q$. Since
  the total rank of $\HFa(Y)$ is an invariant of the underlying
  (unoriented) three-manifold, 
  c.f.
  Equation~\eqref{eq:RanksAreTheSame}, the result holds.
\end{proof}

\begin{prop}
        \label{prop:NuEqualsZero}
        If $K$ is a non-trivial knot with $\nu(K)=0$, then 
        if there are rational numbers $r,s\in\Q$ with $r\neq s$ and
        $S^3_r(K)\cong S^3_s(K)$, then $r=\pm s$.
\end{prop}

\begin{proof}
  According to Proposition~\ref{prop:CalcRanks}, if $\rk
  \HFa(S^3_{p/q}(K))=\rk \HFa(S^3_{p/q'}(K))$ for $q'\neq \pm q$, then
  $\rk H_*(\Aa_s)=1$ for all $s$. Now, in view of
  Proposition~\ref{prop:CharacterizeUnknot}, $K$ is the unknot,
  contrary to our assumption. Thus, it follows that for a cosmetic
  surgery on $K$, $r=-s$.
\end{proof}

\begin{theorem}
\label{thm:GenusOne}
Let $K\subset S^3$ be a knot with Seifert genus equal to one. Then if
$S^3_{r}(K)\cong S^3_s(K)$ as oriented manifolds, then
either $S^3_r(K)$ is an $L$-space or $r=s$.
\end{theorem}

\begin{proof}
In view of Theorem~\ref{thm:OppositeSigns}, we can assume that $r>0$
and $s<0$.
  
As in the proof of Proposition~\ref{prop:CharacterizeUnknot}, it is
clear that
$$\verta_s\colon H_*(\Aa_s(K))\longrightarrow \HFa(S^3)$$
is an isomorphism for all $s>0$. Also, $\nu(K)\leq 1$.

We exclude the possibility that $\nu(K)=1$. 
For $S^3_r(K)$, we have that in any given $\SpinC$ structure,
$\HFa(S^3_r(K),\spinct)$ is described by
$$H_*(\Aa_0(K))^{m} \oplus \Field^{m-1}.$$
Since this group is also described as the Floer homology of $S^3_s(K)$
with negative $s$, it has the form
$$H_*(\Aa_0(K))^{n} \oplus \Field^{n+1}.$$
Since $\rk H_*(\Aa_0(K))\neq 0$ (since its Euler characteristic is $1$),
the equality of these two ranks forces at once that 
$m=n+1$ and $\rk H_*(\Aa_0(K))=1$.
But it is easy to see that a
relatively graded isomorphism
$$H_*(\Aa_0(K))^{n+1} \oplus \Field^{n}_{(-1)}
\cong H_*(\Aa_0(K))^{n}[1] \oplus \Field^{n+1}_{(0)}$$
cannot possibly hold. 

In the case where $\nu(K)=0$, we force a relatively graded isomorphism
$$K^n \oplus \Field_{(0)} \cong K[1]^m \oplus \Field_{(0)},$$ where here $K$
is the kernel of the map on homology $$(\verta_0\oplus
\hora_0)_*\colon H_*(\Aa_0)\longrightarrow \Field\oplus \Field.$$ Such a
relatively graded isomorphism can hold only if the rank of $K$ is
zero. In turn, this forces $\rk H_*(\Aa_0)=1$. From
Proposition~\ref{prop:CharacterizeUnknot}, it follows now that $K$ is
a trivial knot, contradicting our hypothesis.
\end{proof}

By Proposition~\ref{prop:CalcRanks} (or
alternatively~\cite{NoteLens}), a genus one knot with $L$-space
surgeries is easily seen to have the knot Floer homology groups of the
trefoil $T$. In particular, if $S^3_r(K)\cong S^3_s(K)$, then it
follows that $\HFp(S^3_r(T))\cong \HFp(S^3_r(T))$ as $\Q$-graded
groups. Of course, $\HFp(S^3_r(T))$ can be calculated explicitly (for
example using Theorem~\ref{thm:RatSurgeryLSpace}), and the authors
know of no pair of distinct rational numbers $r$ and $s$ for which
$\HFp(S^3_r(T))\cong \HFp(S^3_s(T))$. 

In a different direction, the algebra can be used in some cases to
exclude cosmetic surgeries with a fixed numerator $p$ (i.e. first
homology of the surgered manifold). We content ourselves here with a
discussion of the case where $p=3$.

\begin{theorem}
\label{thm:PEqualsThree}
Suppose that $K$ is a non-trivial knot.
Then, if $r,
s\in \Q$ both with numerators having absolute value $3$, and 
with $r\neq s$, we have that $S^3_r(K)\not \cong S^3_s(K)$ as
oriented manifolds.
\end{theorem}

\begin{proof}
By Lemma~\ref{lemma:Reflect}, we can arrange that $\nu(K)\geq 0$.
The possibility that $\nu(K)\geq 3$ is excluded by counting ranks.
According to
Proposition~\ref{prop:CalcRanks},
we have that
$$\HFa(S^3_{-3/q}(K))=3+
q \cm \left(2(2\nu(K)-1)+\sum_{s} \left(\Rk H_*(\Aa_s)-1\right)\right). 
$$
When $\nu(K)\geq 3$, then for all $q$,
we have that
\begin{equation}
\label{eq:QInequality}
6\leq (2\nu-1)q,
\end{equation}
and hence
$$\HF(S^3_{3/q}(K))=
-3+q\cm\left(2(2\nu(K)-1)+\sum_s\left(rk H_*(\Aa_s)-1\right)\right).
$$
It follows that
\begin{eqnarray*}
\rk\HFa(S^3_{3/q}(K))&<&
        \rk\HFa(S^3_{-3/q}(K)) \\
        &<&
        \rk\HFa(S^3_{3/{q+1}}(K)),
\end{eqnarray*}
the latter inequality following from the fact that
\begin{equation}
\label{eq:Monotonicity}
6<2(2\nu(K)-1)+\sum_s\left(rk H_*(\Aa_s)-1\right).
\end{equation}
It follows that three-manifolds $\{S^3_{3/q}(K)\}_{q\in\Z}$ are all distinct.

Next, we turn our attention to excluding $\nu(K)=2$.
In this case, the above argument
can fail if Inequality~\eqref{eq:Monotonicity} fails.
(Note that when $\nu(K)=2$
and $q=1$, Inequality~\eqref{eq:QInequality} fails, but it is
still the case that $\rk \HFa(S^3_3(K))<\rk \HFa(S^3_{-3}(K))$.)
Thus, it remains to exclude the 
possibility
$S^3_{-3/q}(K)\cong S^3_{3/{q+1}}(K)$.
Now, $\HFa(S^3_{-3/q}(K))\cong \HFa(S^3_{3/{q+1}}(K))$ forces
equality, rather than the inequality of Equation~\eqref{eq:Monotonicity}.
In particular, this forces $\rk H_*(\Aa_s)=1$ for all $s$ and, since
$\chi(\Aa_s)=1$, each $H_*(\Aa_s)$ is supported in even degree.
The condition that $\nu(K)=2$ and $\rk H_*(\Aa)=1$ forces both
$H_*(\Aa_0)$ and $H_*(\Aa_1)$ to be supported in negative degrees 
(c.f.~\cite{NoteLens}).
On the one hand, $H_*(S^3_{-3/(q+1)}(K),0)$ is non-trivial in degree zero,
while the even degree part of
$H_*(S^3_{3/{q}}(K),0)$ is carried 
by elements from $H_*(\Aa_0)$ and $H_*(\Aa_1)$,
which in turn is supported in negative degrees.
Observe that our use of absolute degree in this case is justified
by the fact that $L(3,2)\cong -L(3,1)$. 

We turn our attention now to the case where $\nu(K)=1$. Let
$$C=\sum_s\left(rk H_*(\Aa_s)-1\right).$$
Inequality~\eqref{eq:Monotonicity} holds (and hence leads to no
cosmetic surgeries, as above) except if $C=0$, $2$, or $4$. (Note that
$C$ is even.)  In the case where $C=2$, it follows from symmetry that
$\rk H_*(\Aa_0)=3$, and $\rk H_*(\Aa_s)=1$ for all $s\neq 0$. Thus, in
the cases where $C=0$ or $2$,

Proposition~\ref{prop:CharacterizeUnknot} ensures that $K$ is a knot
with Seifert genus equal to one.  According to
Theorem~\ref{thm:GenusOne}, $S^3_r(K)$ must be an $L$-space.
But if $K$ is a knot with Seifert genus one,
some positive surgery on $K$ gives an $L$-space, then
by considering Proposition~\ref{prop:CalcRanks}, we can conclude
$S^3_{3/q}(K)$ $q=1$ or $2$. But it the $\Q/\Z$-valued linking
form ensures that $S^3_{3/1}(K)\not\cong S^3_{3/2}(K)$.

When $\nu(K)=1$, we are left with the remaining case that $C=4$.  Once
again, since $\chi(\Aa_s)=1$, there are two ways in which the total
rank $C$ can in principle distribute over the $H_*(\Aa_s)$: either
$H_*(\Aa_0)$ has rank $5$ (and $H_*(\Aa_s)$ has rank one for all other
$s$) or there is exactly one positive $s$ with $\rk
H_*(\Aa_s)=H_*(\Aa_{-s})=3$ and for all $t$ with $|t|\neq |s|$, $\rk
H_*(\Aa_t)=1$.  In the first case,
Proposition~\ref{prop:CharacterizeUnknot} again ensures that $g(K)=1$,
and hence we can apply Theorem~\ref{thm:GenusOne}.

We exclude $S^3_{-3/q}(K)\cong S^3_{3/{q+1}}(K)$ as follows. 
It is easy to see that in this case, in the model for
the spin structure has the form 
$$\HFa(S^3_{-3/q}(K),0)\cong H_*(\Aa_0)^{n}[1] \oplus \Field^{n+1}_{(0)}\oplus G'$$ 
while $$\HFa(S^3_{3/{q+1}}(K),0)\cong H_*(\Aa_0)^{m+1}\oplus
\Field^{m}_{(-1)}\oplus G.$$ Here, $G$ consists of some even number of copies
of $M$,  the kernel of
$$(\verta)_*\colon H_*(\Aa_s) \longrightarrow \Z,$$
each of which is shifted up some even amount;
while $G'$ consists
of direct sum of some even number copies of $M$, each of which is
shifted by an odd degree. 
Except in the case where $n=0$,
which coincides with the case where $q=2$, which we need not consider,
these shifts in degree for $G'$ are downward.

Now, when $q$ is even, then
\begin{eqnarray*}
n=2\left(\lfloor \frac{q-4}{3}\rfloor+1\right)&{\text{and}}&
m=n-1,
\end{eqnarray*}
while if $q$ is odd then 
\begin{eqnarray*}
n=2\left(\lfloor \frac{q}{3}\rfloor+1\right) &{\text{and}}&
m=n+1.
\end{eqnarray*}
In either case, we have that $n\not\equiv m\pmod{2}$.
It follows from the assumption that
$\HFa(S^3_{3/q}(K),0)\cong \HFa(S^3_{-3/{q+1}}(K),0)$ readily that
$H_*(\Aa_0)\cong \Field_{(-2)}$. Thus, we have that
$$\Field^{n}_{(-1)}\oplus \Field^{n+1}_{(0)}\oplus G' \cong
\Field^{m}_{(-1)} \oplus \Field^{m+1}_{(-2)} \oplus G.$$
Assuming that
$n\geq 0$, we have the maximal degree of any element of $G'$ is
smaller than the maximal degree of any element of $G$. It follows at
once that $H_*(\Aa_s)$ consists of $\Z$ in degree $0$, and also the
minimal degree of any element of $G'$ is $-1$. In fact, it follows
that the rank of the degree zero part of $G$ is ${n+1}$, while 
the rank of the degree
$-1$ part of $G$ is $m+1$.  However, since $G$ and $G'$ both have even
rank in each degree, this contradicts the above observation that $n$
and $m$ have different parity.

Finally, we turn attention to the case where $\nu(K)=0$.  Since $K$ is
a non-trivial knot, Proposition~\ref{prop:NuEqualsZero} ensures that
the only possibility is that $S^3_{3/q}(K)\cong S^3_{-3/q}(K)$. But
this is easily excluded by the $\Q/\Z$-valued linking form.
\end{proof}

\subsection{Proofs of theorems in Subsection~\ref{subsec:CosmeticSurgeries}.}

Theorem~\ref{intro:GenusOne} is a restatement of Theorem~\ref{thm:GenusOne}
proved above. Theorem~\ref{intro:OppositeSigns} is stated and proved
as Theorem~\ref{thm:OppositeSigns} above. 
Theorem~\ref{intro:PEqualsThree} is
a restatement of Theorem~\ref{thm:PEqualsThree} stated and proved above.

\section{Seifert fibered spaces}
\label{sec:Seiferts}

Methods from this paper lead to the calculation of the Heegaard Floer
homology groups of a large class of Seifert fibered spaces. The
primary ingredients here are the calculation of the knot Floer
homology of the ``Borromean knot'' (c.f.
Section~\ref{Knots:sec:FiberedExamples} of~\cite{Knots}), the knot
$O_{q/r}$ considered in Section~\ref{sec:RatSurg}, and
Theorem~\ref{thm:SurgeryFormula}. Our aim here is to state and prove
these results.

Let $h\colon \Z \longrightarrow \Z$ be a function with the property
that
$$\lim_{s\goesto \pm \infty} h(s) = +\infty.$$  We describe here a
natural $\Z[U]$-module associated to $X$. It is interesting
to compare the following construction with a construction of
N{\'e}methi~\cite{Nemethi}, see also~\cite{SomePlumbs}.

A {\em well at height $n\geq 0$} is a pair of integers $(i,j)$ with
$i<j-2$ and the property that $h(k)\leq n$ for all $i\leq k \leq j$,
while $h(i)>n$ and $h(j)>n$.  Let $M_n(h)$ denote the free Abelian
group generated by $W_n(h)$.

If $x\in W_n(h)$ and $y\in W_{n-1}(h)$, we write $x>y$ if $x=(i,j)$
and $y=(i',j')$ with $i \leq i'<j'\leq j$.  Define $U\colon W_n(h)
\longrightarrow W_{n-1}(h)$ by the formula
$$U\cm x = \sum_{\{y\in W_{n-1}(h)\big| x>y\}} y,$$
and let 
$$\FormHFp(h)=\bigoplus_{n\in \Z} W_n(h)$$
be the induced module over $\Z[U]$. Indeed, we can view this as
a graded $\Z[U]$ module by the grading which sends $W_n(h)$ to $2n$.

In the language of N{\'e}methi, the set of wells forms a root, and
$\FormHFp(h)$ is the associated $\Z[U]$-module.

Let $Y$ be a Seifert fibered space over a genus $g$ orbifold with
Seifert invariants $(a,r_1/q_1,...,r_n/q_n)$ over a genus $g$ base,
c.f.~\cite{Seifert}, \cite{Scott}, \cite{FurutaSteer}.  The {\em
  orbifold degree} is the quantity
$$\deg(Y)=a+\sum_i \frac{r_i}{q_i}.$$
Recall that $b_1(Y)$ is even if
and only if $\deg(Y)\neq 0$. By reversing the orientation on $Y$, we
can arrange for $\deg(Y)>0$.

There is a presentation of the first homology of $Y$ as
\begin{eqnarray}
\label{eq:PresentHomology}
H_1(Y;\Z)
&{\cong}&
\frac{H_1(\Sigma;\Z)\oplus \Z m_0 \oplus \Z m_1 \oplus ... \oplus \Z m_n}
{\left(\begin{array}{ll}
a\cm  m_0 + \sum_{i=1}^n r_i\cm m_i  = 0, & \\
 r_i \cm m_0 - q_i \cm m_i = 0, & i=1,...,n
\end{array}\right)}
\end{eqnarray}
(c.f. below for more specifics).

\begin{theorem}
        \label{thm:SeifertFiberedSpaces}
        Let $Y$ be a Seifert fibered space over a genus $g$ orbifold
        with positive degree, and Seifert invariants
        $(a,r_1/q_1,...,r_n/q_n)$. There is an affine identification
        $\SpinC(Y)\cong H_1(Y;\Z)$ with the following properties.
        \begin{itemize}
        \item The $\HFp(Y,\spinc)$ is non-trivial only to
        those $\SpinC$ structures which are supported in the span of
        $m_0$, ..., $m_n$ (in the notation of
        Equation~\eqref{eq:PresentHomology}). 
        \item Let $\spinc$ be some $\SpinC$ structure
          over $Y$, and let 
        $\xi_0\cm m_0 + \xi_1\cm m_1 + ... + \xi_n \cm m_n$
        be a representative with
        with $0\leq \xi_i \leq q_i$.
        For integers $-g\leq t \leq g$, let $\delta_t\colon \Z
        \longrightarrow \Z$ be the function defined by
        $$\delta_t(s) = (-1)^{s+1} t + \left(\xi_0 + a \cm s + \sum_{i=1}^n
        \lfloor \frac{\xi_i + r_i \cm s}{q_i}\rfloor\right), $$ 
        and $h_t\colon \Z \longrightarrow \Z$
        be the function 
        $$h_t(s) = 
        \left\{\begin{array}{ll}
                \sum_{i=0}^{s-1} \delta_t(i) & {\text{if $s\geq 0$}} \\
                -\sum_{i=s}^{-1}\delta_t(i) & {\text{if $s<0$.}}
                \end{array}\right.
        $$
        Then, there is a relatively $\Z$-graded isomorphism of $\Z[U]$-modules:
        \begin{equation}
                \label{eq:SeifertIsomorphism}
                \HFp(Y,\spinc) \cong \bigoplus_{-g\leq t\leq g} \Wedge^{g+t} H_1(\Sigma;\Z)
                \otimes_{\Z} \FormHFp(h_t).
        \end{equation}
        \end{itemize}
\end{theorem}

In the above statement, the tensor product is to be taken in the
graded sense: $\Wedge^{g+t} H_1(\Sigma;\Z)$ is supported in grading
$t$.  Note also that the right-hand-side is graded only in the
relative sense, corresponding to our choice of $\spinc$.

Note that in the case where $g=0$, this recaptures (as a relatively
$\Z$-graded group) the description of $\HFp$ for rational homology
Seifert fibered spaces given by N{\'e}methi in terms of his
``computational sequences''.  We prove
Theorem~\ref{thm:SeifertFiberedSpaces} in
Subsection~\ref{subsec:ProveSeifs} below, after giving some sample
calculations.

Although we have described here the Floer homology only as a
relatively graded group, the absolute grading can be obtained by
comparing the summand corresponding to $t=0$ with the calculation of
genus zero Seifert fibered spaces from~\cite{SomePlumbs}, see
also~\cite{Nemethi}.

\subsection{Sample calculations}

We begin with some generalities.  Suppose $h\colon \Z\longrightarrow
\Z$ is a function with $\lim_{s\goesto \pm\infty} h(s)=+\infty$. Let
$\delta(s)=h(s)-h(s-1)$.  Clearly, the rank of $\Ker U\subset
\FormHFp(h)$ agrees with the number of pairs of integers $(i,j)$ with
$i<j$, $\delta(i)<0$, $\delta(j)>0$ and $\delta(k)=0$ for all $i<k<j$
(i.e. these are the local minima of $h$).

In the case of Theorem~\ref{thm:SeifertFiberedSpaces}, explicitly
finding these local minima is a straightforward matter: in the
statement of the theorem, $\delta_t(s)$ differs from a linear function
of $s$ by a periodic function whose period is  the least common multiple
of $2$ (when $t\neq 0$) and the integers $\{q_i\}_{i=1}^n$.

We use Theorem~\ref{thm:SeifertFiberedSpaces} to calculate the Seifert
fibered space over a genus one base, and Seifert invariants $(-1,
\frac{1}{2}, \frac{2}{3})$. Note that $H_1(Y;\Z) \cong \Z^2$,
and there is a unique $\SpinC$ structure with non-trivial $\HFp$.

To calculate it, we proceed as follows. 
Observe that 
$$ \delta_t(s+6)=1 + \delta_t(s).$$
Moreover, the sequence $\{\delta_0(s)\}_{s=0}^{5}$ is
\[
         \{ 0,-1,0,0,0,0\}
\]

Set $t=0$. The sequence of integers $\{h_0(i)\}$ clearly has a unique
local minimum. It follows that the corresponding summand of $\HFp(Y)$
is isomorphic to 
$$H_1(T^2)\otimes_\Z \InjMod{}\cong (\InjMod{})^2.$$ 
In our subsequence descriptions, we will fix an absolute
grading lifting the relative grading, with the additional convention
that this corresponding summand of $\HFp(Y)$ is isomorphic to
$(\InjMod{(0)})^2$.  Indeed, this convention corresponds to the naturally
induced absolute grading on $\HFp(Y)$, as can be seen by comparing
against the case where $g=0$, and observing that the Seifert fibered space
with corresponding invariants is $S^3$.

\begin{figure}
\mbox{\vbox{\epsfbox{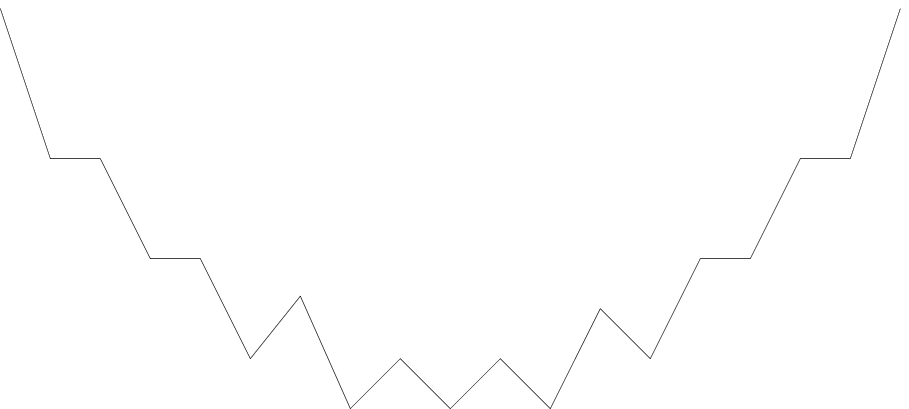}}}
\caption{\label{fig:SeifertEx}
{\bf Height function for the Seifert fibered space with invariants (-1,1/2,2/3)
at $t=-1$.}}
\end{figure}

When $t=-1$, all of the minima of the
sequence $\{h_{-1}(i)\}_{i\in\Z}$
clearly occur for $-5\leq i \leq 12$, where 
$\{\delta_{-1}(s)\}_{s=-5}^{12}$ takes the form
\[ \{ -3,0,-2,0,-2,1,-2,1,-1,1,-1,2,-1,2,0,2,0,3\} .\]
We plot the corresponding function 
$\{h_{-1}(s)\}_{s=-5}^{12}$ in
Figure~\ref{fig:SeifertEx}. It is clear from this description that the
corresponding summand of $\HFp(Y)$ takes the form 
$$\InjMod{(-1)}\oplus
\Z^2_{(-1)} \oplus \Z^2_{(1)},$$ under the normalization convention established
in the previous previous paragraph.

In the case where $t=1$, we have that
$\{\delta_{1}(i)\}_{i=-12}^5$ takes the form
$$\{ -1,-2,0,-2,0,-1,0,-1,1,-1,1,0,1,0,2,0,2,1\},$$
which clearly contains all the minima of $\{h_{-1}(i)\}$ for $i\in\Z$.
It follows that the corresponding summand of $\HFp(Y)$ is of the form
$\InjMod{(-1)}\oplus\Z_{(-1)}$.

Putting this together, if $Y$ denotes
the Seifert fibered space over a genus one base with Seifert invariants
$(-1,1/2,2/3)$, then its Heegaard Floer homology 
can be described as a graded $\Z[U]$ module by
$$\HFp(Y)\cong \Z_{(1)}^2 \oplus \Z_{(-1)}^3 \oplus(\InjMod{(0)})^2
\oplus (\InjMod{(-1)})^2.$$

\subsection{Proof of Theorem~\ref{thm:SeifertFiberedSpaces}.}
\label{subsec:ProveSeifs}

It is useful to have a mild generalization of the rational surgeries
formula.

\begin{defn}
Let $K\subset Y$ be a null-homologous knot in a
three-manifold $Y$.  Fix an integer $a$ and an $n$-tuple of rational
numbers $\{\frac{q_i}{r_i}\}_{i=1}^n$.  Consider an $n$-tuple of
unknotted circles $O_{i}$ each of which links $K$ once, and which are
pairwise mutually unlinked.  Let $Y(K,a,\{q_i/r_i\}_{i=1}^n)$ denote
the three-manifold obtained as $a$-surgery on $K$, followed by
$-\frac{q_i}{r_i}$ surgery on each $O_i$.  This three-manifold
$Y(K,a,\{q_i/r_i\}_{i=1}^n)$ is said to be obtained as a {\em generalized
rational surgery} on $K\subset Y$, with Seifert invariants
$(a,\{q_i/r_i\}_{i=1}^n)$.
\end{defn}

For fixed $K\subset Y$,
let
\begin{eqnarray*}
Y'=Y\#\left(\#_{i=1}^n L(q_i,r_i)\right) &{\text{and}}& 
K'=K\#\left(\#_{i=1}^n O_{q_i/r_i}\right)
\end{eqnarray*}
in the notation of Section~\ref{sec:RatSurg}.
Of course, $Y(K,a,\{q_i/r_i\}_{i=1}^n)$ can be thought of as
the three-manifold gotten by $a$-surgery on $K'\subset Y'$.

\begin{figure}
\mbox{\vbox{\epsfbox{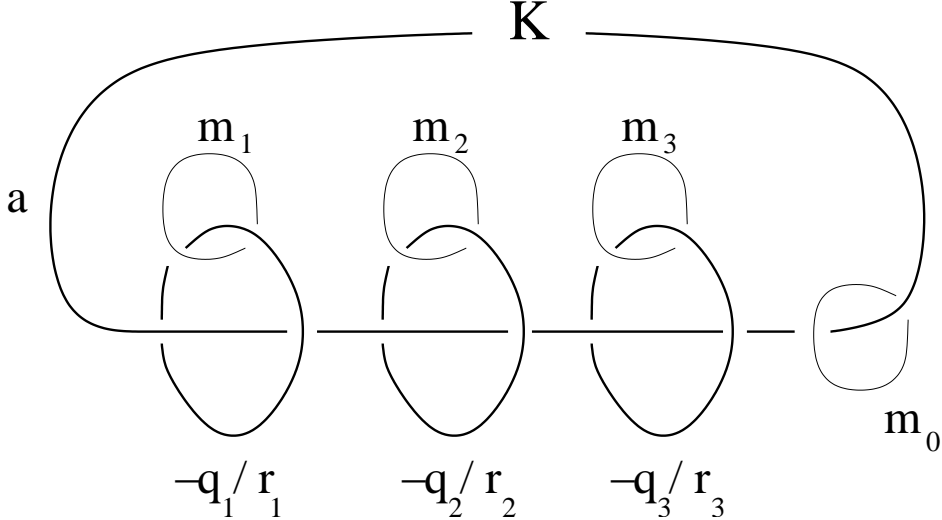}}}
\caption{\label{fig:Seifert}
{\bf A schematic illustration of generalized rational surgery.}
We take here $n=3$. $K$ represents some initial knot, and 
$\{q_i/r_i\}_{i=1}^3$ represent surgery instructions on the
unknots, while $K$ is framed with framing $a$. The lightly drawn
circles represent generators of the homology of the complement of
the dark link (they are meridians).}
\end{figure}

As an example, if we start with the unknot $K\subset S^3$ and form the
three-manifold $S^3(K,a,\{q_i/r_i\}_{i=1}^n)$, we obtain the Seifert
fibered space whose base has genus zero, $n$ singular fibers, and
Seifert invariants $(a,\{q_i/r_i\}_{i=1}^n)$, with the standard
conventions. More generally, if we
start with the Borromean knot $B_g\subset \#^{2g}(S^1\times S^2)$
(this is the knot obtained by taking zero surgery on two of the
components of the Borromean rings to obtain a knot $B_1$ in $\#^2
(S^1\times S^2)$, and then taking the connected sum of $g$ copies of
$B_1$) , the resulting three-manifold $\#^{2g}(S^1\times
S^2)(B_g,a,\{q_i/r_i\}_{i=1}^n)$ is the Seifert space over a genus $g$
base orbifold with Seifert invariants $(a,\{q_i/r_i\}_{i=1}^n)$.

Of course, 
$$H_1(Y';\Z)\cong 
\frac{H_1(Y-K;\Z)\oplus \Z m_1 \oplus ... \oplus \Z m_n}
{\left(\begin{array}{ll}
a\cm  m_0 + \sum_{i=1}^n m_i  = 0, & \\
 r_i \cm m_0 - q_i \cm m_i = 0, & i=1,...,m
\end{array}\right)},$$
where here $m_0\in H_1(Y-K;\Z)$ denotes the homology class
of the meridian of $K$.
We can change basis, letting $g_i = a_i\cm  m_0 + b_i\cm  m_i$,
where here $a_i\cm q_i + r_i \cm b_i = 1$, to obtain a presentation
$$H_1(Y';\Z)\cong 
\frac{H_1(Y-K;\Z)\oplus \Z g_1 \oplus ... \oplus \Z g_n}
{\left(\begin{array}{ll}
a\cm  m_0 + \sum_{i=1}^n r_i \cm g_i  = 0, & \\
m_0 - q_i \cm g_i = 0, & i=1,...,m
\end{array}\right)}.$$

Consider the map
$$ \beta\colon H_1(Y'-K';\Z) \longrightarrow H_1(Y-K;\Z) $$ 
defined by
$$\beta(\xi_0 + \xi_1 \cm g_1 + ... + \xi_n \cm g_n)
= \xi_0 + \left(\sum_{i=1}^n \lfloor \frac{\xi_i}{q_i}\rfloor\right).$$

Letting $N=a \cm m_0 + \sum_{i=1}^n r_i\cm g_i$ be the
push-off of $K$, we have that
$$H_1(Y(K,a,\{q_i/r_i\}_{i=1}^n);\Z)\cong H_1(Y'-K')/\Z\cm N.$$
Given $E\in H_1(Y-K)$, 
define 
\begin{eqnarray*}
\BigAp_{[E]}=\bigoplus_{s\in\Z} (s,\Ap_{\beta(E+s \cm N)}(Y,K)) &{\text{and}}&
\BigBp_{[E]}=\bigoplus_{s\in\Z} (s,\Bp_{\beta(E+s\cm N)}(Y)),
\end{eqnarray*}
and
$$\Dp_{[E]}\colon \BigAp_{[E]}\longrightarrow \BigBp_{[E]}$$
with
$$\Dp_{[E]} \{a_s\}_{s\in \Z} = \{b_s\}_{s\in \Z}$$
where $a_s \in \Ap_{\deg(E+s\cm N)}$ and $b_s\in \Bp_{\deg(E+s\cm N)}$
by
$$b_s = \horp_{\beta(E+(s-1)\cm N)}(a_{s-1}) + 
\vertp_{\beta(E+s\cm N)}(a_{s}).
$$

\begin{prop}
\label{prop:GenRatSurg}
Let $K\subset Y$ be a null-homologous knot, and fix an $E\in
\SpinC(Y'-K')$. Suppose moreover that $b_1(Y(K,a,\{q_i/r_i\}))=b_1(Y)$.
There is a map $f\colon H_1(Y'-K')\longrightarrow
\SpinC(Y(K,a,\{q_i/r_i])$ with the property that for each $E\in Y'-K'$
for which $c_1(f(E))$ is a torsion class, we have that
$\HFp(Y(K,a,\{q_i/r_i\}_{i=1}^n),f(E))$ is identified with the homology
of the mapping cone
$\Dp_{[E]}\colon
\BigAp_{[E]} \longrightarrow
\BigBp_{[E]}$.
\end{prop}

\begin{proof}
Like Theorem~\ref{thm:RationalSurgeries},
this follows from  a direct application of Theorems~\ref{thm:SurgeryFormula},
~\ref{thm:Kunneth}, and the calculation of the invariant for
$O_{q/r}\subset L(q,r)$ of Lemma~\ref{lemma:LensSpaceUKnots}.
\end{proof}

Theorem~\ref{thm:SeifertFiberedSpaces} follows quickly from
Proposition~\ref{prop:GenRatSurg}, together with the calculation of
the knot Floer homology of the Borromean knot (c.f.
Section~\ref{Knots:sec:FiberedExamples} of~\cite{Knots}),
which we summarize here:

\begin{lemma}
        \label{lemma:BorrKnot}
        Let $\BK\subset \#^{2g}(S^2\times S^1)$ be the Borromean knot.
        The $\SpinC$ structure over $\#^{2g}(S^2\times S^1)$ with
        trivial first Chern class is the only one whose induced knot
        filtration consists of non-zero groups.  For that chain
        complex $C$, we have a splitting
        \begin{equation}
        \label{eq:BorrKnot}
        C \cong \bigoplus_{t} \Wedge^{t+g}
        H_1(\Sigma;\Z)\otimes_{\Z} \Z[U,U^{-1}],
        \end{equation} where this
        splitting corresponds to the various summands $C\{i-j=t\}$.
        In particular, there are identifications \begin{eqnarray*}
        \Ap_s(\#^{2g}(S^2\times S^1),\BK) &\cong& \bigoplus_t \Wedge^{t+g}
        H_1(\Sigma;\Z) \otimes_\Z \InjMod{} \\ C\{i\geq 0\}(\#^{2g}(S^2\times S^1),\BK)
        &\cong& \bigoplus_t \Wedge^{t+g} H_1(\Sigma;\Z) \otimes_\Z
        \InjMod, \end{eqnarray*} Moreover, the following squares
        commute: 
        $$\begin{CD} \Wedge^{t+g} H_1(\Sigma;\Z) \otimes_\Z \InjMod{}
        @>{U^{\max(0,t-s)}}>> \Wedge^{t+g} H_1(\Sigma;\Z) \otimes_\Z \InjMod{} \\
        @VVV @VVV \\
        \Ap_s @>{\vertp}>> \Bp
        \end{CD}
        $$
        $$\begin{CD} \Wedge^{t+g} H_1(\Sigma;\Z) \otimes_\Z \InjMod{}
        @>{U^{\max(0,s-t)}}>> \Wedge^{t+g} H_1(\Sigma;\Z) \otimes_\Z \InjMod{} 
        @>{\cong}>> \Wedge^{-t+g} H_1(\Sigma;\Z) \otimes_\Z \InjMod{} \\
        @VVV && @VVV \\
        \Ap_s @>{\horp}>>\Bp @>{=}>> \Bp,
        \end{CD}
        $$
        where here all the vertical maps are induced by inclusions in the identification
        of Equation~\eqref{eq:BorrKnot}.
\end{lemma}

\begin{proof}
Most of the above statements are a direct result of
Proposition~\ref{Knots:prop:KnotHomology} of~\cite{Knots},
which, together with the K\"unneth principle, gives that $$C\{i,j\} =
U^{-i}\otimes \Wedge^{g-i+j} H^1(\Sigma;\Z),$$ with no
differentials. The second square involves the identification between
$C\{i\geq 0\}$ and $C\{j\geq 0\}$, which is sends the subset
$C\{i-j=t\}$ to $C\{i-j=-t\}$, as explained in
Proposition~\ref{IntSurg:prop:HomEq} of~\cite{IntSurg}.
\end{proof}

\begin{lemma}
\label{lemma:Bamboo}
Given a map $\delta \colon \Z \longrightarrow \Z$, consider the chain
map
$$D_\delta \colon \bigoplus_{s} A_s \longrightarrow \bigoplus_{s}
B_s,$$
where all $A_s \cong \InjMod{} \cong B_s$, defined by
$$D_\delta(\{a_s\}_{s\in\Z}) = \{b_s\}_{s\in\Z}$$
where
$$b_s = U^{\max(-\delta(s),0)} a_s + U^{\max(\delta(s+1),0)}
a_{s+1}.$$
Letting
$$h(s)=\sum_{i=0}^s \delta(i),$$
we have that the homology of the
mapping cone of $D_{\delta}$ is $\FormHFp(h)$, provided that
\begin{equation}
\label{eq:BehaviourAtInfinity}
\lim_{s\goesto \pm \infty} h(s)=\pm \infty.
\end{equation}
\end{lemma}
        
\begin{proof}
        The map $D_\delta$ is surjective. It remains to identify its kernel.
        Given $(i,j)\in W_n(h)$, consider the element $\{a_s\}_{s\in\Z}$ 
        defined by the property that 
        $$a_s = \left\{\begin{array}{ll}
                        (-1)^{k} U^{h(k)-n} & {\text{if $i<k<j$}}  \\
                        0 & {\text{otherwise.}}
                        \end{array}\right.
        $$
        By linearity, we can extend this to a homomorphism
        $W_n(h) \longrightarrow \Ker D_\delta$.
        It is straightforward to verify that this extends to an isomorphism
        $\FormHFp(h) \longrightarrow \Ker D_\delta$.
\end{proof}

\vskip.2cm
\noindent{\bf{Proof of Theorem~\ref{thm:SeifertFiberedSpaces}.}}
According to the adjunction inequality
(c.f. Theorem~\ref{HolDiskTwo:thm:Adjunction}
of~\cite{HolDiskTwo}), since the Thurston norm of $Y$ is trivial, it follows 
that $\HFp(Y,\spinc)$ is trivial for all $\SpinC$ structures with
non-torsion $c_1(\spinc)$.  According to the combination of
Proposition~\ref{prop:GenRatSurg} and Lemma~\ref{lemma:BorrKnot},
given $E$, the Floer homology $\HFp(Y,[E])$ splits as
$$\HFp(Y,[E])\cong \bigoplus_{t\in\Z} \Wedge^{g+t} H^1(\Sigma;\Z)
\otimes_\Z X(t),$$ where here $X(t)$ is the homology of a chain
complex satisfying the hypotheses of Lemma~\ref{lemma:Bamboo}, for the
function $\delta_t\colon \Z \longrightarrow \Z$ as in the statement of
the theorem.  The theorem now follows from the calculation in
Lemma~\ref{lemma:Bamboo}.  Note that
Equation~\eqref{eq:BehaviourAtInfinity} holds, since the orbifold has
positive degree.
\qed
\vskip.2cm

\commentable{ 
\bibliographystyle{plain} 
\bibliography{biblio} }

\begin{thebibliography}{10}

\bibitem{Berge}
J.~O. Berge.
\newblock Some knots with surgeries giving lens spaces.
\newblock Unpublished manuscript.

\bibitem{BleilerHodgsonWeeks}
S.~A. Bleiler, C.~D. Hodgson, and J.~R.~Jeffrey Weeks.
\newblock Cosmetic surgery on knots.
\newblock In {\em Proceedings of the Kirbyfest (Berkeley, CA, 1998)}, volume~2
  of {\em Geom. Topol. Monogr.}, pages 23--34 (electronic). Geom. Topol. Publ.,
  Coventry, 1999.

\bibitem{Snap}
David Coulson, Oliver~A. Goodman, Craig~D. Hodgson, and Walter~D. Neumann.
\newblock Computing arithmetic invariants of 3-manifolds.
\newblock {\em Experiment. Math.}, 9(1):127--152, 2000.

\bibitem{Eliashberg}
Y.~Eliashberg.
\newblock A few remarks about symplectic filling.
\newblock {\em Geom. Topol.}, 8:277--293, 2004.

\bibitem{Etnyre}
J.~B. Etnyre.
\newblock On symplectic fillings.
\newblock {\em Algebr. Geom. Topol.}, 4:73--80, 2004.

\bibitem{Froyshov}
K.~A. Fr{\o}yshov.
\newblock The {S}eiberg-{W}itten equations and four-manifolds with boundary.
\newblock {\em Math. Res. Lett}, 3:373--390, 1996.

\bibitem{FurutaSteer}
M.~Furuta and B.~Steer.
\newblock Seifert fibered homology 3-spheres and the {Y}ang-{M}ills equations
  on {R}iemann surfaces with marked points.
\newblock {\em Advances in Mathematics}, 96:38--102, 1993.

\bibitem{GordonConjecture}
C.~McA. Gordon.
\newblock {\em Some aspects of classical knot theory}, pages pp. 1--60.
\newblock Number 685 in Lecture Notes in Math. Springer-Verlag, 1978.

\bibitem{KMOSz}
P.~B. Kronheimer, T.~S. Mrowka, P.~S. Ozsv{\'a}th, and Z.~Szab{\'o}.
\newblock Monopoles and lens space surgeries.
\newblock math.GT/0310164.

\bibitem{Lackenby}
M.~Lackenby.
\newblock Dehn surgery on knots in {$3$}-manifolds.
\newblock {\em J. Amer. Math. Soc.}, 10(4):835--864, 1997.

\bibitem{Mathieu}
Y.~Mathieu.
\newblock Closed {$3$}-manifolds unchanged by {D}ehn surgery.
\newblock {\em J. Knot Theory Ramifications}, 1(3):279--296, 1992.

\bibitem{Nemethi}
A.~N{\'e}methi.
\newblock On the {Ozsv{\'a}th-Szab{\'o}} invariant of negative definite plumbed
  $3$-manifolds.
\newblock math.GT/0310083.

\bibitem{Neumann}
W.~Neumann.
\newblock Private communication.

\bibitem{HolDiskFour}
P.~S. Ozsv{\'a}th and Z.~Szab{\'o}.
\newblock Holomorphic triangles and invariants for smooth four-manifolds.
\newblock math.SG/0110169.

\bibitem{UnknotOne}
P.~S. Ozsv{\'a}th and Z.~Szab{\'o}.
\newblock Knots with unknotting number one and {H}eegaard {F}loer homology.
\newblock math.GT/0401426.

\bibitem{NoteLens}
P.~S. Ozsv{\'a}th and Z.~Szab{\'o}.
\newblock On knot {F}loer homology and lens space surgeries.
\newblock math.GT/0303017, to appear in {\em Topology}.

\bibitem{AbsGraded}
P.~S. Ozsv{\'a}th and Z.~Szab{\'o}.
\newblock Absolutely graded {F}loer homologies and intersection forms for
  four-manifolds with boundary.
\newblock {\em Advances in Mathematics}, 173(2):179--261, 2003.

\bibitem{4BallGenus}
P.~S. Ozsv{\'a}th and Z.~Szab{\'o}.
\newblock Knot {F}loer homology and the four-ball genus.
\newblock {\em Geom. Topol.}, 7:615--639, 2003.

\bibitem{SomePlumbs}
P.~S. Ozsv{\'a}th and Z.~Szab{\'o}.
\newblock On the {F}loer homology of plumbed three-manifolds.
\newblock {\em Geometry and Topology}, 7:185--224, 2003.

\bibitem{BrDCov}
P.~S. Ozsv{\'a}th and Z.~Szab{\'o}.
\newblock On the {H}eegaard {F}loer homology of branched double-covers.
\newblock math.GT/0309170, to appear in {\em Adv. Math.}, 2003.

\bibitem{GenusBounds}
P.~S. Ozsv{\'a}th and Z.~Szab{\'o}.
\newblock Holomorphic disks and genus bounds.
\newblock {\em Geom. Topol.}, 8:311--334 (electronic), 2004.

\bibitem{Knots}
P.~S. Ozsv{\'a}th and Z.~Szab{\'o}.
\newblock Holomorphic disks and knot invariants.
\newblock {\em Adv. Math.}, 186(1):58--116, 2004.

\bibitem{HolDiskTwo}
P.~S. Ozsv{\'a}th and Z.~Szab{\'o}.
\newblock Holomorphic disks and three-manifold invariants: properties and
  applications.
\newblock {\em Ann. of Math. (2)}, 159(3):1159--1245, 2004.

\bibitem{HolDisk}
P.~S. Ozsv{\'a}th and Z.~Szab{\'o}.
\newblock Holomorphic disks and topological invariants for closed
  three-manifolds.
\newblock {\em Ann. of Math. (2)}, 159(3):1027--1158, 2004.

\bibitem{IntSurg}
P.~S. Ozsv{\'a}th and Z.~Szab{\'o}.
\newblock Knot {F}loer homology and integral surgeries.
\newblock math.GT/0410300, 2004.

\bibitem{calcKT}
P.~S. Ozsv{\'a}th and Z.~Szab{\'o}.
\newblock Knot {F}loer homology, genus bounds, and mutation.
\newblock {\em Topology Appl.}, 141(1-3):59--85, 2004.

\bibitem{RasmussenThesis}
J.~A. Rasmussen.
\newblock {\em Floer homology and knot complements}.
\newblock PhD thesis, Harvard University, 2003.

\bibitem{Scott}
P.~Scott.
\newblock The geometries of 3-manifolds.
\newblock {\em Bull. London Math. Soc.}, 15:401--487, 1983.

\bibitem{Seifert}
H.~Seifert.
\newblock Topologie dreidimensionaler gefaserter {R\"a}ume.
\newblock {\em Acta Math.}, 60:147--238, 1932.

\bibitem{Turaev}
V.~Turaev.
\newblock Torsion invariants of {S}pin{$^c$}-structures on $3$-manifolds.
\newblock {\em Math. Research Letters}, 4:679--695, 1997.

\bibitem{SnapPea}
J.~Weeks.
\newblock Snappea, a program for creating and studying hyperbolic volume.
\newblock http://www.geometrygames.org/SnapPea/index.html.

\end{thebibliography}
\end{document}